\documentclass[review,17pt,reqno]{Elsevier}
\usepackage{amssymb,amsmath}
\usepackage{cases}
\usepackage{exscale}
\usepackage{relsize}
\usepackage{subfigure}
\usepackage{natbib}
\usepackage[colorlinks=true,linkcolor=blue,citecolor=blue,urlcolor=blue,]{hyperref}
\usepackage{caption}
\usepackage{booktabs}
\usepackage{makecell}
\usepackage{graphicx}
\usepackage{pythonhighlight}
\usepackage{array, threeparttable}
\usepackage{algorithmicx}
\usepackage[ruled]{algorithm2e}    
\usepackage{algpseudocode}
\usepackage{mathrsfs} 
\usepackage{arydshln}
\usepackage{multirow}
\usepackage{graphicx} 
\usepackage{mathtools}
\usepackage{epsfig}                          
\usepackage{epstopdf}                        
\usepackage{appendix}

\begin{document}
	\newtheorem{mytheorem}{Theorem}[section]
	\newtheorem{myproposition}{Proposition}[section]
	\newtheorem{mylemma}{Lemma}[section]
	\newtheorem{mycorollary}[mytheorem]{Corollary}
	\newtheorem{myprop}[myproposition]{Proposition}
	\newtheorem{myremark}{Remark}[section]
	\renewcommand{\baselinestretch}{1.1}
	\newcommand{\figref}[1]{Fig.~\ref{#1}}
	\makeatletter
	\newcommand\figcaption{\def\@captype{figure}\caption}
	\newcommand\tabcaption{\def\@captype{table}\caption}
	\makeatother
\begin{frontmatter}
\title{Efficient primal--dual splitting methods for a Poisson-constrained JKO scheme for Poisson-Nernst-Planck models}
		
\author[label1]{Wei Wu}
\ead{wuwei837037@163.com}

\author[label1]{Jin Zeng}
\ead{202511110511@std.uestc.edu.cn}

\author[label2]{Zhen Zhang}
\ead{zhangz@sustech.edu.cn}

\author[label1]{Chaozhen Wei\corref{cor1}}
\ead{cwei4@uestc.edu.cn}\cortext[cor1]{Corresponding author.}

\address[label1]{School of Mathematical Sciences, University of Electronic Science and Technology of China, Chengdu, Sichuan 611731, China}

\address[label2]{Department of Mathematics, National Center for Applied Mathematics (Shenzhen), Southern University of Science and Technology (SUSTech), Shenzhen 518055, China}			
\begin{abstract}
The Poisson--Nernst--Planck (PNP) equations strongly couple ionic transport and electrostatic interactions through the Poisson equation, posing substantial numerical challenges under small permittivity and complex potential boundary conditions. Underlying these equations is a natural Wasserstein gradient-flow structure, in which the Poisson equation serves as a local realization of the nonlocal electrostatic interaction energy. Exploiting this structure, we formulate each time step as a constrained convex minimization problem where the ionic continuity equations and the Poisson equation are incorporated as linear constraints, allowing the concentrations, fluxes, and electrostatic potential to be updated simultaneously. The variational structure of the scheme intrinsically guarantees the dissipation of the original free energy, mass conservation, and nonnegativity of ionic concentrations under general electrostatic boundary conditions. Moreover, the framework is structurally modular: extending from classical to modified PNP models with steric interactions and concentration-gradient corrections requires only modifying the energy functional, while all structure-preserving properties are automatically retained. To efficiently solve the resulting large-scale constrained problems, we develop preconditioned and  transformed primal--dual algorithms equipped with tailored fast dual solvers, namely DCT-based direct and Schur-complement iterative methods, that exploit the coupled block structure of the PDE constraints. Numerical experiments on classical and modified PNP systems demonstrate the accuracy and structure-preserving properties of the scheme, and show that the proposed algorithms converge reliably in strongly coupled small-permittivity regimes without significant growth in computational cost.
\end{abstract}
		
\end{frontmatter}

\section{Introduction}

In this paper, we study numerical approximations for Poisson--Nernst--Planck (PNP) models. The models have been widely used to describe charge transport phenomena in electrochemistry \cite{Bazant2004Diffuse,Gillespie2002Coupling,Weber2004Modeling,Latz2011Thermodynamic} and biology \cite{Nonner1998IonPNP,Noskov2004Control,Eisenberg1998Ionic,Im2002Ion}. The PNP model consists of the Nernst--Planck equations for ionic transport and the Poisson equation for electrostatic interactions. We consider Nernst--Planck equations on a bounded domain $\Omega\subset\mathbb{R}^{d}$ ($d=2,3$):
\begin{align}\label{eq:NPequation}
	\dfrac{\partial c_{i}}{\partial t}
	=\nabla\cdot\left(D_{i}c_{i}\nabla\left(\log c_{i}
	+\dfrac{z_{i}e(\phi+\phi_{e})}{k_{B}T}\right)\right),
	\quad i=1,\ldots,N,
\end{align}
where $c_{i}$ denotes the concentration of the $i$-th ionic species, $D_{i}$ is the diffusion coefficient, $z_{i}$ is the valence, $e$ is the elementary charge, $\phi$ is the electrostatic potential, $\phi_{e}$ is the external electric potential, $k_{B}$ is the Boltzmann constant, and $T$ is the absolute temperature. The electrostatic potential $\phi$ in \eqref{eq:NPequation} is determined by the following Poisson equation with nonhomogeneous mixed boundary conditions:
\begin{align}\label{eq:Poisson}
	\left\{
	\begin{aligned}
		&-\nabla\cdot(\epsilon\nabla\phi)=\psi\quad&&(\boldsymbol{x},t)\in\Omega\times\left[0,T\right],\\
		&\alpha\phi+\beta\epsilon\dfrac{\partial\phi}{\partial\mathbf{n}}=\phi^{bc}&&(\boldsymbol{x},t)\in\partial\Omega\times\left[0,T\right],
	\end{aligned}
	\right.
\end{align}
where $\psi=\sum_{i=1}^{N}z_{i}ec_{i}+\psi^{0}$ denotes the total charge density, consisting of the mobile ionic charge density and the fixed charge density $\psi^{0}$. The coefficient $\epsilon$ represents the dielectric permittivity, $\mathbf{n}$ is the unit outward normal vector on $\partial\Omega$, and $\phi^{bc}$ denotes the prescribed boundary data in the generalized boundary condition with coefficients $\alpha$ and $\beta$. Under suitable boundary conditions, the total free energy of classical PNP model is defined as the sum of the entropic energy and electrostatic energy:
\begin{align}\label{eq:totalenebc}
	\mathcal{E}_{pnp}\left\{c_{i},\phi\right\}=\int_{\Omega}k_{B}T\sum_{i=1}^{N}c_{i}\log c_{i}\mathrm{d}\boldsymbol{x}+\mathcal{P}\left\{c_{i},\phi\right\},
\end{align}
where $\mathcal{P}(c_i,\phi)$ denotes the electrostatic energy, consisting of the bulk electrostatic energy and the boundary correction energy $\mathcal{E}_{bc}(\phi)$:
\begin{align}
    \mathcal{P}\left\{c_{i},\phi\right\}=\dfrac{1}{2}\int_{\Omega}\left(\sum_{i=1}^{N}z_{i}ec_{i}+\psi^{0}\right)(\phi+\phi_{e})\mathrm{d}\boldsymbol{x}+\mathcal{E}_{bc}\left\{\phi\right\}.
\end{align}
In addition to the above energies, modified PNP models enrich the classical mean-field free energy by adding steric or excluded-volume interaction energies, and concentration-gradient terms. These additions represent finite ion sizes, short-range ionic correlations, solvent occupancy, and strong spatial variations of the ionic concentrations that are absent from the classical model \cite{Jiang2014DCTPNP,Lu2011DCTPNP,Gavish2020LJPPNP,Qian2021PositivePNPCH}. Both the classical and modified PNP systems can be viewed as Wasserstein gradient flows for the corresponding free energy \cite{Liu2023DynamicMassPNP}, in which the Poisson equation arises naturally as a local realization of the nonlocal electrostatic interaction energy.

Accurate and robust numerical simulation of PNP systems is challenging due to the nonlinear coupling between the Nernst--Planck equations and the electrostatic potential determined by the Poisson equation. Reliable numerical methods must control this coupling while keeping every ionic concentration nonnegative, conserving the mass of each species, and preserving the dissipation of the original free energy at the discrete level. An early second-order finite-difference method conserved ionic mass exactly and ensured positivity under step-size restrictions \cite{Flavell2014ConservativePNPFDM}; a companion discretization reproduced a discrete energy law \cite{Flavell2017EnergyFDMPNP}. Hu and Huang combined a Scharfetter--Gummel reformulation with semi-implicit stepping to prove mass conservation, unconditional positivity, and energy dissipation \cite{Hu2020SSIPNP}. Shen and Xu instead used the logarithmic Wasserstein gradient-flow form with semi-implicit discretization to construct first- and second-order positivity-preserving and mass conservative schemes, while unconditional energy dissipation was only proved for the first-order scheme \cite{Shen2021PNP}. High-order linear SAV schemes preserve positivity, mass, and unconditional modified energy stability \cite{Huang2021SAVPNP}. A third-order direct DG method based on a nonlogarithmic Landau transformation and a scaling limiter preserves positivity, mass conservation, and steady states \cite{Liu2022PositivityDGPNP}. Tong and Cai instead combined a second-order Crank--Nicolson discretization with an $L^2$ projection to enforce positivity and mass conservation \cite{Tong2024PositivityPNP}. For modified PNP models, nonlinear or singular steric terms and concentration-gradient energies add stiffness and, in the latter case, fourth-order operators. Ding {\it et al.} used harmonic-mean approximations of Slotboom variables to construct mass-conservative and positivity-preserving schemes for steric PNP, with unconditional positivity under backward Euler \cite{Ding2019PNPSteric}. For PNP--Cahn--Hilliard systems, Qian et al. proposed a nonlinear semi-implicit scheme formulated as convex minimization; the singular logarithmic entropy ensures positivity, while the scheme conserves mass and dissipates the discrete free energy. Their subsequent work established optimal-rate convergence \cite{Qian2021PositivePNPCH,Qian2023convergencePNPCH}. Ding and Zhou later combined second-order time discretization with a multislope finite-volume reconstruction of positive mobilities on unstructured meshes, preserving positivity, mass, dissipation of the original energy, and steady states \cite{Ding2024SecondorderPNP}. Despite these advances, several PNP-specific computational bottlenecks remain: conservation depends on flux and boundary discretizations, positivity requires control of the logarithmic entropy and mobility or an additional limiter or projection, and dissipation of the original energy requires a compatible coupled update of the Nernst--Planck and Poisson equations. Although some schemes attain all three structures in specific settings, doing so simultaneously with general electrostatic boundary conditions, multi-species coupling, and affordable nonlinear solves remains difficult; concentration-gradient corrections intensify this difficulty.

In this work, we develop a unified Poisson-constrained JKO scheme for classical and modified PNP models. In contrast to the traditional PDE-based discretization approaches that enforce individual structures through term-specific stabilization or post-processing, the proposed method leverages the underlying Wasserstein gradient flow structure of PNP models \cite{Liu2023DynamicMassPNP} to ensure structure-preserving properties, and provides numerical realization of the Jordan-Kinderlehrer-Otto ({\it JKO}) minimizing-movement scheme \cite{Jordan1998Variational}. At each time step, the Nernst-Planck equation of ionic concentrations are converted to an equivalent variational formulation that reduces to a series of convex minimization problems, where the ionic concentrations and electrostatic potential are treated as independent variables, and their couplings by the Poisson equation with boundary conditions are imposed as linear constraints.
Building upon our previous structure-preserving primal--dual JKO method for Wasserstein gradient flows \cite{Carrillo2022PrimalDual,deng2025PDFB,Wu2026PDPFSMCL}, we develop preconditioned primal--dual (PrePD) method \cite{Carrillo2024StructurePD} and variable-preconditioned transformed primal--dual (VPTPD) method \cite{Zeng2026variable} for the resulting Poisson-constrained minimization problems. Tailored fast dual solvers exploit the coupled block structure of the PDE constraints under various boundary conditions, exhibiting robust and efficient performance in a suite of numerical experiments.

This formulation, in particular, has several essential differences compared to previous works \cite{Hu2020SSIPNP,Shen2021PNP,Qian2021PositivePNPCH}, which circumvented the direct computation of the Wasserstein distance by approximating the Wasserstein gradient flow with a weighted $H^{-1}$ gradient flow, together with a semi-implicit discretizations of the mobility and chemical potential. In contrast, the JKO approach directly works with the Wasserstein metric and can be viewed as an implicit Euler scheme in Wasserstein space, thereby naturally inheriting the associated energy-dissipation structure and global bounds of solutions with the mobility treated implicitly. The advantages are three-fold. 
First, the JKO scheme discretizes the original gradient-flow structure itself, so that dissipation of the original energy and nonnegativity of ionic concentrations are intrinsic consequences of the variational formulation. In particular, original energy dissipation follows from the discrete gradient-flow structure, without introducing energy modifications that may alter the long-time asymptotic behavior \cite{Xu2019CMAME}; positivity is enforced intrinsically by the admissible set of the Wasserstein transport action, rather than through variable transformations, post-processing projections or relying on the effective domain of the logarithmic potential to impose the physical bounds. 
Second, the primal--dual splitting algorithms developed for the resulting constrained optimization problems, together with the tailored fast dual solvers that exploit the coupled block structure of the PDE constraints, converge reliably in strongly coupled small-permittivity regimes without significant growth in computational cost. 
Third, the JKO framework is structurally modular: additional physical effects can be incorporated through either additional free-energy functional or coupling constraints without compromising the underlying structure-preserving properties. This modularity makes the formulation readily extensible to multi-species, modified, and multiphysics coupled systems without redesigning the underlying transport structure. 

The rest of the paper is organized as follows. In Sec.~\ref{sec:2}, we introduce the classical PNP model and modified PNP models, and discuss their fundamental properties. In Sec.~\ref{sec:3}, we first review the JKO scheme based on the Wasserstein gradient-flow structure, derive a variational formulation of PNP models, and prove that the proposed fully discrete variational scheme preserves the essential structural properties. In Sec.~\ref{sec:4} and Sec.~\ref{sec:5}, we elaborate the primal--dual splitting method for the resulting convex optimization problems, and the fast dual solvers for the coupled dual subproblem with various boundary conditions. In Sec.~\ref{sec:6}, we present a series of numerical experiments to demonstrate the effectiveness and efficiency of the proposed methods. Finally, we conclude the paper with a brief summary and outlook.

\section{Poisson-Nernst-Planck equations and its extensions}\label{sec:2}

\subsection{Poisson-Nernst-Planck equations}\label{sec:PNP}
We first recall the classical PNP model on a bounded connected domain $\Omega$. The Nernst-Planck equations~\eqref{eq:NPequation} can be rewritten in the form of continuity equations with concentration-dependent mobilities:
\begin{align}\label{eq:pnpgradient}
	\dfrac{\partial c_{i}}{\partial t}=\nabla\cdot\Big(\dfrac{D_{i}}{k_{B}T}c_{i}\nabla\mu_{i}\Big)\quad i=1,\cdots,N,
\end{align}
which can be regarded as coupled Wasserstein gradient flows for the free energy \eqref{eq:totalenebc} in the metric space of probability measures on $\Omega$~\cite{Jordan1998Variational,Liu2023DynamicMassPNP}. In this interpretation, the mobility is given by $\dfrac{D_{i}}{k_{B}T}c_{i}$, which depends on the ion concentration, and the driving force is given by the chemical potential $\mu_{i}$, which is the variational derivative of the free energy \eqref{eq:totalenebc} with respect to $c_{i}$:
\begin{align}\label{eq:chemical}
	\mu_{i}=\dfrac{\delta\mathcal{E}_{pnp}}{\delta c_{i}}=k_{B}T\log c_{i}+z_{i}e(\phi+\phi_{e}).
\end{align}
The coupling between the ionic concentrations and the electrostatic potential through the Poisson equation~\eqref{eq:Poisson} adds an additional layer of complexity to this gradient flow structure.


To close the system, we consider consistent initial values for the ionic concentrations and the electrostatic potential, and appropriate boundary conditions. Typically, one imposes no-flux boundary conditions for the ionic concentrations in the Nernst-Planck equations
\begin{align}\label{eq:cinoflux}
	D_{i}c_{i}\nabla\dfrac{\delta\mathcal{E}_{pnp}}{\delta c_{i}}\cdot\mathbf{n}=0.
\end{align}
For regular domains, periodic boundary conditions may also be used.
The choice of boundary conditions for the electrostatic potential depends on the physical setting of specific problems. In this work, we consider the following general form for the boundary condition of the Poisson equation:
\begin{align}\label{eq:phibc}
	\alpha\phi+\beta\epsilon\dfrac{\partial\phi}{\partial\mathbf{n}}=\phi^{bc},\quad(\boldsymbol{x},t)\in\partial\Omega\times\left[0,T\right],
\end{align}
where $\partial\Omega=\Gamma_{D}\cup\Gamma_{N}\cup\Gamma_{R}$ denotes the decomposition of the boundary into the Dirichlet, Neumann, and Robin parts, respectively. We normalize the coefficients in \eqref{eq:phibc} by taking $\alpha=1$ and $\beta=0$ on $\Gamma_{D}$, $\alpha=0$ and $\beta=1$ on $\Gamma_{N}$, and $\alpha=\alpha_{R}>0$ and $\beta=\beta_{R}>0$ on $\Gamma_{R}$. The wellposedness of the Poisson equation requires the compatibility condition
\begin{align}
	\int_{\Omega}\left(\psi^{0}+\sum_{i=1}^{N}z_{i}ec_{i}\right)
	\mathrm{d}\boldsymbol{x}
	+\int_{\partial\Omega}\frac{\phi^{bc}}{\beta}\,\mathrm{d}s=0.
\end{align}
In this case, the electrostatic potential is determined only up to an additive constant, and an additional gauge condition is needed for uniqueness.

In summary, we obtain the following classical PNP model
\begin{align}\label{eq:classicalPNP}
\left\{
\begin{aligned}
		&\dfrac{\partial c_{i}}{\partial t}=\nabla\cdot\Bigg(D_{i}c_{i}\nabla\Big(\log c_{i}+\dfrac{z_{i}e(\phi+\phi_{e})}{k_{B}T}\Big)\Bigg)&&(\boldsymbol{x},t)\in\Omega\times\left[0,T\right]\quad i=1,\cdots,N,\\
		&-\nabla\cdot(\epsilon\nabla\phi)=\sum_{i=1}^{N}z_{i}ec_{i}+\psi^{0}&&(\boldsymbol{x},t)\in\Omega\times\left[0,T\right],
	\end{aligned}
	\right.
\end{align}
subject to the following initial and boundary conditions
\begin{align}\label{eq:classicalbc}
\left\{
    \begin{aligned}
		&c_{i}(\mathbf{x},0)=c^{0}_{i}(\mathbf{x})&&(\boldsymbol{x},t)\in\Omega\times\left\{0\right\}\quad i=1,\cdots,N,\\
		&D_{i}\Big(\nabla c_{i}+\dfrac{ z_{i}ec_{i}}{k_{B}T}\nabla(\phi+\phi_{e})\Big)\cdot\mathbf{n}=0&&(\boldsymbol{x},t)\in\partial\Omega\times\left[0,T\right],\quad i=1,\cdots,N,\\
		&\alpha\phi+\beta\epsilon\dfrac{\partial\phi}{\partial\mathbf{n}}=\phi^{bc}&&(\boldsymbol{x},t)\in\partial\Omega\times\left[0,T\right].
	\end{aligned}
	\right.
\end{align}

We next summarize several energy representations and basic structural properties of \eqref{eq:classicalPNP}. For convenience of discussion, we set the elementary charge $e=1$, the external potential $\phi_{e}=0$, and $k_{B}T=1$.
\begin{myremark}\label{rmk:Energy}
    For mixed boundary conditions on $\Gamma_D$, $\Gamma_N$, and $\Gamma_R$ with $\phi^{bc}\neq 0$, the electrostatic part of the free energy must include appropriate boundary correction terms to compensate for the boundary contributions arising from integration by parts so that the energy is still dissipating along the solution \cite{Liu2023DynamicMassPNP}. Assume that the boundary data are independent of time, and denote by $\phi_{D}^{bc}$, $\phi_{N}^{bc}$, and $\phi_{R}^{bc}$ the prescribed data on the Dirichlet, Neumann, and Robin boundaries, respectively. Then the free energy of the classical PNP model can be written as
	\begin{align}\label{eq:totalEnePNP}
		\mathcal{E}_{pnp}\left\{c_{i},\phi\right\}=\int_{\Omega}\Big(\sum_{i=1}^{N}c_{i}\log c_{i}+\dfrac{1}{2}\psi\phi\Big)\mathrm{d}\boldsymbol{x}-\dfrac{1}{2}\int_{\Gamma_{D}}\epsilon\phi^{bc}_{D}\partial_{\mathbf{n}}\phi\mathrm{d}s+\dfrac{1}{2}\int_{\Gamma_{N}}\phi^{bc}_{N}\phi\mathrm{d}s+\dfrac{1}{2\beta_{R}}\int_{\Gamma_{R}}\phi^{bc}_{R}\phi\mathrm{d}s.
	\end{align}
    Here $\psi=\sum_{i=1}^{N}z_{i}c_i+\psi^{0}$ under the normalization $e=1$. Substituting the Poisson equation into \eqref{eq:totalEnePNP} and applying Green's formula give the equivalent form
	\begin{align}
		\mathcal{E}_{pnp}\left\{c_{i},\phi\right\}=\int_{\Omega}\Big(\sum_{i=1}^{N}c_{i}\log c_{i}+\dfrac{1}{2}\epsilon|\nabla\phi|^{2}\Big)\mathrm{d}\boldsymbol{x}-\int_{\Gamma_{D}}\epsilon\phi^{bc}_{D}\partial_{\mathbf{n}}\phi\mathrm{d}s+\dfrac{\alpha_{R}}{2\beta_{R}}\int_{\Gamma_{R}}\phi^{2}\mathrm{d}s,
	\end{align}
    which is often more convenient for analysis and numerical discretization \cite{Liu2023DynamicMassPNP}. In the pure Neumann case, the electrostatic potential is determined only up to an additive constant; a gauge condition, such as prescribing the spatial average of $\phi$, is therefore needed for uniqueness.
\end{myremark}

\begin{myremark}
	In the whole-space setting $\Omega=\mathbb{R}^{d}$, the Poisson equation can be solved through the fundamental solution of $-\Delta$. Formally, under suitable decay and neutrality assumptions, the electrostatic potential is represented as \cite{Kinderlehrer2017WassersteinPNP}
	\begin{align}
		\phi=W*\dfrac{\psi}{\epsilon}, \quad \mathrm{where} \quad 
        W(x):=\left\{
		\begin{aligned}
			&-\dfrac{1}{2\pi}\mathrm{ln}|x|,\quad&&x\in\mathbb{R}^{2}/\left\{0\right\},\\
			&\dfrac{1}{d(d-2)\omega_{d}}|x|^{d-2},\quad &&x\in\mathbb{R}^{d}(d\geq 3),
		\end{aligned}
		\right.
	\end{align}
where $\omega_{d}$ is the volume of the unit ball. Eliminating $\phi$ in the free energy yields a nonlocal interaction term
	\begin{align}\label{eq:aggrefree}
		\mathcal{E}_{pnp}\left\{c_{i},\phi\right\}=\int_{\mathbb{R}^{d}}\sum_{i=1}^{N}c_{i}\log c_{i}\mathrm{d}\boldsymbol{x}+\dfrac{1}{2\epsilon}\int_{\mathbb{R}^{d}\times\mathbb{R}^{d}}W(x-y)\psi(x)\psi(y)\mathrm{d}\boldsymbol{x}\mathrm{d}\boldsymbol{y}.
	\end{align}
This entropy-interaction energy is the standard energy for aggregation-diffusion type equations: the entropy term drives diffusion, while the interaction kernel describes the electrostatic attraction or repulsion. Related gradient-flow structures also appear in models for biological swarming \cite{Burger2007Aggregation,Burger2008Large} and chemotaxis \cite{Burger2006Keller,Yoon2017Global}. In this setting, the classical PNP model can be interpreted as a Wasserstein gradient flow with respect to \eqref{eq:aggrefree} in the metric space of probability measures on $\mathbb{R}^{d}$; see \cite{Kinderlehrer2017WassersteinPNP} for the wellposedness of the JKO scheme and our previous work \cite{Carrillo2022PrimalDual} for the numerical realization.
\end{myremark}

\begin{myremark}
    Assume that the boundary data $\phi_{D}^{bc}$, $\phi_{N}^{bc}$, and $\phi_{R}^{bc}$ are independent of time and that the solution of \eqref{eq:classicalPNP} is sufficiently smooth. Then the following properties hold.
	\begin{itemize}
		\item Energy dissipation:
		\begin{align}
			\dfrac{\mathrm{d}\mathcal{E}_{pnp}}{\mathrm{d}t}=-\int_{\Omega}\sum_{i=1}^{N}D_{i}c_{i}|\nabla\mu_{i}|^{2}\mathrm{d}\boldsymbol{x}=-\int_{\Omega}\sum_{i=1}^{N}D_{i}c_{i}\Big|\nabla\Big(\log c_{i}+z_{i}\phi\Big)\Big|^2\mathrm{d}\boldsymbol{x}\leq0.
		\end{align}
        This identity follows from substituting $\partial_t c_i$ and applying integration by parts, where the boundary terms associated with the ionic chemical potential vanish due to the no-flux boundary conditions while the boundary terms associated with the electrostatic potential cancel with the boundary correction terms introduced in Remark.~\ref{rmk:Energy} \cite{Liu2023DynamicMassPNP}.
  
		\item Positivity: if the initial concentrations are nonnegative and sufficiently regular, then
		\begin{align}
			c_{i}(\boldsymbol{x},0)\geq0\implies c_{i}(\boldsymbol{x},t)\geq0,\quad\forall t>0,\quad i=1,\cdots,N.
		\end{align}
		\item Mass conservation: under periodic or no-flux boundary conditions~\eqref{eq:cinoflux} for $c_i$, each ionic species satisfies
		\begin{align}
			\int_{\Omega}c_{i}(\boldsymbol{x},t)\mathrm{d}\boldsymbol{x}=\int_{\Omega}c_{i}(\boldsymbol{x},0)\mathrm{d}\boldsymbol{x},\quad\forall t>0,\quad i=1,\cdots,N.
		\end{align}
	\end{itemize}
\end{myremark}

\subsection{Modified Poisson--Nernst--Planck Models}\label{sec:2.3}

The classical PNP model does not include steric interactions or concentration-gradient effects and therefore may be insufficient for describing concentrated electrolytes with finite ion sizes and strong spatial variation of ionic concentrations. Following modified PNP formulations with steric interactions and concentration-gradient corrections \cite{Qian2021PositivePNPCH,Ding2019PNPSteric,Gavish2018PNPSteric,Siddiqua2017ModifiedPNP},
we consider the free energy
\begin{align}\label{eq:Empnp_mu}
	\mathcal{E}_{mpnp}
	=\mathcal{E}_{pnp}
	+\dfrac{1}{2}\int_{\Omega}\boldsymbol{c}^{\mathrm{T}}G\boldsymbol{c}
	\mathrm{d}\boldsymbol{x}
	+\sum_{i=1}^{N}
	\int_{\Omega}\dfrac{\sigma_i}{2}|\nabla c_i|^2
	\mathrm{d}\boldsymbol{x},
	\quad
	\mu_i=\dfrac{\delta\mathcal{E}_{mpnp}}{\delta c_i}.
\end{align}
Here $\boldsymbol{c}=(c_1,\ldots,c_N)^{\mathrm{T}}$ and $G=(g_{ij})_{N\times N}$ is a symmetric positive-semidefinite matrix whose entry $g_{ij}$ measures the steric interaction between the $i$-th and $j$-th ionic species. The coefficient $\sigma_i\geq 0$ measures the strength of the concentration-gradient
correction for the $i$th ionic species.

The corresponding modified chemical potentials are
\begin{align}\label{eq:mu_modified}
	\mu_i=k_BT\log c_i+z_ie(\phi+\phi_e)
	+\sum_{j=1}^{N}g_{ij}c_j-\sigma_i\Delta c_i ,
\end{align}
where the additive constant from the entropy term has been omitted.

Using the gradient-flow form of the Nernst--Planck equation
\eqref{eq:NPequation}, we obtain the modified PNP model
\begin{align}\label{eq:modifiedPNPsecond}
	\dfrac{\partial c_i}{\partial t}
	=\nabla\cdot
	\Bigg(
	D_i\Big(
	\nabla c_i
	+\dfrac{z_ie c_i}{k_BT}\nabla(\phi+\phi_e)
	+\dfrac{c_i}{k_BT}\sum_{j=1}^{N}g_{ij}\nabla c_j
	-\dfrac{\sigma_i c_i}{k_BT}\nabla\Delta c_i
	\Big)
	\Bigg).
\end{align}
The modified model is supplemented with the same Poisson equation and electrostatic boundary conditions as the classical PNP model, together with the natural boundary conditions associated with the concentration-gradient term. Since the modification only changes the free-energy functional and the chemical potentials, the classical and modified PNP models can be treated in the same constrained JKO framework developed below.



\section{Variational schemes for Poisson-Nernst-Planck models}\label{sec:3}
This section develops a variational scheme for PNP models based on optimal transport. Following the dynamic JKO formulation and the primal--dual framework in \cite{Carrillo2022PrimalDual,deng2025PDFB,Wu2026PDPFSMCL,Carrillo2024StructurePD}, the scheme is designed to preserve energy dissipation, positivity of ionic concentrations, and mass conservation. For clarity, we use the non-dimensionalized setting for PNP models: $k_{B}T=1$, $z_{p}=1$, $z_{n}=-1$, $e=1$, $D_{p}=D_{n}=1$, and the external potential $\phi_{e}=0$. The extension to the general setting is straightforward.

\subsection{Semi-discrete variational scheme for PNP models}\label{sec:semi_discrete}
We first recall the variational formulation of a Wasserstein gradient flow. For a single density $\rho$, the gradient flow of an energy functional $\mathcal{E}$ can be written in the continuity-equation form
\begin{align}\label{eq:general_WGF}
	\dfrac{\partial\rho}{\partial t}=-\nabla\cdot(\rho\boldsymbol{v}),\quad
	\boldsymbol{v}=-\nabla\dfrac{\delta\mathcal{E}}{\delta\rho}.
\end{align}
Given a time step $\tau>0$, the classical JKO scheme constructs a sequence $\left\{\rho^{k}(x)\right\}$ approximating $\rho(x,t_{k})$ with $t_{k}=k\tau$ through
\begin{align}\label{eq:JKO_scalar}
	\rho^{k+1}\in\mathop{\text{arg min}}\limits_{\rho}
	\Bigl\{d^{2}_{\mathcal{W}}\big(\rho,\rho^{k}\big)+2\tau\mathcal{E}(\rho)\Bigr\}.
\end{align}
Under suitable assumptions, such variational schemes are well posed, and the time-discrete sequences converge weakly to the solutions to the gradient flows as $\tau\rightarrow0$ \cite{Jordan1998Variational,Carrillo2010NonlinearMobility,Lisini2012CahnHilliard}.

We now specialize the above construction to the two-species PNP model, which can be rewritten as
\begin{align}\label{eq:dtpn}
	\left\{
		\begin{aligned}
		&\dfrac{\partial}{\partial t}
		\begin{bmatrix}
			p \\ 
			n
		\end{bmatrix}  
		=\nabla\cdot \Bigg(
		\begin{bmatrix}
			p & 0\\
			0 & n
		\end{bmatrix}
		\begin{bmatrix}
			\nabla\dfrac{\delta\mathcal{E}}{\delta p} \\ \nabla\dfrac{\delta\mathcal{E}}{\delta n}
		\end{bmatrix}\Bigg),\\
		&-p+n-\epsilon\Delta\phi=\psi^{0}.
		\end{aligned}
		\right.
\end{align}
The essential difference of PNP model from the scalar Wasserstein gradient flow is the Poisson coupling: the transport distance acts on the ionic concentrations, while the electrostatic potential involved in the energy must be determined simultaneously from the Poisson equation. Thus, the semi-discrete JKO step for PNP models can be formulated as a variational problem with the Poisson equation imposed as an additional constraint:
\begin{align}\label{eq:JKO_PNP}
	\begin{aligned}
		&(\boldsymbol{\rho}^{k+1},\phi^{k+1})\in
		\mathop{\text{arg min}}\limits_{\boldsymbol{\rho},\phi}
		\Bigl\{d^{2}_{\mathcal{W}}\big(\boldsymbol{\rho},\boldsymbol{\rho}^{k}\big)
		+2\tau\mathcal{E}(\boldsymbol{\rho},\phi)\Bigr\},\\
		&\text{s.t.}
		\left\{\begin{aligned}
			&-p+n-\epsilon\Delta\phi=\psi^{0}\quad&&\text{in $\Omega$},\\
			&\alpha\phi+\beta\epsilon\dfrac{\partial\phi}{\partial\mathbf{n}}=\phi^{bc}\quad&&\text{on $\partial\Omega$},
		\end{aligned}
		\right.
	\end{aligned}
\end{align}
where we define the vector density $\boldsymbol{\rho}=(p,n)^{\mathrm{T}}$. 
Following the Benamou--Brenier dynamic characterization of optimal transport distances~\cite{Benamou2000Computational}, we introduce an auxiliary transport time variable $s\in[0,1]$ and define the transport distance for the vector density by
\begin{align}\label{eq:dwpn}
    d_{\mathcal{W}}(\boldsymbol{\rho}_{0},\boldsymbol{\rho}_{1})
    =\mathop{\mathrm{inf}}\limits_{(\boldsymbol{\rho},\textbf{m}_{\boldsymbol{\rho}})}
    \Bigl\{\int_{0}^{1}\int_\Omega\Big(\mathcal{G}(p,\boldsymbol{m}_{p})+\mathcal{G}(n,\boldsymbol{m}_{n})\Big)\mathrm{d}\boldsymbol{x}\mathrm{d}s\Bigr\}^{\frac{1}{2}},
\end{align}
where $\textbf{m}_{\boldsymbol{\rho}}=(\boldsymbol{m}_{p},\boldsymbol{m}_{n})^{\mathrm{T}}$ denotes the momentum variable and the action function is defined as
\begin{align}\label{eq:Wdistance}
	\mathcal{G}(c,\boldsymbol{m})=\left\{\begin{aligned}
		&\dfrac{|\boldsymbol{m}|^{2}}{c}\quad&&c>0,\\
		&0 \quad\quad&&(c,\boldsymbol{m})=(0,\boldsymbol{0}),\\
		&+\infty \quad&&\text{otherwise}.
	\end{aligned}
	\right.
\end{align}
The admissible paths for \eqref{eq:dwpn} satisfy
\begin{align}\label{eq:constraint}
	\left\{\begin{aligned}
&\partial_{s}\boldsymbol{\rho}+\nabla\cdot\textbf{m}_{\boldsymbol{\rho}}=0\quad&&(\boldsymbol{x},s)\in\Omega\times[0,1],\\
		&\textbf{m}_{\boldsymbol{\rho}}\cdot\mathbf{n}=0\quad&&(\boldsymbol{x},s)\in\partial\Omega\times[0,1],\\
		&\boldsymbol{\rho}(\boldsymbol{x},0)=\boldsymbol{\rho}_{0},\quad \boldsymbol{\rho}(\boldsymbol{x},1)=\boldsymbol{\rho}_{1}\quad&&\boldsymbol{x}\in\Omega.
	\end{aligned}
	\right.
\end{align}
The above constraints of continuity equations define only the transport distance between two concentration states. For PNP models, the Poisson equation should be imposed, where the endpoint concentrations and the electrostatic potential are coupled. We thus arrive at the following variational formulation for the JKO step of PNP models.

\textbf{Problem 1 (Semi-discrete Poisson-constrained JKO scheme).} Given $(p^{k},n^{k},\phi^{k})$, find $u^{k+1}=(p,\boldsymbol{m}_{p};n,\boldsymbol{m}_{n};\phi)$ by solving
\begin{align}\label{eq:semiJKO}
	\begin{aligned}
		&u^{k+1}\in\mathop{\text{arg min}}\limits_{u}\int_{\Omega}\Big(\mathcal{G}(p,\boldsymbol{m}_{p})+\mathcal{G}(n,\boldsymbol{m}_{n})\Big)\mathrm{d}\boldsymbol{x}+2\tau\mathcal{E}(p,n,\phi),\\
		&\text{s.t.}\left\{\begin{aligned}
			&p+\nabla\cdot\boldsymbol{m}_{p}=p^{k},\quad n+\nabla\cdot\boldsymbol{m}_{n}=n^{k},\quad&&\text{in $\Omega$}\\
			&\boldsymbol{m}_{p}\cdot\mathbf{n}=0,\quad \boldsymbol{m}_{n}\cdot\mathbf{n}=0,\quad&&\text{on $\partial\Omega$}\\
			&-p+n-\epsilon\Delta\phi=\psi^{0}\quad&&\text{in $\Omega$}\\
			&\alpha\phi+\beta\epsilon\dfrac{\partial\phi}{\partial\mathbf{n}}=\phi^{bc}\quad&&\text{on $\partial\Omega$}.
		\end{aligned}
		\right.
	\end{aligned}
\end{align}
This formulation combines the Benamou--Brenier transport representation with the Poisson constraint and uses a one-step discretization in the auxiliary variable $s$. Namely, the action is approximated by the right endpoint rule and the continuity equation is discretized with $\Delta s=1$. This is consistent with the first-order accuracy of the outer JKO step and does not reduce the temporal accuracy of the scheme \cite{LiWuchen2020Fisher}. The resulting variational formulation~\eqref{eq:semiJKO} has a structure of convex objectives with linear PDE constraints, which is the basis of implementing of efficient primal--dual methods developed in \cite{Carrillo2024StructurePD,Wu2026PDPFSMCL,Zeng2026variable}.

\begin{myremark}
The existence and uniqueness of minimizers for the variational scheme of PNP models have been established under suitable assumptions in \cite{Kinderlehrer2017WassersteinPNP,Liu2023DynamicMassPNP}. However, for multi-species PNP models with anisotropic diffusion and complex boundary conditions, the convergence of the JKO solutions to solutions of the continuous PDE system remains an open problem.
\end{myremark}

\subsection{Fully discrete variational scheme for PNP models}\label{sec:full_discrete}
In this subsection, we construct the fully discrete variational scheme for PNP models, which applies to both the classical PNP model and the modified models by choosing the steric-interaction matrix $G$ and the coefficients $\sigma_{p},\sigma_{n}\geq 0$ appropriately. We first introduce the spatial grid and the discrete operators used in the following formulation. Let $\Omega=[a,b]\times[c,d]$ and divide it into a uniform cell-centered grid with $N_{x}\times N_{y}$ cells. The mesh sizes are
\begin{align}
	\Delta x=\dfrac{b-a}{N_{x}},\quad\Delta y=\dfrac{d-c}{N_{y}}.
\end{align}
The cell centers take the form
\begin{align}
	x_{i}=a+\left(i-\dfrac{1}{2}\right)\Delta x,\quad
	y_{j}=c+\left(j-\dfrac{1}{2}\right)\Delta y,
	\quad
	1\leq i\leq N_{x},\quad 1\leq j\leq N_{y}.
\end{align}
The cell $C_{i,j}$ is given as
\begin{align}
	C_{i,j}=
	\left[x_{i}-\dfrac{\Delta x}{2},x_{i}+\dfrac{\Delta x}{2}\right]
	\times
	\left[y_{j}-\dfrac{\Delta y}{2},y_{j}+\dfrac{\Delta y}{2}\right].
\end{align}
For a cell-centered grid function $\left\{\rho_{i,j}\right\}$, we introduce the following average and difference operators:
\begin{align}
    \left\{\begin{aligned}
	&\mathcal{A}_{x}\rho_{i+\frac{1}{2},j}=\dfrac{\rho_{i+1,j}+\rho_{i,j}}{2},\quad\mathcal{A}_{y}\rho_{i,j+\frac{1}{2}}=\dfrac{\rho_{i,j+1}+\rho_{i,j}}{2},\\
    &\mathcal{D}_{x}\rho_{i+\frac{1}{2},j}=\dfrac{\rho_{i+1,j}-\rho_{i,j}}{\Delta x},\quad
	\mathcal{D}_{y}\rho_{i,j+\frac{1}{2}}=\dfrac{\rho_{i,j+1}-\rho_{i,j}}{\Delta y}.
    \end{aligned}
	\right.
\end{align}

For the concentration-gradient energy terms, we impose the following
homogeneous natural boundary conditions:
\begin{align}\label{eq:discretebcpn}
    &\mathcal{D}_{x}\rho_{\frac{1}{2},j}=\mathcal{D}_{x}\rho_{N_{x}+\frac{1}{2},j}=\mathcal{D}_{y}\rho_{i,\frac{1}{2}}=\mathcal{D}_{y}\rho_{i,N_{y}+\frac{1}{2}}=0,\quad \text{for $\rho=p,n$}.
\end{align}
The no-flux condition associated with the JKO transport variables will be imposed directly on the discrete fluxes in the transport constraints below.

For the electrostatic potential, the contributions induced by nonhomogeneous boundary conditions must be retained when discretizing the free energy (see the boundary integral terms in \eqref{eq:totalEnePNP} in Remark~\ref{rmk:Energy}). For illustrative purposes, we consider Dirichlet boundary conditions on the left and right boundaries and Neumann boundary conditions on the bottom and top boundaries, namely,
\begin{align}
    &\left.\phi\right|_{x=a}=\phi_a,\quad
    \left.\phi\right|_{x=b}=\phi_b,\quad
    \left.\epsilon\partial_{\mathbf{n}}\phi\right|_{y=c}=\phi_c,\quad
    \left.\epsilon\partial_{\mathbf{n}}\phi\right|_{y=d}=\phi_d,
\end{align}
where $\phi_a,\phi_b,\phi_c,$ and $\phi_d$ denote the right-hand-side terms of the corresponding boundary conditions. 

The classical free energy \eqref{eq:totalEnePNP} and the local steric-interaction energy in \eqref{eq:Empnp_mu} are discretized by the midpoint rule. For the concentration-gradient energy, the derivatives $\mathcal{D}_{x}\rho$ and $\mathcal{D}_{y}\rho$ are evaluated at cell edges and the integral of $|\nabla\rho|^{2}$ is approximated by the trapezoidal rule over these edge values. Because the homogeneous Neumann condition sets the boundary normal derivatives to zero, the boundary edge contributions vanish, and only the interior edge sums appear below. Write the symmetric steric-interaction matrix as $G=(g_{\rho\eta})_{\rho,\eta=p,n}$. This leads to the following discrete energy $\mathcal{E}^{h}$:
\begin{eqnarray}\label{eq:discreteEnergy}
	\begin{aligned}
\mathcal{E}^{h}&=\sum_{i=1}^{N_{x}}\sum_{j=1}^{N_{y}}\Bigg(\sum_{\rho=p,n}\Big(\rho_{i,j}\log \rho_{i,j}+\dfrac{1}{2}\sum_{\eta=p,n}g_{\rho\eta}\rho_{i,j}\eta_{i,j}\Big)+\dfrac{1}{2}(p_{i,j}-n_{i,j}+\psi^{0}_{i,j})\phi_{i,j}\Bigg)\Delta x\Delta y\\
		&-\dfrac{1}{2}\sum_{j=1}^{N_{y}}\epsilon\Big(-\phi_{a}\mathcal{D}_{x}\phi_{\frac{1}{2},j}+\phi_{b}\mathcal{D}_{x}\phi_{N_{x}+\frac{1}{2},j}\Big)\Delta y+\dfrac{1}{2}\sum_{i=1}^{N_{x}}\Big(\phi_{c}\mathcal{A}_{y}\phi_{i,\frac{1}{2}}+\phi_{d}\mathcal{A}_{y}\phi_{i,N_{y}+\frac{1}{2}}\Big)\Delta x\\
		&+\sum_{\rho=p,n}\dfrac{\sigma_{\rho}}{2}\Big(\sum_{j=1}^{N_{y}}\sum_{i=1}^{N_{x}-1}(\mathcal{D}_{x}\rho_{i+\frac{1}{2},j})^{2}+\sum_{i=1}^{N_{x}}\sum_{j=1}^{N_{y}-1}(\mathcal{D}_{y}\rho_{i,j+\frac{1}{2}})^{2}\Big)\Delta x\Delta y.
	\end{aligned}
\end{eqnarray}
When treating $(p_{i,j},n_{i,j},\phi_{i,j})$ as independent variables and Poisson equation as a constraint in the proposed JKO scheme, the corresponding energy gradient is given by:
\begin{eqnarray}\label{eq:discretegradient}
	\begin{aligned}
		&\nabla\mathcal{E}^{h}_{i,j}:=\Big(\partial_{p_{i,j}}\mathcal{E}^{h},\mathbf{0},\mathbf{0},
		\partial_{n_{i,j}}\mathcal{E}^{h},\mathbf{0},\mathbf{0},\partial_{\phi_{i,j}}\mathcal{E}^{h}\Big)^{\mathrm{T}},\\
		&\partial_{p_{i,j}}\mathcal{E}^{h}=\Big(\log p_{i,j}+1+g_{pp}p_{i,j}+g_{pn}n_{i,j}+\dfrac{1}{2}\phi_{i,j}-\sigma_{p}\mathcal{D}^2p_{i,j}\Big)\Delta x\Delta y,\\
		&\partial_{n_{i,j}}\mathcal{E}^{h}=\Big(\log n_{i,j}+1+g_{pn}p_{i,j}+g_{nn}n_{i,j}-\dfrac{1}{2}\phi_{i,j}-\sigma_{n}\mathcal{D}^2n_{i,j}\Big)\Delta x\Delta y,\\
		&\partial_{\phi_{i,j}}\mathcal{E}^{h}=\Big(\dfrac{1}{2}(p_{i,j}-n_{i,j}+\psi^{0}_{i,j})\Big)\Delta x\Delta y+W^{y}_{i,j}\Delta y+W^{x}_{i,j}\Delta x.
	\end{aligned}
\end{eqnarray}
where $\mathcal{D}^{2}\rho_{i,j}=\mathcal{D}^{2}_{x}\rho_{i,j}+\mathcal{D}^2_{y}\rho_{i,j}$ represents the Laplace operator discretization with embedded boundary conditions \eqref{eq:discretebcpn}:
\begin{eqnarray}\label{eq:disLaplace}
	\mathcal{D}^2_{x}\rho_{i,j}=\left\{
	\begin{aligned}
		&\dfrac{\rho_{2,j}-\rho_{1,j}}{\Delta x^{2}}&&i=1,\\
		&\dfrac{\rho_{i+1,j}-2\rho_{i,j}+\rho_{i-1,j}}{\Delta x^{2}}&&\text{others},\\
		&\dfrac{\rho_{i-1,j}-\rho_{i,j}}{\Delta x^{2}}&&i=N_{x}.
	\end{aligned}
	\right.\qquad
	\mathcal{D}^2_{y}\rho_{i,j}=\left\{
	\begin{aligned}
		&\dfrac{\rho_{i,2}-\rho_{i,1}}{\Delta y^{2}}&&j=1,\\
		&\dfrac{\rho_{i,j+1}-2\rho_{i,j}+\rho_{i,j-1}}{\Delta y^{2}}&&\text{others},\\
		&\dfrac{\rho_{i,j-1}-\rho_{i,j}}{\Delta y^{2}}&&j=N_{y}.
	\end{aligned}
	\right.
\end{eqnarray}
The quantities $W^{y}_{i,j}$ and $W^{x}_{i,j}$ in \eqref{eq:discretegradient} collect the contributions from the boundary energy of $\phi$:
\begin{eqnarray}\label{eq:Wij}
	\begin{aligned}
		W^{y}_{i,j}=\left\{
		\begin{aligned}
			&\epsilon\dfrac{\phi_{a}}{\Delta x}\quad&&i=1,\\ 
			&\epsilon\dfrac{\phi_{b}}{\Delta x}\quad&&i=N_{x},\\
			&0\quad&&\text{otherwise},
		\end{aligned}
		\right.\quad
		W^{x}_{i,j}=\left\{
		\begin{aligned}
			&\dfrac{1}{2}\phi_{c}\quad&&j=1,\\ 
			&\dfrac{1}{2}\phi_{d}\quad&& j=N_{y},\\
			&0\quad&&\text{otherwise}.
			\end{aligned}
		\right.
	\end{aligned}
\end{eqnarray}

The integral of action function $\mathcal{G}(\rho,\boldsymbol{m})$ in \eqref{eq:semiJKO} is approximated by the midpoint rule, and the continuity equation is discretized with a centered difference scheme as in \cite{Carrillo2024StructurePD,Wu2026PDPFSMCL}. The resulting fully discrete variational formulation is given below.

\textbf{Problem 2 (Fully discrete Poisson-constrained JKO scheme).} Given the current numerical solution $\left\{p^{k}_{i,j},n^{k}_{i,j},\phi^{k}_{i,j}\right\}$ at $t=k\tau$, we solve $\left\{p^{k+1}_{i,j},n^{k+1}_{i,j},\phi^{k+1}_{i,j}\right\}$ at next time step through the following scheme, for $1\leq i\leq N_x$ and $1\leq j\leq N_y$:
\begin{eqnarray}\label{eq:fullJKO}
	\begin{aligned}
		&\{u^{k+1}_{i,j}\}\in\mathop{\text{arg min}}\limits_{u}\sum_{i=1}^{N_{x}}\sum_{j=1}^{N_y}\Big(\mathcal{G}(p_{i,j},{\boldsymbol{m}_{p}}_{i,j})+\mathcal{G}(n_{i,j},{\boldsymbol{m}_{n}}_{i,j})\Big)\Delta x\Delta y +2\tau\mathcal{E}^{h}(p,n,\phi)\\
		&\text{s.t.}\left\{
		\begin{aligned}
			&p_{i,j}+\dfrac{1}{2\Delta x}\big((m^{x}_{p})_{i+1,j}-(m^{x}_{p})_{i-1,j}\big)+\dfrac{1}{2\Delta y}\big((m^{y}_{p})_{i,j+1}-(m^{y}_{p})_{i,j-1}\big)=p^{k}_{i,j},\\
			&n_{i,j}+\dfrac{1}{2\Delta x}\big((m^{x}_{n})_{i+1,j}-(m^{x}_{n})_{i-1,j}\big)+\dfrac{1}{2\Delta y}\big((m^{y}_{n})_{i,j+1}-(m^{y}_{n})_{i,j-1}\big)=n^{k}_{i,j},\\
			&(m^{x}_{p})_{0,j}=-(m^{x}_{p})_{1,j},(m^{x}_{p})_{N_{x}+1,j}=-(m^{x}_{p})_{N_{x},j},(m^{y}_{p})_{i,0}=-(m^{y}_{p})_{i,1},(m^{y}_{p})_{i,N_{y}+1}=-(m^{y}_{p})_{i,N_{y}},\\
			&(m^{x}_{n})_{0,j}=-(m^{x}_{n})_{1,j},(m^{x}_{n})_{N_{x}+1,j}=-(m^{x}_{n})_{N_{x},j},(m^{y}_{n})_{i,0}=-(m^{y}_{n})_{i,1},(m^{y}_{n})_{i,N_{y}+1}=-(m^{y}_{n})_{i,N_{y}},\\
			&-p_{i,j}+n_{i,j}+\epsilon L_{\phi}\phi_{i,j}=\psi^{0}_{i,j}+F^{bc}_{i,j}.
		\end{aligned}
		\right.
	\end{aligned}
\end{eqnarray}
Here $L_{\phi}$ denotes the discrete Laplace operator with the prescribed boundary conditions for the electrostatic potential, and $F^{bc}$ collects the corresponding boundary terms. 

In the following theorem, we prove that the proposed fully discrete JKO scheme preserves the desired structural properties at the discrete level.
\begin{mytheorem}\label{Th:EMB}
	The full-discrete variational scheme has the following structure-preserving properties:
	
	(i) Original energy dissipation structure;
	
	(ii) Mass conservation of $p$ and $n$;
	
	(iii) Positivity preserving of $p$ and $n$ in the sense that $p^{k+1}_{i,j}\geq 0$ and $n^{k+1}_{i,j}\geq 0$.
    
Proof: (i) Let
$u^{k+1}=\left\{p^{k+1},\boldsymbol{m}_{p}^{k+1};n^{k+1},\boldsymbol{m}_{n}^{k+1};\phi^{k+1}\right\}$
be the optimal solution of the minimization problem \eqref{eq:fullJKO}. Since the previous state
$\left\{p^{k},\boldsymbol{0};n^{k},\boldsymbol{0};\phi^{k}\right\}$ is an admissible competitor, provided that they satisfy the discrete Poisson constraint, the optimality of $u^{k+1}$ yields
\begin{align}
	\sum_{i=1}^{N_{x}}\sum_{j=1}^{N_y}\Big(\mathcal{G}(p^{k+1}_{i,j},{\boldsymbol{m}^{k+1}_{p_{i,j}}})+\mathcal{G}(n^{k+1}_{i,j},{\boldsymbol{m}^{k+1}_{n_{i,j}}})\Big)\Delta x\Delta y+2\tau\mathcal{E}^{h}(p^{k+1},n^{k+1},\phi^{k+1})\leq2\tau\mathcal{E}^{h}(p^{k},n^{k},\phi^{k}).
\end{align}
This gives the discrete counterpart of the original energy dissipation law in Sec.~\ref{sec:PNP}. Since the action terms are nonnegative, we obtain
\begin{align}
	\mathcal{E}^{h}(p^{k+1},n^{k+1},\phi^{k+1})\leq\mathcal{E}^{h}(p^{k},n^{k},\phi^{k}).
\end{align}	

(ii) We prove the mass conservation of $p$, and the proof for $n$ is identical. Summing the discrete continuity constraint for $p$ in \eqref{eq:fullJKO} over $i$ and $j$ gives
\begin{align}
	&\sum_{i=1}^{N_x}\sum_{j=1}^{N_y}\Big(p_{i,j}+\dfrac{1}{2\Delta x}\big((m^{x}_{p})_{i+1,j}-(m^{x}_{p})_{i-1,j}\big)+\dfrac{1}{2\Delta y}\big((m^{y}_{p})_{i,j+1}-(m^{y}_{p})_{i,j-1}\big)\Big)=\sum_{i=1}^{N_x}\sum_{j=1}^{N_y}p^{k}_{i,j}.
\end{align}
By summing over the flux variables $m^{x}_{p}$, $m^{y}_{p}$ and applying the no-flux boundary conditions in \eqref{eq:fullJKO}, we obtain
\begin{eqnarray}\label{eq:discrete_m}
	\begin{aligned}
		&\sum_{i=1}^{N_{x}}\big((m^{x}_{p})_{i+1,j}-(m^{x}_{p})_{i-1,j}\big)
		=-\big((m^{x}_{p})_{0,j}+(m^{x}_{p})_{1,j}\big)+\big((m^{x}_{p})_{N_{x},j}+(m^{x}_{p})_{N_{x}+1,j}\big)=0,\quad 1\leq j\leq N_{y},\\
		&\sum_{j=1}^{N_{y}}\big((m^{y}_{p})_{i,j+1}-(m^{y}_{p})_{i,j-1}\big)
		=-\big((m^{y}_{p})_{i,0}+(m^{y}_{p})_{i,1}\big)+\big((m^{y}_{p})_{i,N_{y}}+(m^{y}_{p})_{i,N_{y}+1}\big)=0,\quad 1\leq i\leq N_{x}.
	\end{aligned}
\end{eqnarray}
Thus, the discrete mass of $p$ is conserved:
\begin{align}
	\sum_{i=1}^{N_x}\sum_{j=1}^{N_y}p_{i,j}\Delta x\Delta y
	=
	\sum_{i=1}^{N_x}\sum_{j=1}^{N_y}p^{k}_{i,j}\Delta x\Delta y.
\end{align}
The same argument gives the mass conservation of $n$.

(iii) By the definition of the action function $\mathcal{G}$, the admissible set of \eqref{eq:fullJKO} requires $p_{i,j}\geq 0$ and $n_{i,j}\geq 0$. Hence, the minimizer of the JKO scheme remains in the nonnegative admissible set, and the scheme is positivity-preserving in the sense that $p^{k+1}_{i,j}\geq 0$ and $n^{k+1}_{i,j}\geq 0$.
\qed
\end{mytheorem}

\section{Primal--dual splitting methods}\label{sec:4}
In this section we briefly introduce two primal--dual splitting methods for solving Problem~2~\eqref{eq:fullJKO}. The full-discrete JKO scheme can be formulated as a linearly constrained convex problem in the following form:
	\begin{eqnarray}\label{eq:minproblem}
		\begin{aligned}
			\mathop{\mathrm{min}}\limits_{u}\Psi(u)+E(u),\quad \text{s.t.}\quad Au=b.
		\end{aligned}
	\end{eqnarray}
where we define
\begin{eqnarray}
	\begin{aligned}
		\left\{\begin{aligned}
			&u=\left\{p_{i,j},{m^{x}_{p}}_{i,j}, {m^{y}_{p}}_{i,j},n_{i,j},{m^{x}_{n}}_{i,j},{m^{y}_{n}}_{i,j},\phi_{i,j}\right\}^{1\leq j\leq N_{y}}_{1\leq i\leq N_{x}},\\
			&\Psi(u)=\sum_{i=1}^{N_{x}}\sum_{j=1}^{N_y}\Big(\mathcal{G}(p_{i,j},{\boldsymbol{m}_{p}}_{i,j})+\mathcal{G}(n_{i,j},{\boldsymbol{m}_{n}}_{i,j})\Big)\Delta x\Delta y,\\
			&E(u)=2\tau\mathcal{E}^{h}(p,n,\phi),
		\end{aligned}
		\right.\\
	\end{aligned}
\end{eqnarray}
and the linear system $Au=b$ corresponds to the discretized constraints in \eqref{eq:fullJKO}, whose explicit form will be given in Section~\ref{sec:5}. 
By introducing a relaxed penalty term for the constraints, the problem~\eqref{eq:minproblem} can be transformed into an unconstrained optimization:
\begin{eqnarray}\label{eq:unconstrained}
	\begin{aligned}
		\mathop{\mathrm{min}}_{u}\Psi(u)+E(u)+\iota_\delta(Au),\quad
		\iota_\delta(z)=\left\{
		\begin{aligned}
			&0\quad\Vert z-b\Vert_{2}\leq\delta, \\
			&\infty\quad\text{otherwise},
		\end{aligned}
		\right.
	\end{aligned}
\end{eqnarray}
where $\iota_\delta(z)$ is the indicator function of the set $\{z:\Vert z-b\Vert_{2}\leq\delta\}$, and $\delta\geq 0$ is a relaxation parameter that allows for a small violation of the constraints. This formulation enables the use of efficient proximal algorithms to solve the optimization problem while ensuring that the constraints are approximately satisfied within a prescribed tolerance dependent on the order of truncation errors \cite{Carrillo2022PrimalDual}. Using the duality relation $\iota_{\delta}(Au)=\max_{v}\langle Au,v\rangle-\iota^{*}_{\delta}(v)$, where $\iota^{*}_{\delta}(v)=\langle b,v\rangle+\delta\|v\|_2$, \eqref{eq:unconstrained} can be rewritten as the saddle-point problem
\begin{align}\label{eq:saddlePD}
	\begin{aligned}
		\mathop{\mathrm{min}}_{u}\mathop{\mathrm{max}}_{v}\,
		\Psi(u)+E(u)+\langle Au,v\rangle-\iota^{*}_{\delta}(v).
	\end{aligned}
\end{align}
The first-order optimality system can be viewed as the steady state of the following primal--dual flow:
\begin{align}\label{eq:PDflow}
	\begin{bmatrix}
		u'\\
		v'
	\end{bmatrix}
	\in
	\begin{bmatrix}
		-(\partial\Psi+\nabla E) & -A^{\mathrm T}\\
		A & -\partial\iota^{*}_{\delta}
	\end{bmatrix}
	\begin{bmatrix}
		u\\
		v
	\end{bmatrix},
\end{align}
where 
$\partial\Psi$ and $\partial\iota_\delta^*$ represent the subdifferential of the nonsmooth functionals $\Psi(u)$ and $\iota_\delta^*(v)$, respectively. In practice, we apply proximal gradient methods to solve the saddle-point problem \eqref{eq:saddlePD} by discretizing the primal--dual flow \eqref{eq:PDflow}, where the proximal operator is a fundamental tool in these algorithms. In particular, the generalized proximal operator with respect to a metric induced by a symmetric positive definite matrix $T$ is defined as follows
\begin{align}\label{eq:defprox}
	\mathrm{Prox}^{T}_{h}(w)=
	\mathop{\mathrm{argmin}}\limits_{z}
	\left\{\dfrac{1}{2}\Vert z-w\Vert^{2}_{T}+h(z)\right\},
	\qquad
	\Vert z\Vert^{2}_{T}=z^{\mathrm{T}}Tz.
\end{align}

\subsection{Preconditioned primal--dual method}
By leveraging modern operator splitting techniques, the problem~\eqref{eq:unconstrained} can be solved by primal--dual methods for three operators (PD3O)~\cite{Yan2018NewPDThreeFunctions}. However, direct PD3O-type discretizations of \eqref{eq:PDflow} suffer from severe step-size restrictions and slow convergence for dynamic JKO minimization problems \cite{Carrillo2024StructurePD}. The preconditioned primal--dual (PrePD) method \cite{Carrillo2024StructurePD} is then proposed following the acceleration strategy in \cite{Liu2021AccelerationPrimalDual}, which introduces block-diagonal preconditioners $T_u$ and $T_v$ for the primal and dual flows, respectively. This yields the following iteration:
\begin{eqnarray}\label{eq:PrePD}
	\begin{aligned}
		\begin{cases}
			v^{(\ell+1)}=\mathrm{Prox}^{T_{v}}_{\lambda^{-1}\iota^{*}_{\delta}}
			\Big(v^{(\ell)}+\lambda^{-1}T^{-1}_{v} A\bar{u}^{(\ell)}\Big),\\
			u^{(\ell+1)}=\mathrm{Prox}^{T_{u}}_{\lambda\Psi}
			\Big(u^{(\ell)}-\lambda T^{-1}_{u}\nabla E(u^{(\ell)})-\lambda T^{-1}_{u} A^{\mathrm{T}}v^{(\ell+1)}\Big),\\
			\bar{u}^{(\ell+1)}=2u^{(\ell+1)}-u^{(\ell)}
			+\lambda T^{-1}_{u}\nabla E(u^{(\ell)})-\lambda T^{-1}_{u}\nabla E(u^{(\ell+1)}).
		\end{cases}
	\end{aligned}
\end{eqnarray}
where $T_u$ and $T_v$ are chosen such that $T_{v}\succeq A T^{-1}_{u}A^{\mathrm{T}}$ to guarantee convergence \cite{Liu2021AccelerationPrimalDual}. In the PrePD method, we take $T_u=I$ and $T_v=AA^{\mathrm T}$ with step size $\lambda>0$ to improve the conditioning of the saddle-point system \eqref{eq:PDflow} and accelerates convergence. With this choice, the primal update is a local proximal step for the transport action, while the dual update contains the main global linear solve. The corresponding proximal operators are summarized in Section~\ref{sec:proximal_operator}.

The stopping criteria for the convergence of the iterative update for the primal variable $u$ and the dual variable $v$ consists of the constraints \eqref{eq:fullJKO} and the relative error for energy and transport distance 
\begin{eqnarray}\label{eq:stop}
	\begin{aligned}
 		e^{(\ell+1)}_{A}=\Vert Au^{(\ell+1)}-b\Vert_{2}\leq\delta,\quad e^{(\ell+1)}_{u}\leq\epsilon_{1},\quad e^{(\ell+1)}_{v}\leq\epsilon_{1},\quad e^{(\ell+1)}_{E}\leq\epsilon_{2},\quad e^{(\ell+1)}_{\Psi}\leq\epsilon_{2}.
	\end{aligned}
\end{eqnarray}

\subsection{Variable preconditioned transformed primal--dual method}
To further improve the computational efficiency and make large-scale simulations feasible, we have extended the recently proposed transformed primal--dual (TPD) methodology \cite{Chen2023TPD,Chen2025Transformed} to the variable preconditioned transformed primal--dual (VPTPD) method for the JKO scheme \cite{Zeng2026variable}. Compared with the block-diagonal preconditioning in PrePD, VPTPD first applies a Schur-complement-based block-triangular transformation to the primal--dual flow and then a variable-dependent preconditioning to the transformed system. Formally, the transformed and preconditioned flow associated with \eqref{eq:PDflow} can be written as
\begin{align}\label{eq:VPTPDflow}
	\begin{bmatrix}
		u'\\
		v'
	\end{bmatrix}
	&\in
	\underbrace{
	\begin{bmatrix}
		T_u^{-1} & 0\\
		0 & T_v^{-1}
	\end{bmatrix}}_{\text{preconditioning}}
	\underbrace{
	\begin{bmatrix}
		I & 0\\
		A T_u^{-1} & I
	\end{bmatrix}}_{\text{transformation}}
	\begin{bmatrix}
		-\partial\mathcal H & -A^{\mathrm T}\\
		A & -\partial\iota^{*}_{\delta}
	\end{bmatrix}
	\begin{bmatrix}
		u\\
		v
	\end{bmatrix}\notag\\
	&=
	\begin{bmatrix}
		T_u^{-1} & 0\\
		0 & T_v^{-1}
	\end{bmatrix}
	\begin{bmatrix}
		-\partial\mathcal H & -A^{\mathrm T}\\
		A-A T_u^{-1}\partial\mathcal H & -A T_u^{-1}A^{\mathrm T}-\partial\iota^{*}_{\delta}
	\end{bmatrix}
	\begin{bmatrix}
		u\\
		v
	\end{bmatrix},
\end{align}
where we define $\mathcal H(u)=\Psi(u)+E(u)$ and $\partial\mathcal H=\partial\Psi+\nabla E$. If $T_u$ approximates the local primal gradient operator $\partial\mathcal H$ (or, in the smooth case, the local Hessian of $\mathcal H$), then the lower-left block $A-A T_u^{-1}\partial\mathcal H$ becomes small, and the transformed saddle-point system is nearly upper triangular. Moreover, the lower-right block contains the Schur-complement-type operator $A T_u^{-1}A^{\mathrm T}$, which introduces strong monotonicity in the dual variable and therefore can significantly accelerate convergence. In the present JKO problem, the nonsmooth action $\Psi$ is treated by a proximal step and the smooth energy $E$ is treated explicitly. This gives the following semi-implicit-explicit VPTPD iteration:
\begin{eqnarray}\label{eq:VPTPD}
	\begin{aligned}
		\begin{cases}
			u^{(\ell+1)}=\mathrm{Prox}^{T_{u}}_{\lambda\Psi}
			\Big(\bar{u}^{(\ell)}-\lambda(T_{u})^{-1}
			\big(\nabla E(\bar{u}^{(\ell)})+A^{\mathrm{T}}\bar{v}^{(\ell)}\big)\Big),\\
			v^{(\ell+1)}=\mathrm{Prox}^{T_{v}}_{\sigma\iota^{*}_{\delta}}
			\Big(\bar{v}^{(\ell)}+\sigma(T_{v})^{-1}
			\big(Au^{(\ell+1)}
			-A(T_{u})^{-1}
			\big(\tilde{\nabla}\mathcal{H}(u^{(\ell+1)})
			+A^{\mathrm{T}}\bar{v}^{(\ell)}\big)\big)\Big),\\
			\bar{u}^{(\ell+1)}=(1+\zeta_{1})u^{(\ell+1)}-\zeta_{1}u^{(\ell)},\\
			\bar{v}^{(\ell+1)}=(1+\zeta_{2})v^{(\ell+1)}-\zeta_{2}v^{(\ell)}.
		\end{cases}
	\end{aligned}
\end{eqnarray}
where $\lambda$ and $\sigma$ are the primal and dual step sizes, and $\zeta_{1}$ and $\zeta_{2}$ are extrapolation coefficients. The computable subgradient $\tilde{\nabla}\mathcal{H}$ is obtained from the primal proximal inclusion:
\begin{align}
	\tilde{\nabla}\mathcal{H}(u^{(\ell+1)})
	=\dfrac{T_{u}(\bar{u}^{(\ell)}-u^{(\ell+1)})}{\lambda}
	-A^{\mathrm{T}}\bar{v}^{(\ell)}
	-\nabla E(\bar{u}^{(\ell)})+\nabla E(u^{(\ell+1)}).
\end{align}
Following \cite{Zeng2026variable}, $T_u$ is not taken as the exact Hessian of $\mathcal{H}$. Instead, we construct a regularized objective $\widehat{\mathcal H}_{r}$ by replacing the singular action density $\mathcal G(\rho,\boldsymbol m)$ with $\widehat{\mathcal G}_{r}(\rho,\boldsymbol m)=|\boldsymbol m|^{2}/(\rho+r)$ for a small $r>0$, and take
\begin{align}\label{eq:VPTPD_preconditioner}
	T_{u}=\mathrm{diag}\big(\nabla^{2}\widehat{\mathcal H}_{r}(u^{k})\big),
	\qquad
	T_{v}=A(T_{u})^{-1}A^{\mathrm T},
\end{align}
where $u^{k}$ is the solution from the previous JKO step. This diagonal, PDE-time-dependent choice preserves the separability of the primal proximal step and balances convergence acceleration with the computational cost of each iteration. 

Moreover, the VPTPD method can be implemented with an adaptive step-size strategy that automatically adjusts the primal and dual step sizes based on the current state of the iteration, further enhancing convergence speed. We refer the audience to \cite{Zeng2026variable} for the detailed derivation and implementation of VPTPD and its comparison with other primal--dual-type methods. Here we focus only on the ingredients needed for the PNP-constrained JKO problem.

\subsection{Computation of proximal operators}\label{sec:proximal_operator}

\subsubsection[Proximal operator for primal variables]{Proximal operator for primal variables: $\mathrm{Prox}^{T_{u}}_{\lambda\Psi}$}
Since $p$ and $n$ are independent in the transport action $\Psi(u)$, the primal proximal operator is naturally separated with respect to the two species. More importantly, as long as $T_u$ is diagonal, as in both PrePD and VPTPD, the proximal operator remains pointwise separable and can be computed in parallel over all grid cells. For a local variable $(\rho,\boldsymbol m)$, denote the corresponding local diagonal block of $T_u$ by $T_u^{\mathrm{loc}}=\mathrm{diag}(d_{\rho},d_{m}I)$. Then the componentwise proximal problem for the local action function $\mathcal{G}$ takes the form
\begin{align}\label{eq:local_prox_primal}
	(\rho^{*},\boldsymbol m^{*})
	=\mathrm{Prox}^{T_u^{\mathrm{loc}}}_{\lambda\mathcal{G}}(\hat{\rho},\hat{\boldsymbol m})
	=\mathop{\mathrm{argmin}}\limits_{\rho,\boldsymbol m}
	\left\{
	\dfrac{d_{\rho}}{2}|\rho-\hat{\rho}|^{2}
	+\dfrac{d_{m}}{2}\Vert\boldsymbol m-\hat{\boldsymbol m}\Vert^{2}
	+\lambda\dfrac{\Vert\boldsymbol m\Vert^{2}}{\rho}
	\right\}.
\end{align}
The definition of $\mathcal{G}$ implies that the admissible solution satisfies either $\rho^{*}>0$ or $(\rho^{*},\boldsymbol m^{*})=(0,\boldsymbol0)$. In the positive case, $\rho^{*}$ is the largest positive real root of the scalar cubic polynomial
\begin{align}\label{eq:cubic}
	P(X)=d_{\rho}(X-\hat{\rho})(d_mX+2\lambda)^{2}
	-\lambda d_m^{2}\Vert\hat{\boldsymbol m}\Vert^{2}=0,
\end{align}
and the optimal momentum is given by
\begin{align}
	\boldsymbol m^{*}
	=\dfrac{d_m\rho^{*}}{d_m\rho^{*}+2\lambda}\hat{\boldsymbol m}.
\end{align}
Otherwise, if the cubic equation does not have a positive root, then the optimal solution is $(\rho^{*},\boldsymbol m^{*})=(0,\boldsymbol0)$.
This scalar equation \eqref{eq:cubic} can be solved either by the closed-form formula for cubic equations or by the Newton method with tailored initial guesses developed in \cite{Zeng2026variable}. Thus the primal proximal step is inexpensive, pointwise, and well suited for parallel implementation; see also \cite{Carrillo2024StructurePD,Zeng2026variable} for detailed derivations.

\subsubsection[Proximal operator for dual variables]{Proximal operator for dual variables: $\mathrm{Prox}^{T_{v}}_{\lambda^{-1}\iota^{*}_{\delta}}$}
The dual proximal operator can be written through the generalized Moreau identity:
\begin{align}
	\mathrm{Prox}^{T_{v}}_{\lambda^{-1}\iota^{*}_{\delta}}(y)
	=y-\lambda^{-1}T^{-1}_{v}\mathrm{Prox}_{\iota_{\delta}}^{T_{v}^{-1}}\big(\lambda T_{v}y\big).
\end{align}
For the PrePD update, we set $\delta>0$ to keep the diagonal dominance of the saddle-point system \eqref{eq:PDflow}. The corresponding dual proximal step reduces to a classical trust-region subproblem that can be solved by standard exact solvers or approximately by a projection onto the ball with very small $\delta$ \cite{Wu2026PDPFSMCL}: 
\begin{eqnarray}
	\begin{aligned}
		\mathrm{Prox}^{T_{v}}_{\lambda^{-1}\iota^{*}_{\delta}}(y)&\approx y-\lambda^{-1}T^{-1}_{v}\left\{
		\begin{aligned}
			&\lambda T_{v}y\quad&&\text{$\Vert \lambda T_{v}y-b\Vert_{2}<\delta$},\\
			&\delta\dfrac{\lambda T_{v}y-b}{\Vert \lambda T_{v}y-b\Vert_{2}}+b\quad&&\text{otherwise}.
		\end{aligned}
		\right.
	\end{aligned}
\end{eqnarray}
The approximation above is used as an inexact subproblem solver in the overall primal--dual iteration, which does not affect overall convergence in practice \cite{Wu2026PDPFSMCL,Liu2021AccelerationPrimalDual}. To further avoid redundant matrix inversions in \eqref{eq:PrePD}, one practical realization introduces auxiliary variables
$z^{(\ell)}=\bar{v}^{(\ell)}+r^{(\ell)}$, where
$\bar{v}^{(\ell)}=\lambda T_{v}v^{(\ell)}$ and
$r^{(\ell)}=A\bar{u}^{(\ell)}-b$. The inexact dual proximal step can then be written as
\begin{eqnarray}\label{eq:proxideltaybar}
	\begin{aligned}
		v^{(\ell+1)}=\lambda^{-1}T_{v}^{-1}\bar{v}^{(\ell+1)},\quad\text{where}\quad
		\bar{v}^{(\ell+1)}&=\left\{
		\begin{aligned}
			&\textbf{0}\quad&&\text{$\Vert z^{(\ell)}\Vert_{2}<\delta$},\\
			&\Big(1-\dfrac{\delta}{\Vert z^{(\ell)}\Vert_{2}}\Big)z^{(\ell)} \quad&&\text{otherwise}.
		\end{aligned}
		\right.
	\end{aligned}
\end{eqnarray}

In the VPTPD implementation, we can freely set $\delta=0$ due to the presence of the Schur complement in the transformed saddle-point system \eqref{eq:VPTPDflow}. In this case, $\iota_{0}^{*}(v)=\langle b,v\rangle$, and the dual proximal step reduces to
\begin{align}
	\mathrm{Prox}^{T_{v}}_{\sigma\iota^{*}_{0}}(y)=y-\sigma T_v^{-1}b.
\end{align}
Therefore, the transformed dual update in VPTPD is explicitly given by
\begin{align}\label{eq:VPTPD_dual_delta0}
	v^{(\ell+1)}
	=\bar{v}^{(\ell)}+\sigma(T_{v})^{-1}
	\Big(Au^{(\ell+1)}-b
	-A(T_{u})^{-1}
	\big(\tilde{\nabla}\mathcal{H}(u^{(\ell+1)})
	+A^{\mathrm T}\bar{v}^{(\ell)}\big)\Big).
\end{align}

For scalar Wasserstein gradient flows, the dual proximal step reduces to the inversion of a transport-related elliptic operator, which can be efficiently computed by FFT-based fast algorithms or preconditioned iterative methods. In PNP models, however, the action of $T_v^{-1}$ involves the coupling of the Poisson constraint with the two ionic continuity equations, which becomes the dominant computational bottleneck. Therefore, the decisive question is whether the dual proximal operator can be evaluated rapidly under different potential boundary conditions. This motivates the fast solvers developed in Section~\ref{sec:5}, where we exploit the block structure induced by the PNP constraints and discuss efficient algorithms for both PrePD and VPTPD under Dirichlet and Neumann boundary conditions.

\section{Fast solvers for the coupled dual subproblem}\label{sec:5}
For the Poisson-constrained JKO scheme of PNP models, the coupling between the ionic concentrations and the electrostatic potential induces the following block structure in the linear constraint operator. Considering the one-dimensional case for illustration purpose, the discrete constraint operator $A$ can be written as:
\begin{align}\label{eq:constraintA}
	A=\begin{bmatrix}
		B & 0 & 0\\
		0 & B & 0\\
		-R & R & \epsilon L_{\phi}
	\end{bmatrix},
	\quad
	B=\begin{bmatrix} I & \nabla_h\cdot \end{bmatrix},
	\quad
	R=\begin{bmatrix} I & 0 \end{bmatrix},
\end{align}
where $B$ is the discrete one-step temporal-spatial divergence operator for the continuity equation \eqref{eq:fullJKO}, $R$ is the restriction operator that extracts the density component, and $L_{\phi}$ is the discrete Laplace operator for the electrostatic potential with appropriate boundary conditions. 

For a general diagonal preconditioner for the primal update:
\begin{align}
	T_u=\mathrm{diag}(\Lambda_{p},\Lambda_{n},\Lambda_{\phi}),
\end{align}
with $T_u=I$ for PrePD and $T_u=\mathrm{diag}(\nabla^{2}\widehat{\mathcal H}_{r}(u^{k}))$ for VPTPD \eqref{eq:VPTPD_preconditioner}, the corresponding dual preconditioner is given by
\begin{align}\label{eq:block}
	T_{v}=AT_u^{-1}A^{\mathrm{T}}
	=
	\begin{bmatrix}
		M_{p} & 0 & -Q_{p}\\
		0 & M_{n} & Q_{n}\\
		-Q_{p}^{\mathrm{T}} & Q_{n}^{\mathrm{T}} & P
	\end{bmatrix},
	\quad \text{where }
	\begin{aligned}
		&M_{p}=B\Lambda_{p}^{-1}B^{\mathrm{T}},\quad
		M_{n}=B\Lambda_{n}^{-1}B^{\mathrm{T}},\\
		&Q_{p}=B\Lambda_{p}^{-1}R^{\mathrm{T}},\quad
		Q_{n}=B\Lambda_{n}^{-1}R^{\mathrm{T}},\\
		&P=R\Lambda_{p}^{-1}R^{\mathrm{T}}+R\Lambda_{n}^{-1}R^{\mathrm{T}}
		+\epsilon^{2}L_{\phi}\Lambda_{\phi}^{-1}L_{\phi}^{\mathrm{T}}.
	\end{aligned}
\end{align}
The dual subproblem for both PrePD \eqref{eq:proxideltaybar} and VPTPD \eqref{eq:VPTPD_dual_delta0} reduces to solve the following linear system:
\begin{align}\label{eq:dual_block_system}
	\begin{bmatrix}
		M_{p} & 0 & -Q_{p}\\
		0 & M_{n} & Q_{n}\\
		-Q_{p}^{\mathrm{T}} & Q_{n}^{\mathrm{T}} & P
	\end{bmatrix}
	\begin{bmatrix}
		v_{p}\\
		v_{n}\\
		v_{\phi}
	\end{bmatrix}
	=
	\begin{bmatrix}
		f_{p}\\
		f_{n}\\
		f_{\phi}
	\end{bmatrix}.
\end{align}

Owing to the large size of the above system, direct inversion of the full coefficient matrix can be computationally expensive. We therefore exploit the specific structures induced by PrePD and VPTPD to construct their efficient fast solvers. In the following, we first present two general iterative strategies for solving the block system \eqref{eq:dual_block_system}, and then discuss their realizations under different primal--dual frameworks for different boundary conditions.

\subsection{Block solvers for the coupled dual system}\label{sec:dual_solver}
Before specifying the realizations under different primal--dual frameworks, we first present two general iterative strategies for solving the block system \eqref{eq:dual_block_system}, namely the block Gauss--Seidel method (BGS) and the Schur-complement-transformed preconditioned conjugate gradient method (Schur-PCG).

\subsubsection{Block Gauss-Seidel solver (BGS)}
A natural approach is to apply a block Gauss--Seidel iteration directly to \eqref{eq:dual_block_system} to decouple the variables. Given the current iterate $v_{\phi}^{m}$ in the $m$-th inner iteration, we successively compute
\begin{align}
	v_{p}^{m+1}=M_{p}^{-1}\bigl(f_{p}+Q_{p}v_{\phi}^{m}\bigr),\quad
	v_{n}^{m+1}=M_{n}^{-1}\bigl(f_{n}-Q_{n}v_{\phi}^{m}\bigr),
\end{align}
and then update $v_{\phi}$ by
\begin{align}
	v_{\phi}^{m+1}=P^{-1}\bigl(f_{\phi}+Q_{p}^{\mathrm{T}}v_{p}^{m+1}-Q_{n}^{\mathrm{T}}v_{n}^{m+1}\bigr).
\end{align}
The iteration is terminated when
\begin{align}\label{eq:GS_Stop}
	\dfrac{\Vert v^{m+1}_{\phi}-v^{m}_{\phi}\Vert_{2}}{\max\left\{1,\Vert v^{m}_{\phi}\Vert_{2}\right\}}\leq\epsilon_{tol}.
\end{align}
Once \eqref{eq:GS_Stop} is satisfied, we set $(v_{p},v_{n},v_{\phi})$ as the updated dual variable. 

The BGS method is simple to implement and can be efficient if the inversion of $M_{p}$, $M_{n}$ and $P$ is cheap. In particular, for the PrePD choice $T_u=I$ (and hence $T_v=AA^{\mathrm T}$), the matrices $M_{p}=M_{n}$ and $P$ are all discrete elliptic operators that can be inverted efficiently by FFT-based fast solvers. However, it may converge slowly for large-scale problems, especially when the coupling between the variables is strong. Therefore, we also consider an alternative based on the preconditioned conjugate gradient method applied to the Schur-complement reduced system.

\subsubsection{Schur-reduced PCG solver (Schur-PCG)}
We consider an alternative approach that first reduces the coupled block system \eqref{eq:dual_block_system} to a Schur-complement equation for the electrostatic dual variable and then solves the reduced system by preconditioned conjugate gradients, which we refer to as the Schur-reduced PCG (Schur-PCG) method. Specifically, from the first two equations in \eqref{eq:dual_block_system}, we obtain
\begin{align}\label{eq:elim_ypyn_new}
	v_{p}=M_{p}^{-1}(f_{p}+Q_{p}v_{\phi}),\quad
	v_{n}=M_{n}^{-1}(f_{n}-Q_{n}v_{\phi}).
\end{align}
Substituting \eqref{eq:elim_ypyn_new} into the third equation of \eqref{eq:dual_block_system} yields
\begin{align}\label{eq:schur_eq_new}
	\bigl(P-Q_{p}^{\mathrm{T}}M_{p}^{-1}Q_{p}-Q_{n}^{\mathrm{T}}M_{n}^{-1}Q_{n}\bigr)v_{\phi}
	=f_{\phi}+Q_{p}^{\mathrm{T}}M_{p}^{-1}f_{p}-Q_{n}^{\mathrm{T}}M_{n}^{-1}f_{n}.
\end{align}
We solve the above reduced system $Sv_{\phi}=r_{\phi}$ for $v_{\phi}$ by PCG iteration until the stopping criterion \eqref{eq:GS_Stop} is satisfied, where we define the Schur operator
\begin{align}\label{eq:schur_def_new}
	S:=P-Q_{p}^{\mathrm{T}}M_{p}^{-1}Q_{p}-Q_{n}^{\mathrm{T}}M_{n}^{-1}Q_{n},\quad
	r_{\phi}:=f_{\phi}+Q_{p}^{\mathrm{T}}M_{p}^{-1}f_{p}-Q_{n}^{\mathrm{T}}M_{n}^{-1}f_{n},
\end{align}
and then recover $(v_{p},v_{n})$ by back substitution through \eqref{eq:elim_ypyn_new}. 

The efficiency of Schur-PCG stems from two main aspects. First, compared with BGS applied to the original coupled dual system, the Schur reduction removes the explicit block coupling and allows the reduced system to be solved by PCG, which typically provides a more robust Krylov acceleration than a stationary block iteration, especially for strongly coupled or ill-conditioned systems. Second, the Schur operator is applied in a matrix-free manner: each PCG iteration only requires matrix--vector products with the Poisson-related block ($L_{\phi}$) and inverse actions of the transport-related blocks ($M_p^{-1}$ and $M_n^{-1}$) and the shifted biharmonic block ($P$). These operations are efficient because they involve the inversion of discrete elliptic or Laplace-type operators, which can be accelerated by FFT-based solvers for discrete Laplacian in PrePD, or by sparse Cholesky factorizations, PCG, and multigrid-type solvers for more general structured matrices in VPTPD \cite{Zeng2026variable}. The same fast inverse actions are reused throughout the primal--dual iterations for one-step JKO, leading to an efficient and memory-friendly solver for the Poisson-constrained JKO dual system.

\subsection{FFT-based dual solvers for PrePD}\label{sec:prepd_dual_solver}
For PrePD, $T_u=I$, or equivalently $\Lambda_p=\Lambda_n=\Lambda_\phi=I$ in \eqref{eq:block}, and hence $T_v=AA^{\mathrm T}$ has the blocks
\begin{align}
	M_p=M_n=M:=I+ L_{\mathrm{cont}},\qquad
	Q_p=Q_n=I,\qquad
	P=2I+\epsilon^2L_{\phi}L_{\phi}^{\mathrm T},
\end{align}
where $L_{\mathrm{cont}}=-\nabla_h \cdot (\nabla_h) $ is the discrete (negative) Laplace operator associated with the continuity equations with no-flux boundary conditions, and $L_\phi$ is the discrete (negative) Laplace operator associated with the Poisson equation with the specified boundary conditions. 


On a uniform rectangular grid, $L_{\mathrm{cont}}$, $L_{\phi,N}$ (for Neumann boundary conditions), and $L_{\phi,D}$ (for Dirichlet boundary conditions) have tensor-product structures and are diagonalized by discrete cosine or sine transforms. Consequently, $M$ and $P$ can be inverted by transform-space division (with special care for homogeneous Neumann boundary conditions where the corresponding operator has a null mode). Table~\ref{tab:transforms} lists the one-dimensional transforms and eigenvalues. In multiple dimensions, the eigenvalue associated with a tensor-product mode is the sum of the corresponding one-dimensional eigenvalues. 

\begin{table}[htbp]
	\centering
	\small
	\caption{One-dimensional transforms and eigenvalues for the PrePD block solvers.}
	\begin{tabular}{ccccc}
		\toprule
		Operator & Boundary conditions & Forward & Inverse & Eigenvalue ($i=0,\ldots,N_x-1$)\\
		\midrule
		$L_{\text{cont}}$ & no-flux  & DCT-II & $\frac{1}{2N_x}$ DCT-III & $\mu^i_{\text{cont}}=\frac{1}{2\Delta x^2}\bigl(1-\cos(\frac{2\pi i}{N_x})\bigr)$\\
		$L_{\phi,N}$ & Neumann & DCT-II & $\frac{1}{2N_x}$ DCT-III & $\mu^i_{\phi,N}=\frac{4}{\Delta x^2}\sin^2\bigl(\frac{\pi i}{2N_x}\bigr)$\\
		$L_{\phi,D}$ & Dirichlet & DST-II & $\frac{1}{2N_x}$ DST-III & $\mu^i_{\phi,D}=\frac{4}{\Delta x^2}\sin^2\bigl(\frac{\pi(i+1)}{2N_x}\bigr)$\\
		\bottomrule
	\end{tabular}
	\label{tab:transforms}
\end{table}

\subsubsection{Dirichlet boundary conditions: FFT-based BGS and Schur-PCG solver}
For Dirichlet boundary conditions for $\phi$, $M$ is diagonalized via a DCT, while $P$ is diagonalized by a DST, so the three diagonal blocks in \eqref{eq:dual_block_system} can not be diagonalized by a single transform. Instead, the BGS iteration implements individual fast actions of $M^{-1}$ (by DCT-based fast algorithms) and $P^{-1}$ (by DST-based fast algorithms) without assembling the coupled matrix. Similarly, the Schur-PCG iteration, which involves the matrix--vector multiplication of $S=P-2M^{-1}$ and the computation of $r_\phi=f_\phi+M^{-1}(f_p-f_n)$ and $(v_p,v_n)$, can also be efficiently implemented with the invertible block $P$ as the preconditioner. One BGS sweep or one Schur-PCG iteration therefore costs $\mathcal O(N\log N)$ ($N$ as the total number of grid cells) and does not require storage of assembling the full block matrix.

\subsubsection{Neumann boundary conditions: DCT-based direct solver (DCT-DS)}\label{sec:DCT-DS}
For Neumann boundary conditions, both $L_{\mathrm{cont}}$ and $L_{\phi,N}$ are diagonalized by the same forward-inverse transform pair (Table~\ref{tab:transforms}). Then the full block system \eqref{eq:dual_block_system} after transformation can be decomposed mode by mode into $N$ independent 3-by-3 systems. Specifically, for each mode $i$, the transformed 3-by-3 subsystem is
\begin{align}\label{eq:modewise_new}
	\begin{bmatrix}
		a_i & 0 & -1\\
		0 & a_i & 1\\
		-1 & 1 & c_i
	\end{bmatrix}
	\begin{bmatrix}
		\hat v^i_p\\
		\hat v^i_n\\
		\hat v^i_\phi
	\end{bmatrix}
	=
	\begin{bmatrix}
		\hat f^i_p\\
		\hat f^i_n\\
		\hat f^i_\phi
	\end{bmatrix},
\end{align}
where the components (for 1D case) are given by
\begin{align}
	a_i:=1+\mu^i_{\mathrm{cont}},\qquad
	c_i:=2+\epsilon^2(\mu^i_{\phi,N})^2.
\end{align}
For every nonzero mode $i=1,\ldots,N_x-1$, the subsystem can be explicitly solved by:
\begin{align}\label{eq:neu_phi_explicit_new}
	\hat v^i_p=\dfrac{\hat f^i_p+\hat v^i_\phi}{a_i}, \quad
	\hat v^i_n=\dfrac{\hat f^i_n-\hat v^i_\phi}{a_i}, \quad \hat v^i_\phi=
	\dfrac{\hat f^i_\phi+(\hat f^i_p-\hat f^i_n)/a_i}{c_i-2/a_i}.
\end{align}
In particular, the zero-frequency block (for $i=0$) is singular since $a_0=1$ and $c_0=2$, and requires special treatment. Solvability therefore requires the compatibility condition $\hat f^0_\phi=-\hat f^0_p+\hat f^0_n$, which is equivalent to the zero-mode compatibility condition for the Poisson equation in the primal problem. When this condition holds, the dual solution is determined up to a multiple of $(1,-1,1)^{\mathrm T}$. We select a unique representative by imposing the gauge condition for the electrostatic potential, which yields:
\begin{align}\label{eq:neu_gauge_new}
	\hat v^0_\phi=0,\quad \text{and hence } \hat v^0_p=\hat f^0_p,\quad \hat v^0_n=\hat f^0_n.
\end{align}
The dual variables $(v_p,v_n,v_\phi)$ are then recovered by the inverse DCT from $(\hat v_p,\hat v_n,\hat v_\phi)$.

\begin{myremark}
	BGS and Schur-PCG solve the same coupled dual system from two different perspectives. BGS applies a stationary block iteration to the original system, whereas Schur-PCG eliminates the ionic dual variables and applies Krylov acceleration to the reduced electrostatic system. Thus, BGS has a lower cost per inner iteration, while Schur-PCG is generally more robust in strongly coupled or ill-conditioned regimes (see the case $\epsilon=0.06$ in Table~\ref{tal:exacteps}).
\end{myremark}

\begin{myremark}\label{rmk:inexact}
	In practice, the dual system need not be solved exactly at every primal–dual iteration. With the solution from the previous JKO step as a warm start, one BGS sweep or one PCG iteration is often sufficient for the overall convergence of the primal-dual method \cite{Liu2021AccelerationPrimalDual}. Although such inexact solves may slightly increase the number of outer iterations, they substantially reduce the overall computational cost (see Fig.~\ref{fig:BGSSchurPCG} and Table~\ref{tal:Inexacteps}).
\end{myremark}

\subsection{Sparse block dual solver for VPTPD}\label{sec:vptpd_dual_solver}
For VPTPD, since $T_u$ is a diagonal matrix \eqref{eq:VPTPD_preconditioner} rather than an identity matrix for PrePD, the blocks $M_p$, $M_n$ and $P$ in \eqref{eq:block} cannot be diagonalized by the DCT/DST pairs. Instead, we exploit sparse Cholesky factorization with approximate minimum degree ordering for moderate sizes, or via incomplete Cholesky preconditioned conjugate gradient method for large-size systems \cite{Zeng2026variable} to efficiently realize the inversion of $M_p$, $M_n$, and $P$ in the BGS or Schur-PCG iteration. Given that the tailored preconditioner $T_u$ in \eqref{eq:VPTPD_preconditioner} and the corresponding block matrices ($M_p$, $M_n$, $P$) are fixed during the primal--dual iterations of a given JKO step, their factorizations can be reused in all inner BGS/Schur-PCG iterations, which amortizes the setup cost.

\section{Numerical results}\label{sec:6}
In this section, we validate the convergence of the proposed Poisson-constrained JKO scheme, assess the performance of the primal--dual splitting methods for PNP models, and investigate ionic interaction phenomena through a series of numerical experiments. In Sec.~\ref{sec:6.1}, we verify the accuracy and structure-preserving properties of the fully discrete JKO scheme, compare it with existing numerical and optimization methods, and assess the proposed primal--dual algorithms and dual solvers. In Sec.~\ref{sec:6.2}, we present extended experiments for modified PNP models, illustrating the influence of concentration-gradient energy and spatial ionic interactions.

Unless otherwise specified, the following parameters are used in the numerical experiments:
\begin{eqnarray}
	\begin{aligned}
		&D_{p}=D_{n}=1,&\quad& z_{p}=1,&\quad& z_{n}=-1,&\quad&
		\phi_{e}=0,\quad &e=1,&\quad&\delta=10^{-7},&\quad&\epsilon_{tol}=10^{-5}.
	\end{aligned}
\end{eqnarray}

\subsection{Validation tests}\label{sec:6.1}
\subsubsection{Benchmark experiments for accuracy}\label{sec:6.1.1}
We first consider one-dimensional two-species classical PNP model within the domain $\left[-1,1\right]$
\begin{eqnarray}\label{eq:basicPNP}
	\begin{aligned}
		\begin{cases}
			\partial_{t} p=\partial_x\Big(p\,\partial_x(\log p+\phi)\Big),\\
			\partial_{t} n=\partial_x\Big(n\,\partial_x(\log n-\phi)\Big),\\
			-\epsilon\partial_{xx}\phi=p-n+\psi^{0},
		\end{cases}
	\end{aligned}
\end{eqnarray}
where we impose Dirichlet boundary conditions on the electrostatic potential and no-flux boundary conditions on the ionic concentrations. The initial data and boundary conditions are chosen as:
\begin{align}\label{eq:1Dinitial}
	p^{0}=2-x^{2},\quad n^{0}=2+\sin(\pi x),\quad\psi^{0}=0,\quad\phi(-1,t)=-1,\quad\phi(1,t)=1,\quad\epsilon=1.
\end{align}

We solve the Poisson-constrained JKO scheme \eqref{eq:fullJKO} for the above equation by the PrePD method \eqref{eq:PrePD}, and compare the numerical solution at $t=0.1$ with the reference solution, which is computed with a very small time step $\tau=10^{-5}$ and fine spatial discretization $\Delta x=0.002$. Table~\ref{tab:FirstTimeorder} shows the first-order temporal accuracy of the Poisson-constrained JKO scheme.

\begin{table}[htbp]
	\centering
	\caption{1D test \eqref{eq:basicPNP} and \eqref{eq:1Dinitial} with Dirichlet boundary conditions: Temporal convergence of $p$, $n$, and $\phi$ at $t=0.1$.}
	\label{tab:FirstTimeorder}
	\begin{tabular}{ccccccc}
		\toprule
		\multirow{2}{*}{$\tau$}
		& \multicolumn{2}{c}{$p$}
		& \multicolumn{2}{c}{$n$}
		& \multicolumn{2}{c}{$\phi$}\\
		\cmidrule(lr){2-3}\cmidrule(lr){4-5}\cmidrule(lr){6-7}
		& $\Vert p-p_{\text{ref}}\Vert_{\infty}$ & $\text{Rate}$
		& $\Vert n-n_{\text{ref}}\Vert_{\infty}$ & $\text{Rate}$
		& $\Vert\phi-\phi_{\text{ref}}\Vert_{\infty}$ & $\text{Rate}$\\
		\midrule
		$1/50$  & 3.85E-1 & -    & 3.83E-1 & -    & 1.77E-1 & -\\
		$1/100$ & 1.94E-1 & 0.99 & 1.95E-1 & 0.97 & 8.99E-2 & 0.98\\
		$1/200$ & 9.08E-2 & 1.10 & 9.09E-2 & 1.10 & 4.46E-2 & 1.01\\
		$1/400$ & 4.17E-2 & 1.12 & 4.21E-2 & 1.11 & 2.11E-2 & 1.08\\
		$1/800$ & 1.83E-2 & 1.19 & 1.83E-2 & 1.20 & 9.21E-3 & 1.19\\
		\bottomrule
	\end{tabular}
\end{table}

\begin{figure}[htbp]
	\centering
	\subfigure[Total energy]{
		\begin{minipage}{0.31\textwidth}
			\includegraphics[width=\linewidth]{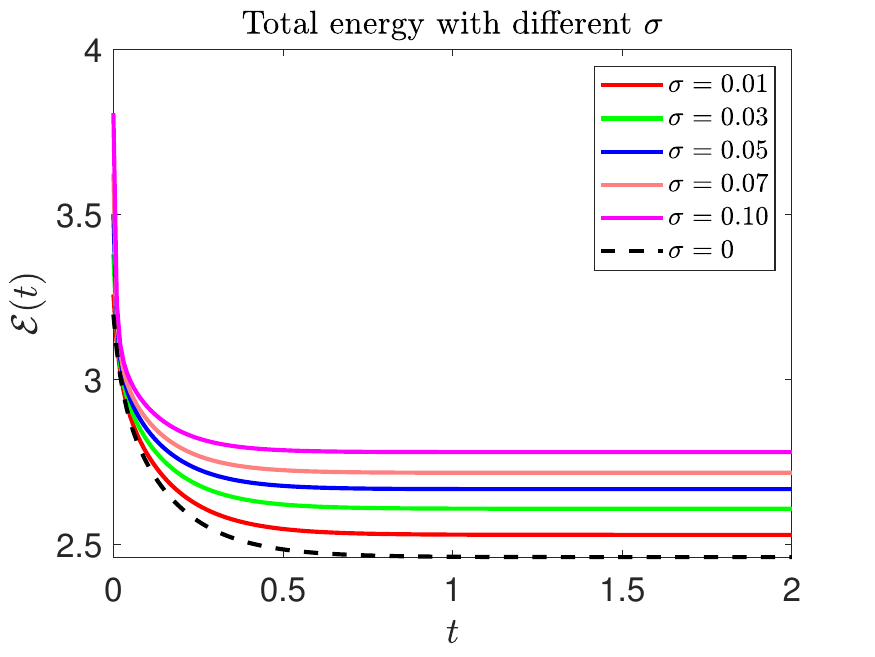}
		\end{minipage}
	}
	\subfigure[Relative mass errors]{
		\begin{minipage}{0.31\textwidth}
			\includegraphics[width=\linewidth]{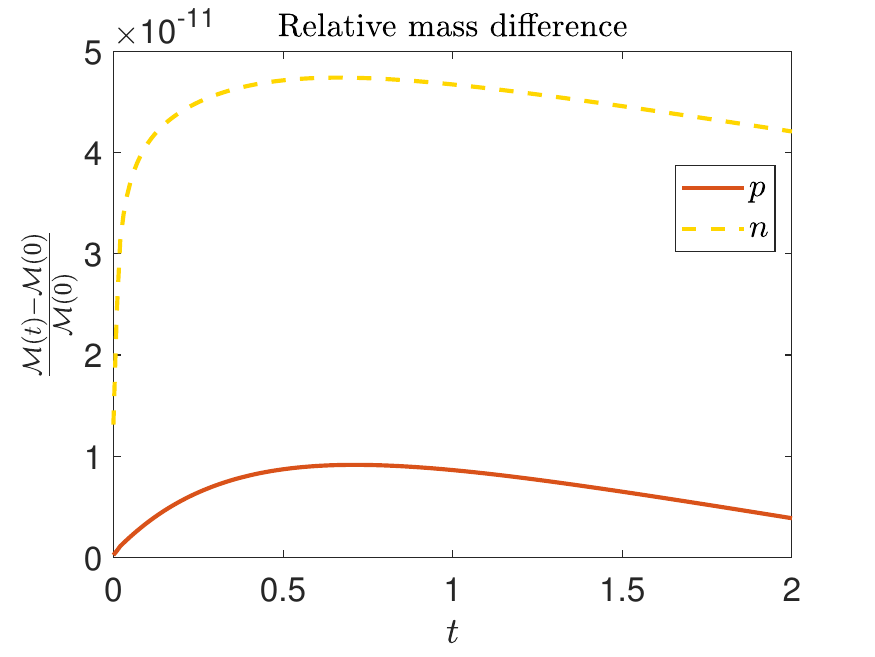}
		\end{minipage}
	}
	\subfigure[Ion concentration extrema]{
		\begin{minipage}{0.31\textwidth}
			\includegraphics[width=\linewidth]{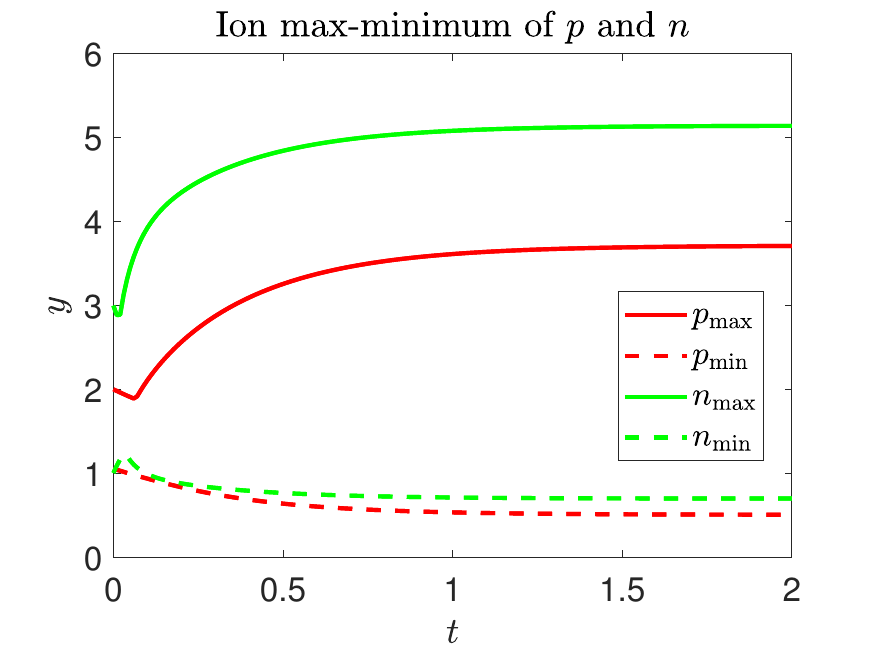}
		\end{minipage}
	}
	\caption{Time evolution of the discrete total energy, relative mass errors, and ionic concentration extrema for 1D test \eqref{eq:basicPNP} and \eqref{eq:1Dinitial} with Dirichlet boundary conditions.}
	\label{fig:EneMassValue}
\end{figure}

\begin{figure}[htbp]
	\centering
	\subfigure{
		\begin{minipage}{0.31\textwidth}
			\includegraphics[width=\linewidth]{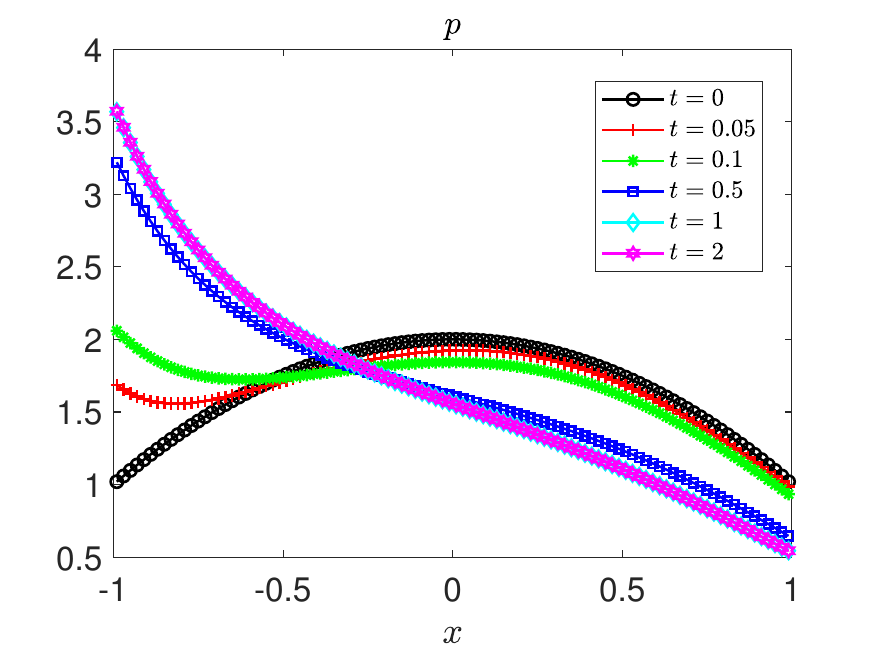}
			\includegraphics[width=\linewidth]{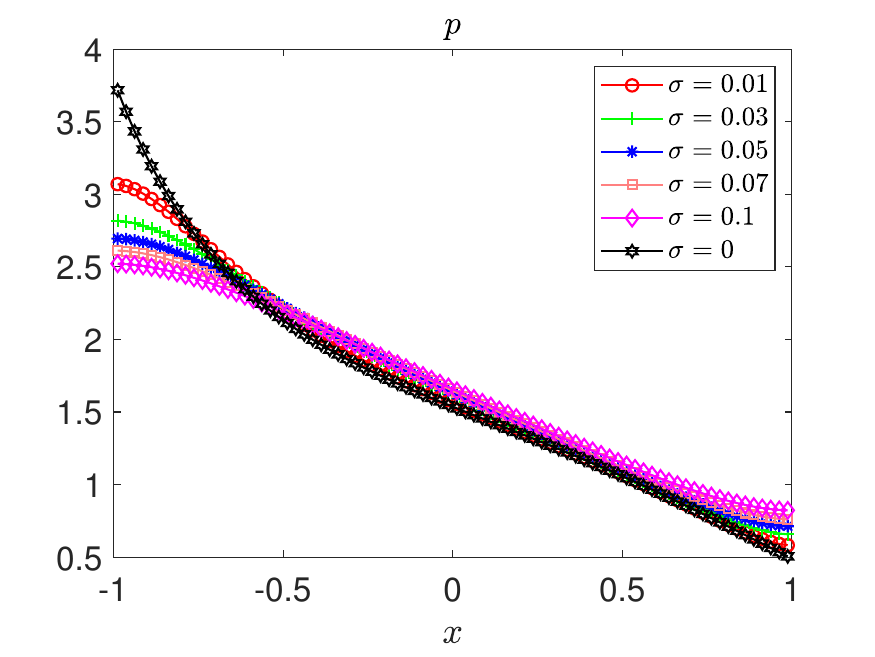}
		\end{minipage}
	}
	\subfigure{
		\begin{minipage}{0.31\textwidth}
			\includegraphics[width=\linewidth]{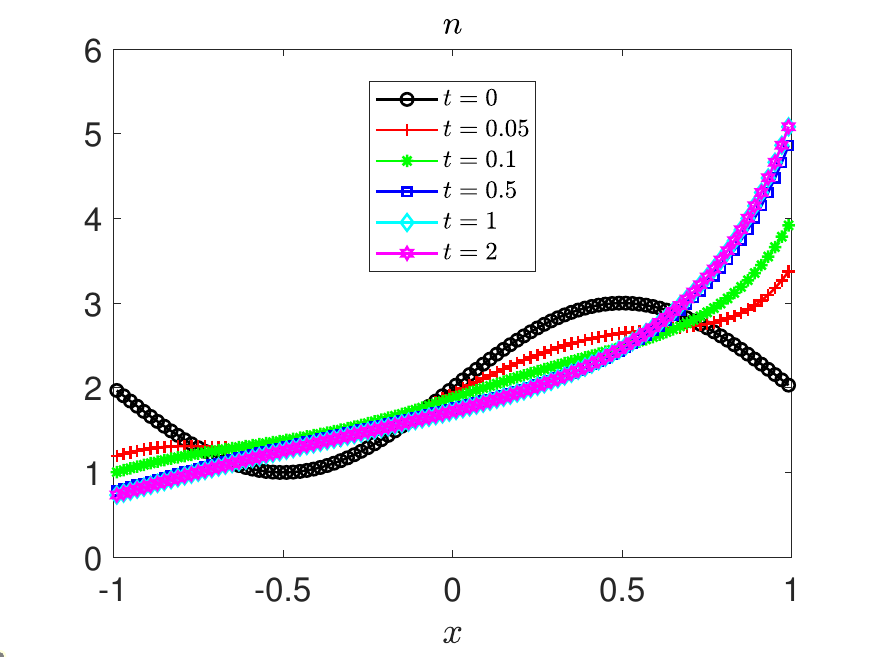}
			\includegraphics[width=\linewidth]{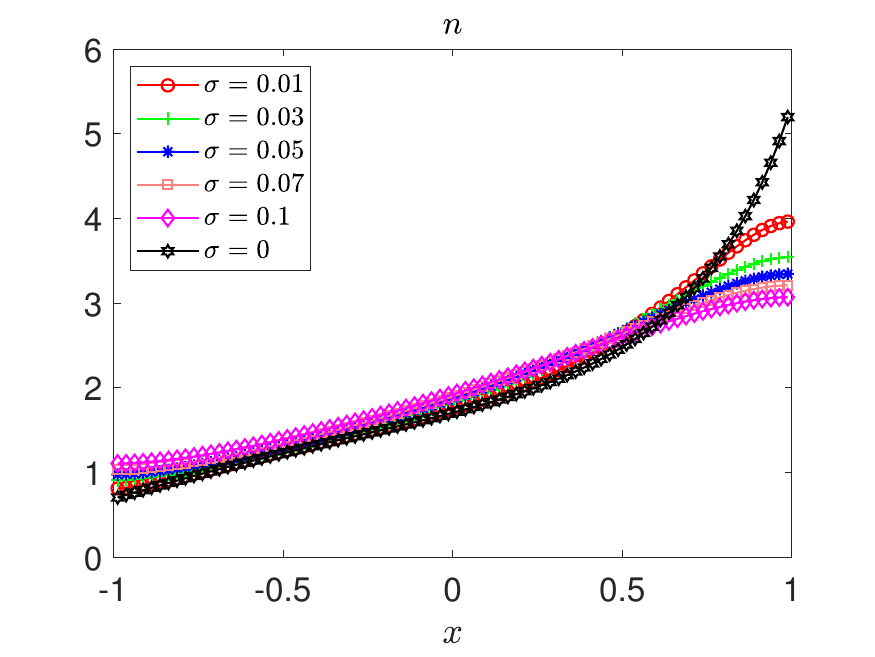}
		\end{minipage}
	}
	\subfigure{
		\begin{minipage}{0.31\textwidth}
			\includegraphics[width=\linewidth]{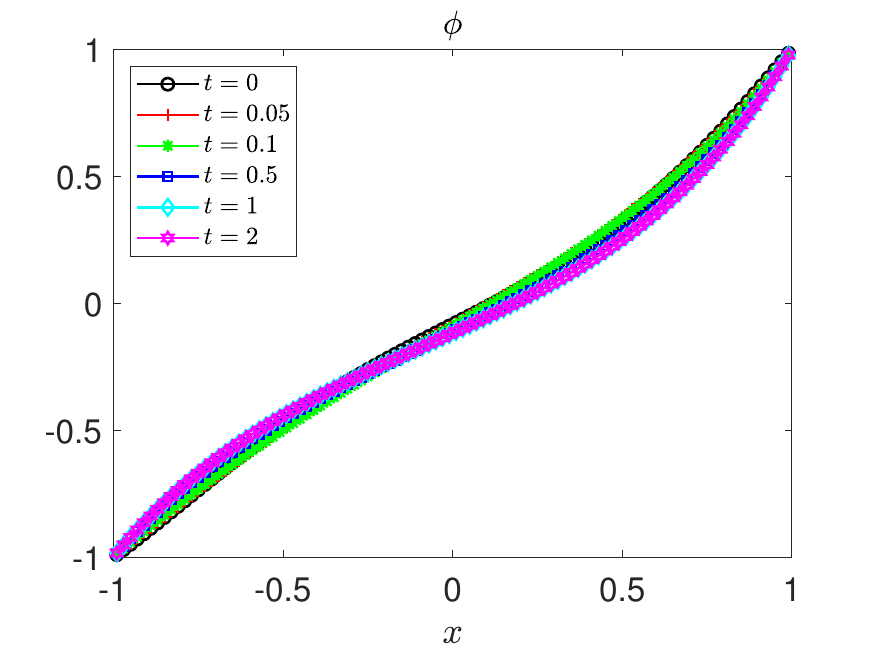}
			\includegraphics[width=\linewidth]{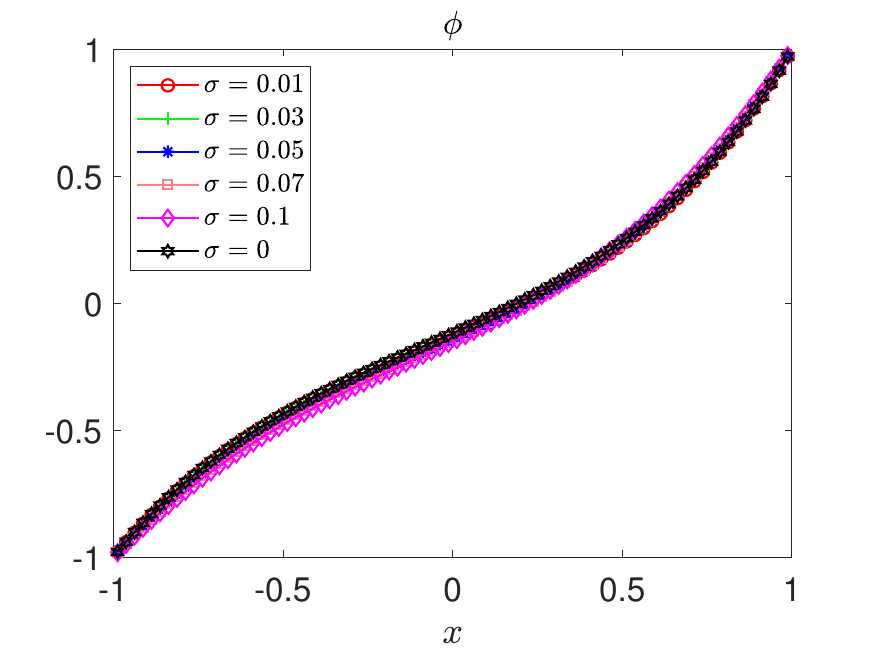}
		\end{minipage}
	}
	\caption{1D test \eqref{eq:basicPNP} and \eqref{eq:1Dinitial} with Dirichlet boundary conditions. Top: Time evolution of $p$, $n$, and $\phi$ for the classical PNP model. Bottom: Final-time solutions of the modified PNP model \eqref{eq:modifiedPNPsecond} for different values of $\sigma$.}
	\label{fig:timepnpsigam}
\end{figure}

We also compute the modified PNP model \eqref{eq:modifiedPNPsecond} with different values of strength of the concentration-gradient correction $\sigma$. Fig.~\ref{fig:EneMassValue} shows the monotone decay of the total energy (a), the relative mass errors (b), and the ionic concentration extrema (c), confirming the energy-dissipation, mass conservation and positivity-preserving properties of the proposed scheme. 
Fig.~\ref{fig:timepnpsigam} (top) presents the evolution of the ionic concentrations ($p,n$) and the electrostatic potential ($\phi$) for the classical PNP model. For the modified PNP model, we observe that the concentration-gradient energy term penalizes sharp spatial variations in the ionic concentrations, leading to smoother profiles; see Fig.~\ref{fig:timepnpsigam} (bottom).

\subsubsection{Diffuse-charge dynamics with small permittivity}

We further assess the proposed Poisson-constrained JKO scheme in the small-permittivity regime using the canonical diffuse-charge relaxation problem of Bazant, Thornton, and Ajdari \cite{Bazant2004Diffuse}, which described the response of a symmetric binary electrolyte between two parallel blocking electrodes to applied voltage. In the thin-double-layer regimes, the system reduces to an one-dimensional PNP model with uniform electroneutral initial state $p(x,0)=n(x,0)=1$ within the interval $\Omega=(-1,1)$. The dimensionless Poisson equation is
\begin{align}
	-\epsilon_{D}^{2}\partial_{xx}\phi
	=\dfrac{p-n}{2},
\end{align}
where $\epsilon_{D}$ denotes the dimensionless Debye length to and it is related to the dimensionless permittivity parameter $\epsilon$ in \eqref{eq:basicPNP} by $\epsilon=2\epsilon_{D}^{2}$. In particular, we consider the initial condition $\phi(x,0)=vx$ and the Dirichlet boundary conditions for the electrostatic potential $\phi$ at the electrodes:
\begin{align}
	\phi(-1,t)=-v,\qquad \phi(1,t)=v,
\end{align}
where $v$ is the applied voltage. To compare the numerical charge with the leading-order linear resistance--capacitance (RC) prediction, we consider the cathodic half-cell charge defined by
\begin{align}
	q_h(t)=\dfrac{1}{2}\int_{-1}^{0}\bigl(p(x,t)-n(x,t)\bigr)\,\mathrm{d}x.
\end{align}
In the thin-double-layer limit, i.e., $\epsilon_{D}\ll 1$, the leading-order  asymptotic approximation was obtained \cite{Bazant2004Diffuse}
\begin{equation}
    q_h/\epsilon_{D}\sim v(1-\exp(-s)),
\end{equation}
where $s=t/\epsilon_{D}$ is the RC charging time scale.

\begin{figure}[htbp]
	\centering
	\subfigure{
		\begin{minipage}{0.31\textwidth}
			\includegraphics[width=\linewidth]{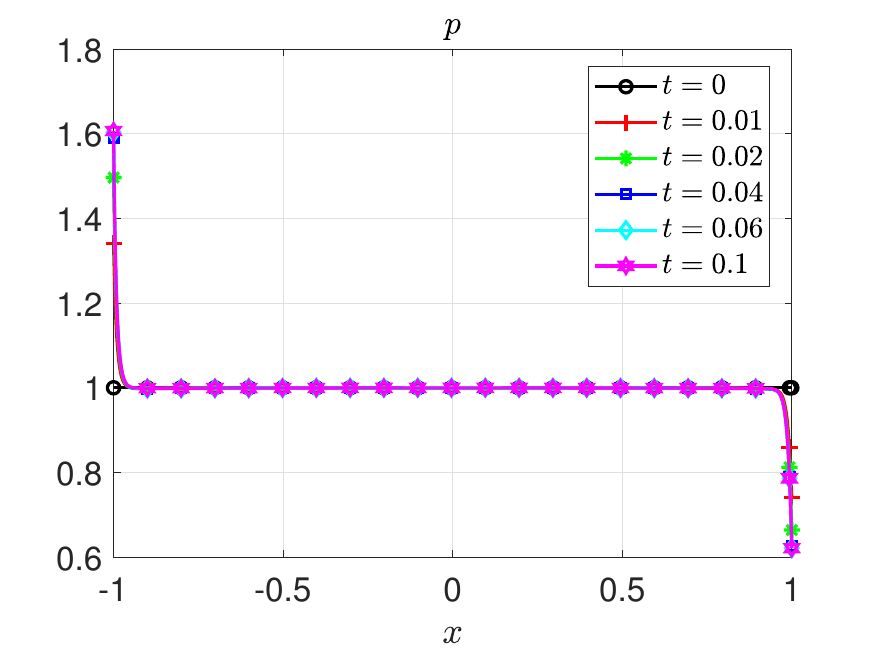}
			\includegraphics[width=\linewidth]{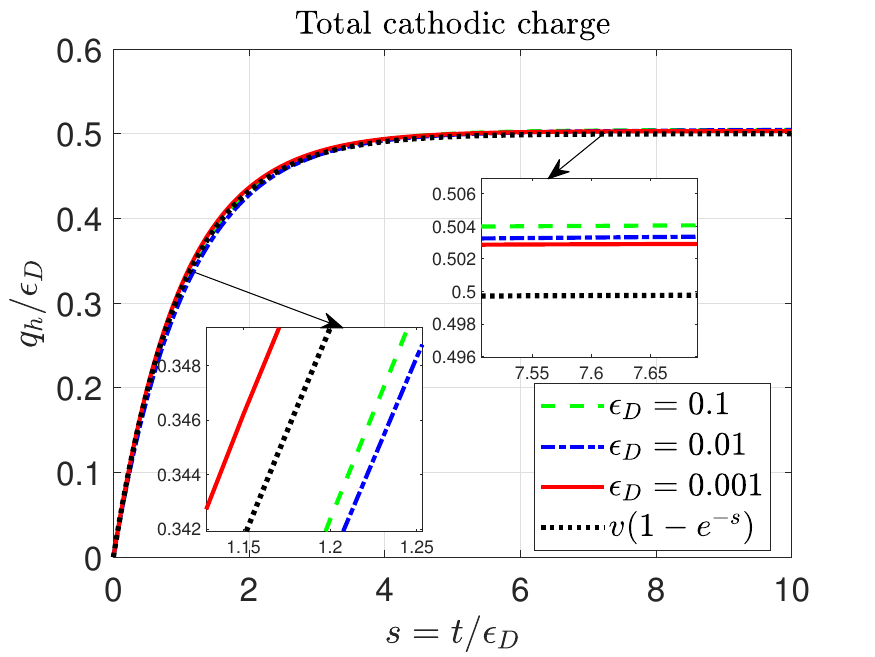}
		\end{minipage}
	}
	\subfigure{
		\begin{minipage}{0.31\textwidth}
			\includegraphics[width=\linewidth]{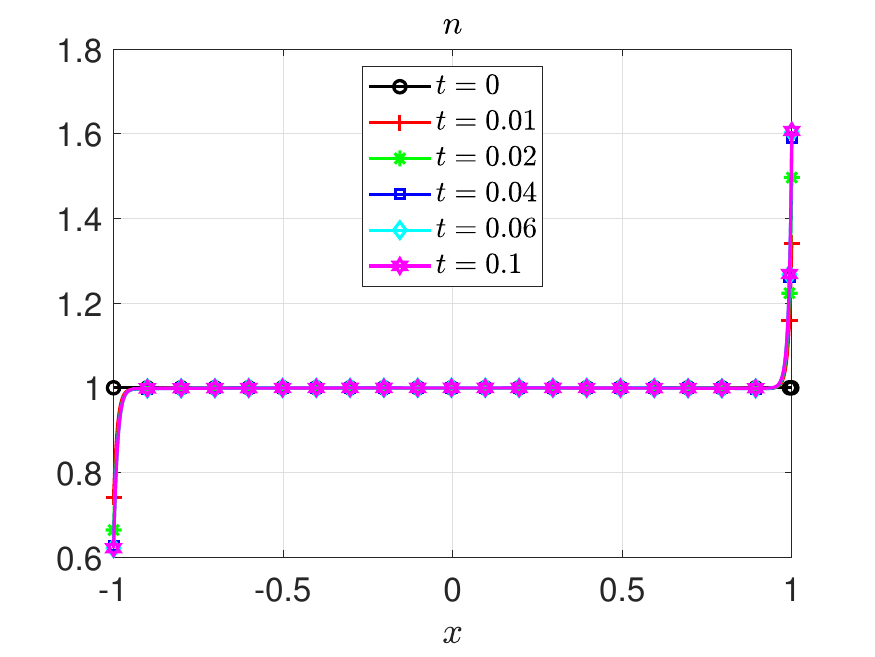}
			\includegraphics[width=\linewidth]{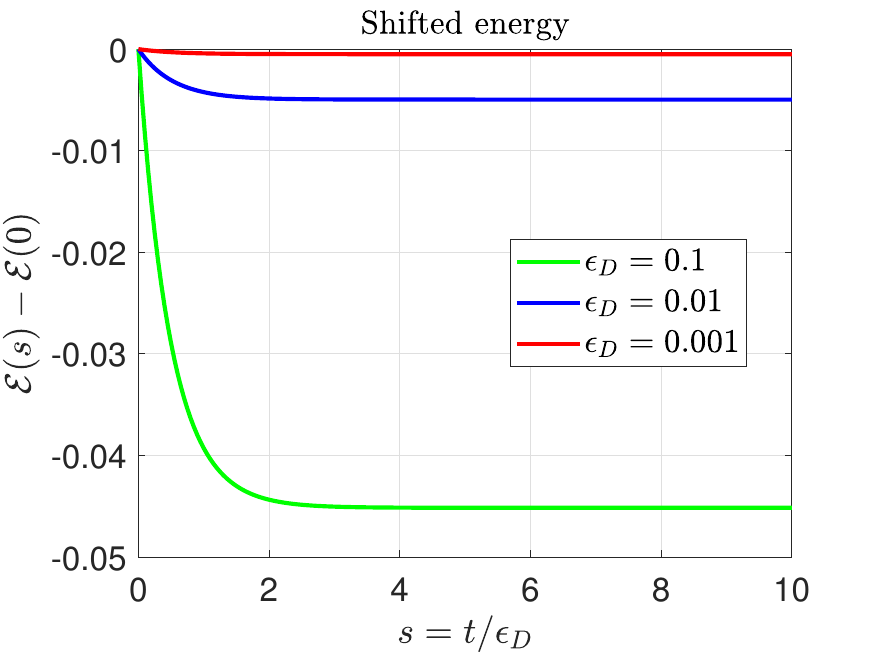}
		\end{minipage}
	}
	\subfigure{
		\begin{minipage}{0.31\textwidth}
			\includegraphics[width=\linewidth]{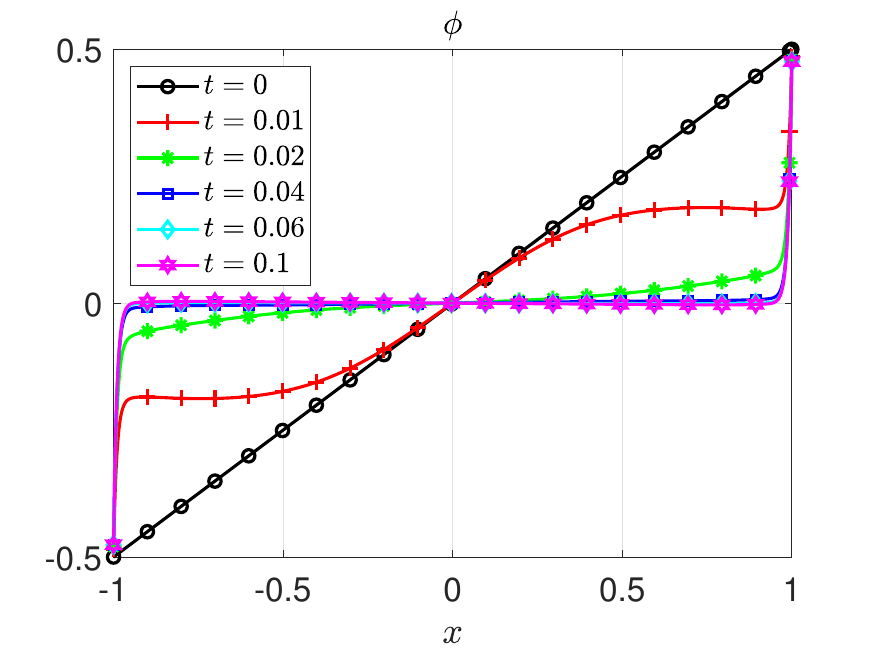}
			\includegraphics[width=\linewidth]{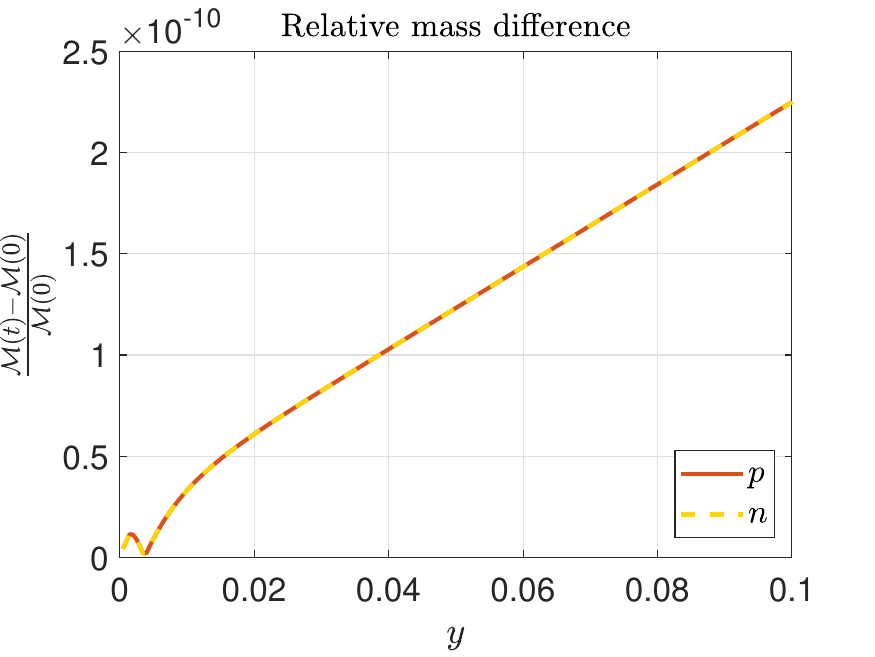}
		\end{minipage}
	}
	\caption{Weak-voltage diffuse-charge dynamics for $v=0.5$. Top: evolution of $p$, $n$, and $\phi$ for $\epsilon_{D}=0.01$. Bottom: normalized cathodic charge and the linear RC asymptote (left), shifted discrete energy (middle), and relative mass changes for $\epsilon_{D}=0.01$ (right).}
	\label{fig:1DCharge_Ex1}
\end{figure}

We first simulate the weak-voltage thin-double-layer test with $v=0.5$ for $0\leq s\leq10$ (with 200 JKO steps). The simulation results for $\epsilon_{D}=0.01$ (with $N_{x}=512$ cells) is shown in Fig.~\ref{fig:1DCharge_Ex1} (top). We observed that the initially uniform concentrations develop opposite enrichment and depletion layers near the electrodes, while remaining close to the electroneutral state in the bulk. Correspondingly, the initially linear potential relaxes toward an almost field-free bulk, with most of the voltage drop confined to the diffuse layers. 
We further compare the numerical result of the cathodic-charge curve with the analytical asymptotic solution for $\epsilon_{D}\in\{0.1,0.01,0.001\}$ in Fig.~\ref{fig:1DCharge_Ex1} (bottom). The numerical results show perfect match with the analytic solutions, and the long-time equilibrium state is better captured by the asymptotic solution as $\epsilon_{D}$ decreases.

Our method exhibits superior robustness for small permittivity ($\epsilon=2\times 10^{-6}$) and still preserves the desired properties of energy dissipation, positivity-preserving and mass-conservation (Fig.~\ref{fig:1DCharge_Ex1} (bottom)). Table~\ref{tab:BTA_PrePD_VPTPD} reports the computational efficiency of two primal--dual methods for the weak-voltage runs with small permittivity. For both methods, the iteration number only mildly increases as $\epsilon_{D}$ decreases, while the CPU time per iteration increases almost linearly with $1/\epsilon_{D}$. VPTPD consistently requires much fewer iterations and less CPU time than PrePD, and its advantage becomes increasingly pronounced as $\epsilon_{D}$ decreases.

\begin{table}[htbp]
    \centering
	\tabcaption{Comparison of PrePD and VPTPD for the weak-voltage diffuse-charge test over 200 JKO steps.}
    \begin{tabular}{ccccccc}
        \toprule
        $\epsilon_{D}$ & $\epsilon$ & $\text{Method}$ & $\text{Mean Iter}$ & $\text{Total Iter}$ & $\text{CPU (s)}$ & $\text{CPU/Iter (s)}$ \\
        \midrule
        \multirow{2}{*}{$0.1$}
        & \multirow{2}{*}{$2\times10^{-2}$}
		& PrePD  & 205 & 41027 & 177.63 & $4.330\times10^{-3}$ \\
		& & VPTPD & 76 & 15117 & 7.14 & $4.725\times10^{-4}$ \\
        \midrule
        \multirow{2}{*}{$0.01$}
        & \multirow{2}{*}{$2\times10^{-4}$}
		& PrePD  & 3804 & 760895 & 6085.20 & $7.997\times10^{-3}$ \\
		& & VPTPD & 272 & 54365 & 107.64 & $1.980\times10^{-3}$ \\
        \midrule
        \multirow{2}{*}{$0.001$}
        & \multirow{2}{*}{$2\times10^{-6}$}
		& PrePD  & 16938 & 3387502 & 181622.89 & $5.362\times10^{-2}$ \\
		& & VPTPD & 564 & 112778 & 1623.74 & $1.440\times10^{-2}$ \\
        \bottomrule
        \label{tab:BTA_PrePD_VPTPD}
    \end{tabular}
\end{table}

\begin{figure}[htbp]
	\centering
	\subfigure{
		\begin{minipage}{0.31\textwidth}
			\includegraphics[width=\linewidth]{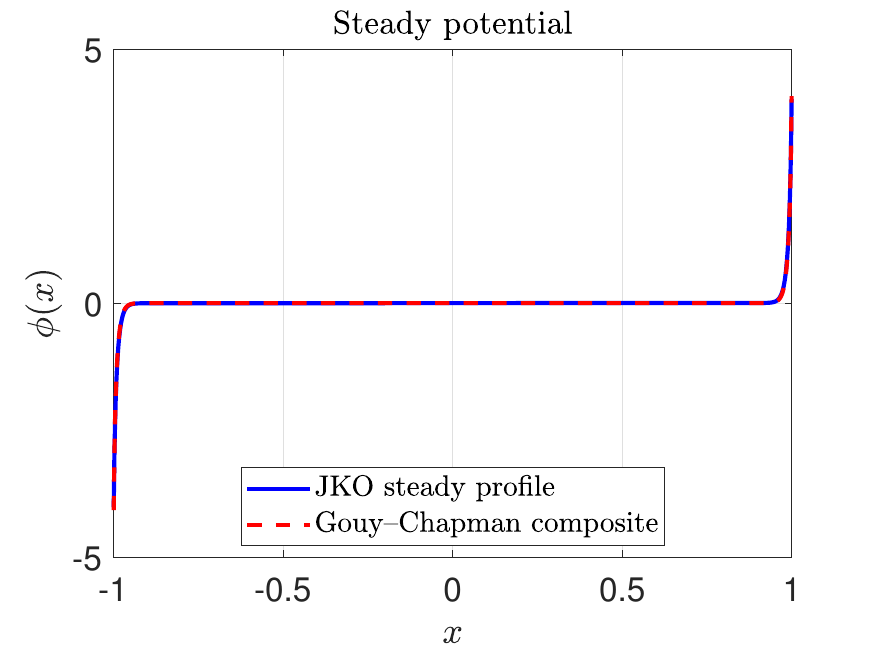}
		\end{minipage}
	}
	\subfigure{
		\begin{minipage}{0.31\textwidth}
			\includegraphics[width=\linewidth]{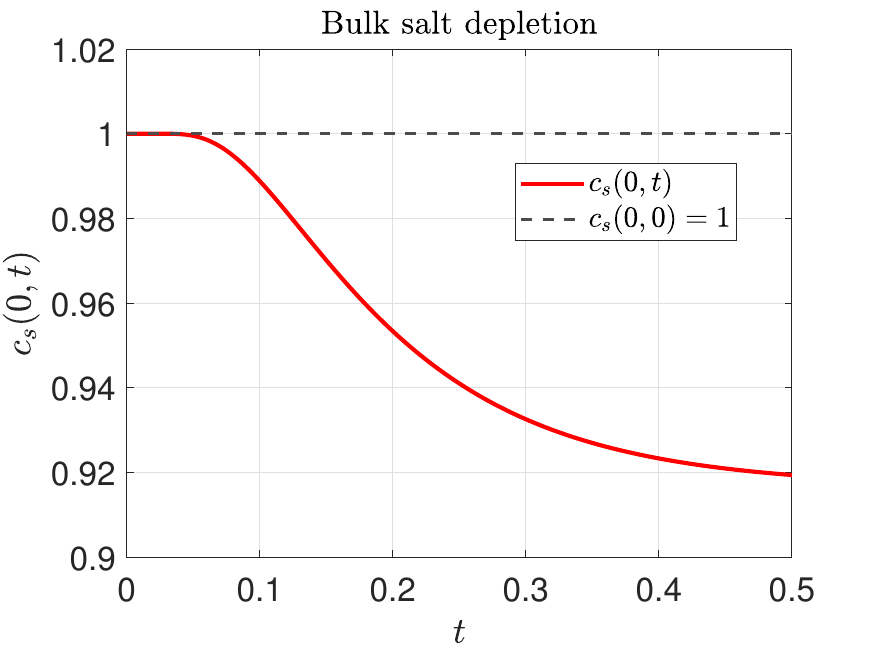}
		\end{minipage}
	}
	\subfigure{
		\begin{minipage}{0.31\textwidth}
			\includegraphics[width=\linewidth]{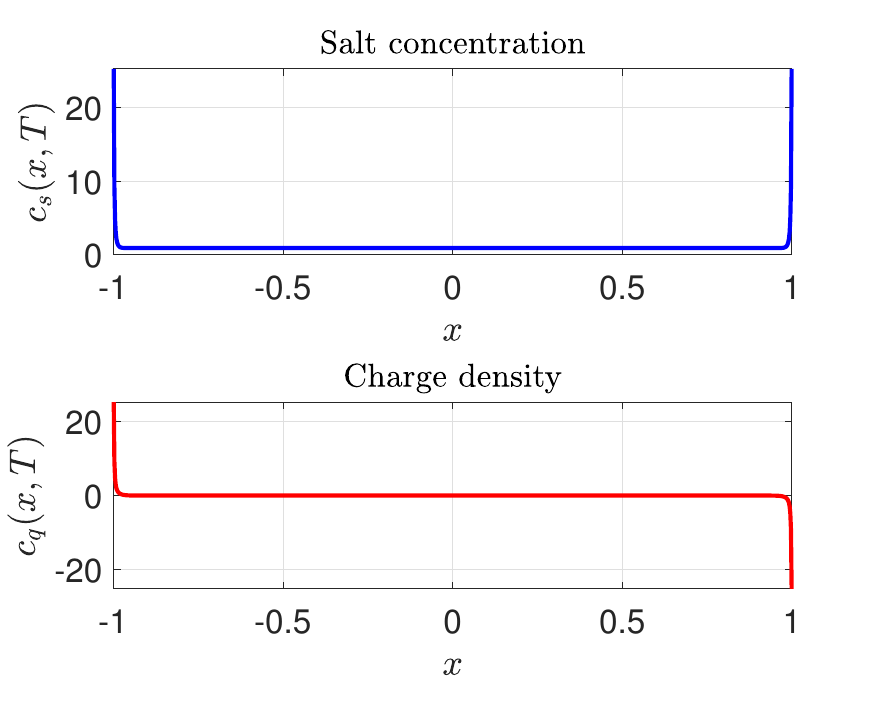}
		\end{minipage}
	}
	\caption{Strongly nonlinear diffuse-charge dynamics for $v=5$ and $\epsilon_{D}=0.01$. Left: final-time numerical potential and the Gouy--Chapman composite approximation. Middle: evolution of the bulk salt concentration $c_{s}(0,t)$. Right: final salt concentration $c_{s}(x,T)$ (top) and charge density $c_{q}(x,T)$ (bottom).}
	\label{fig:1DCharge_Ex2}
\end{figure}

We next consider the strongly nonlinear charging regime with $v=5$ and $\epsilon_{D}=0.01$ ($N_{x}=1024$). In this regime, neutral-salt adsorption by the diffuse layers produces an appreciable bulk response on the diffusion time scale \cite{Bazant2004Diffuse}, where the salt concentration and space-charge density are defined respectively as $c_{s}=(p+n)/2$ and $c_{q}=(p-n)/2$. We perform numerical simulations and monitor the bulk salt concentration $c_{s}(x=0,t)$ in the middle between two electrodes. 
Fig.~\ref{fig:1DCharge_Ex2} shows that the final-time numerical solution of the potential $\phi(x)$ (at $T=0.5$) agrees with the analytic leading-order Gouy--Chapman composite profile \cite{Bazant2004Diffuse}. The decrease of $c_{s}(0,t)$ from its initially uniform level demonstrates neutral-salt depletion in the bulk as predicted by the analysis. At the final time, excess salt and equal-and-opposite diffuse charge are strongly localized near the two electrodes, whereas the interior remains approximately electroneutral. These results are consistent with the leading-order interfacial structure and capture the slower bulk-diffusion response described in \cite{Bazant2004Diffuse}.

\subsection{Comparison tests}\label{sec:6.2}
\subsubsection{Comparison with existing numerical methods}\label{sec:6.1.2}

\begin{figure}[htbp]
	\centering
	\subfigure{
		\begin{minipage}{0.23\textwidth}
			\includegraphics[width=\linewidth]{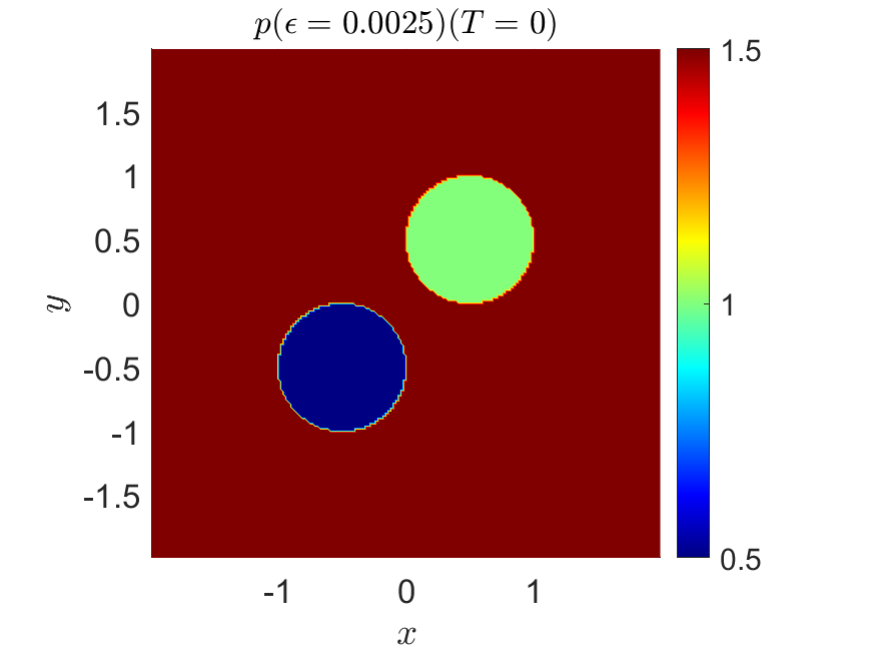}
			\includegraphics[width=\linewidth]{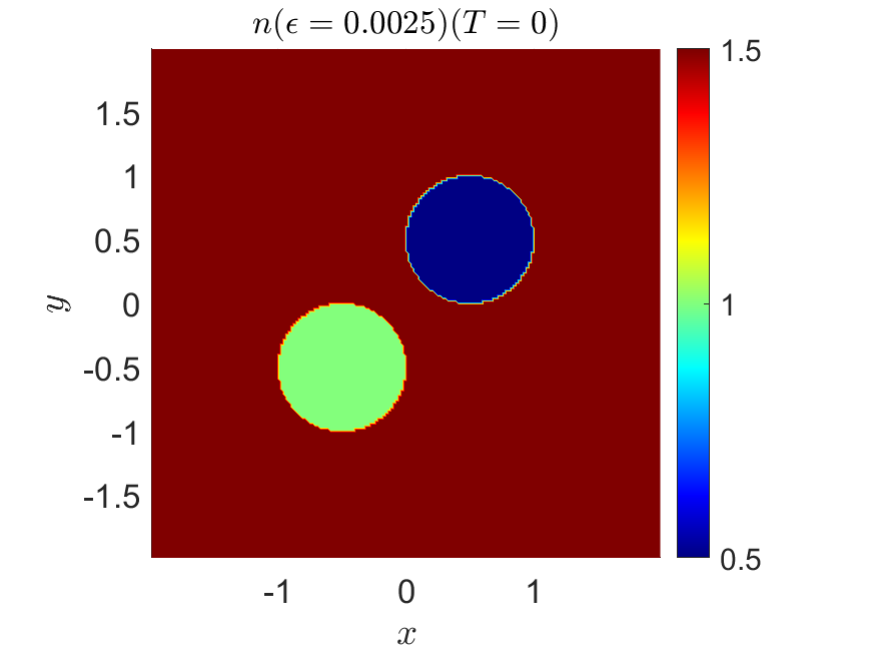}
			\includegraphics[width=\linewidth]{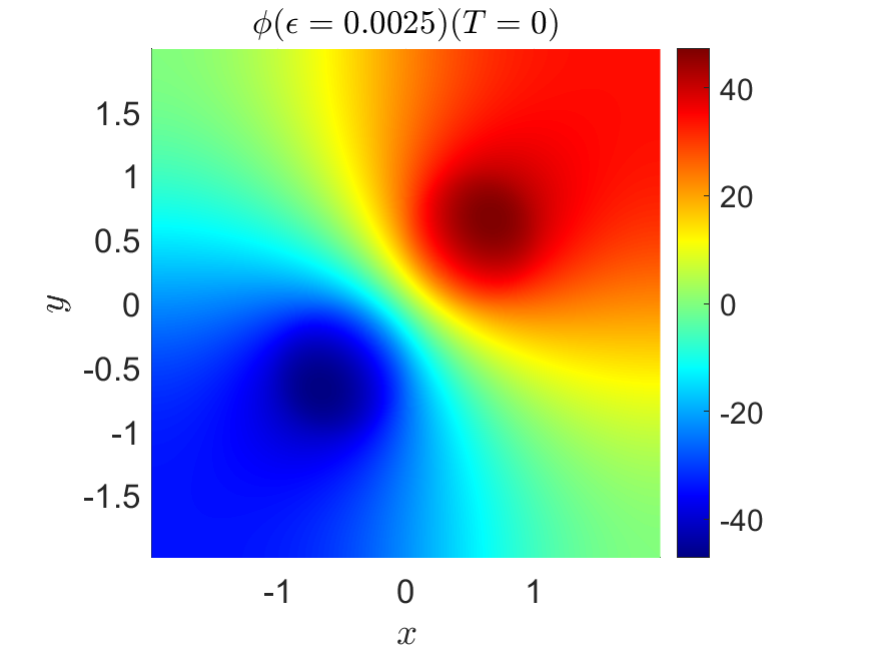}
		\end{minipage}
	}
	\subfigure{
		\begin{minipage}{0.23\textwidth}
			\includegraphics[width=\linewidth]{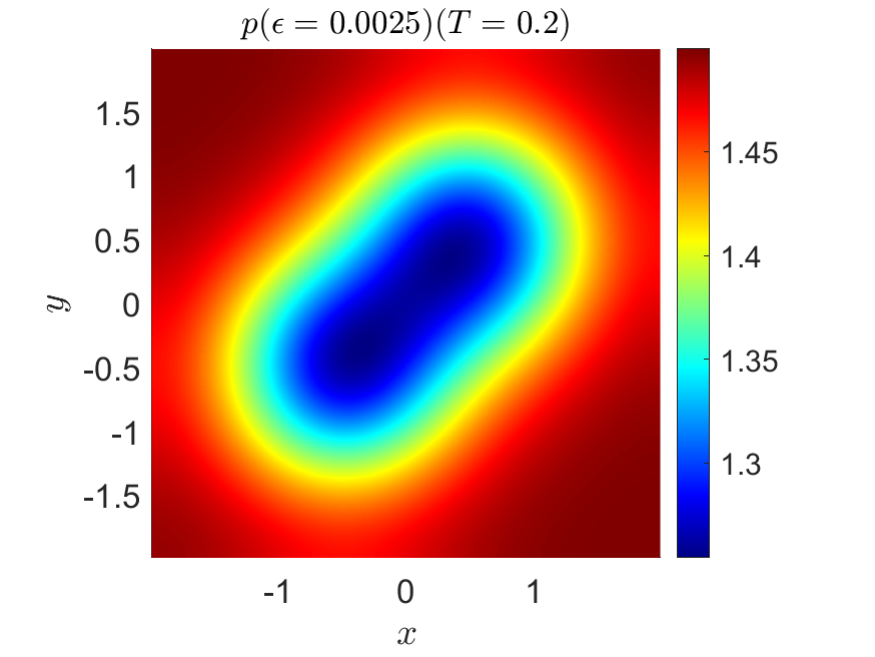}
			\includegraphics[width=\linewidth]{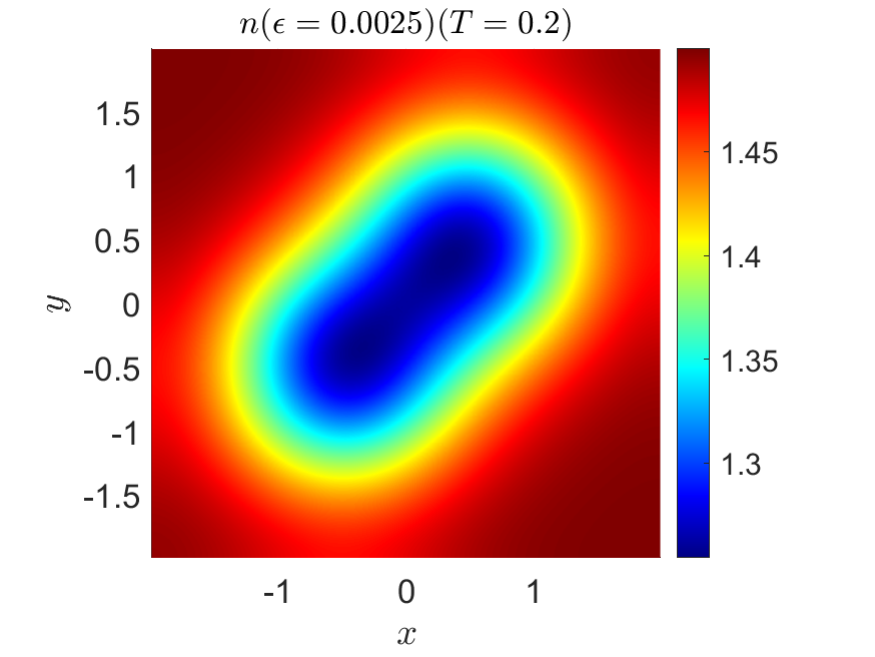}
			\includegraphics[width=\linewidth]{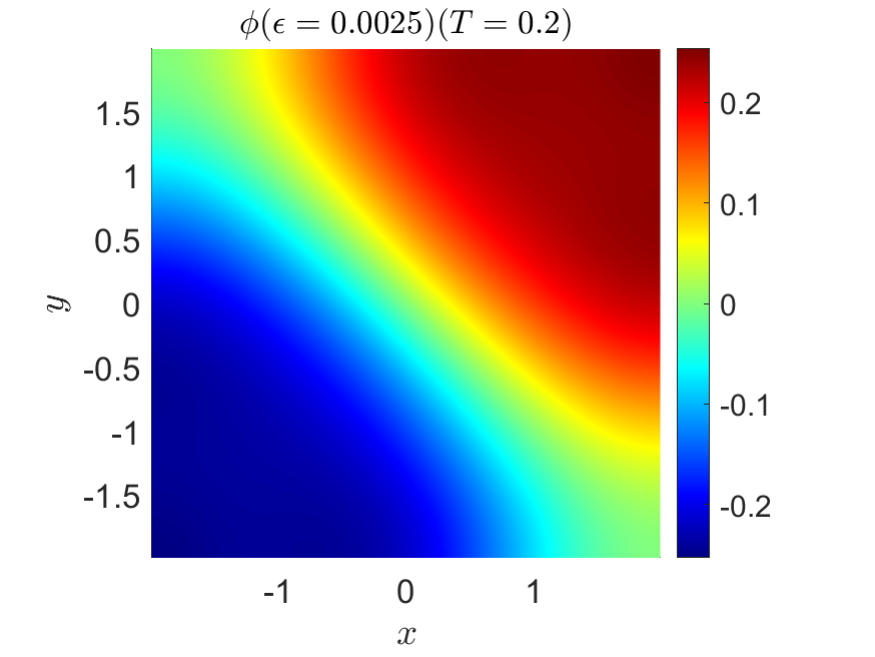}
		\end{minipage}
	}
	\subfigure{
		\begin{minipage}{0.23\textwidth}
			\includegraphics[width=\linewidth]{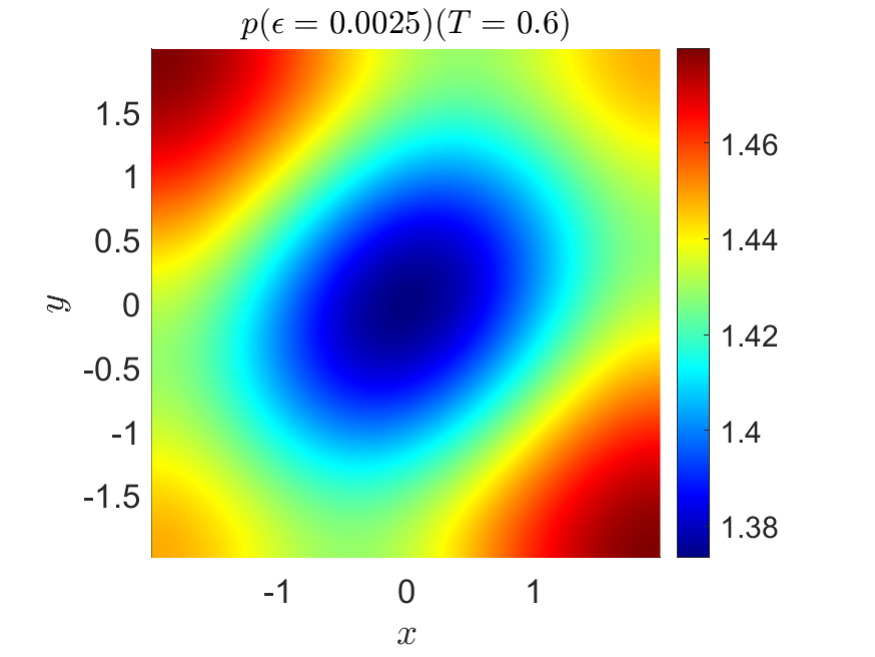}
			\includegraphics[width=\linewidth]{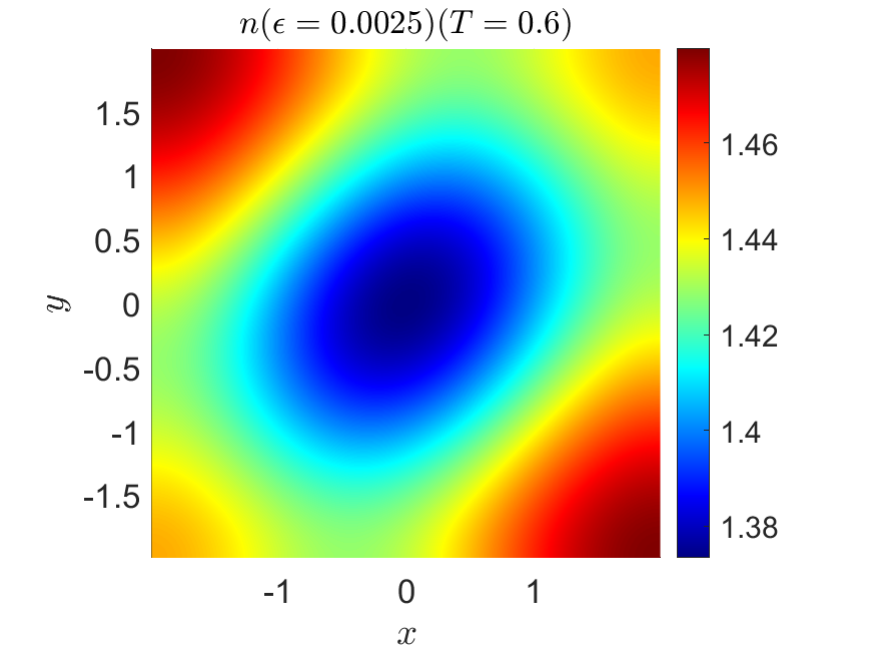}
			\includegraphics[width=\linewidth]{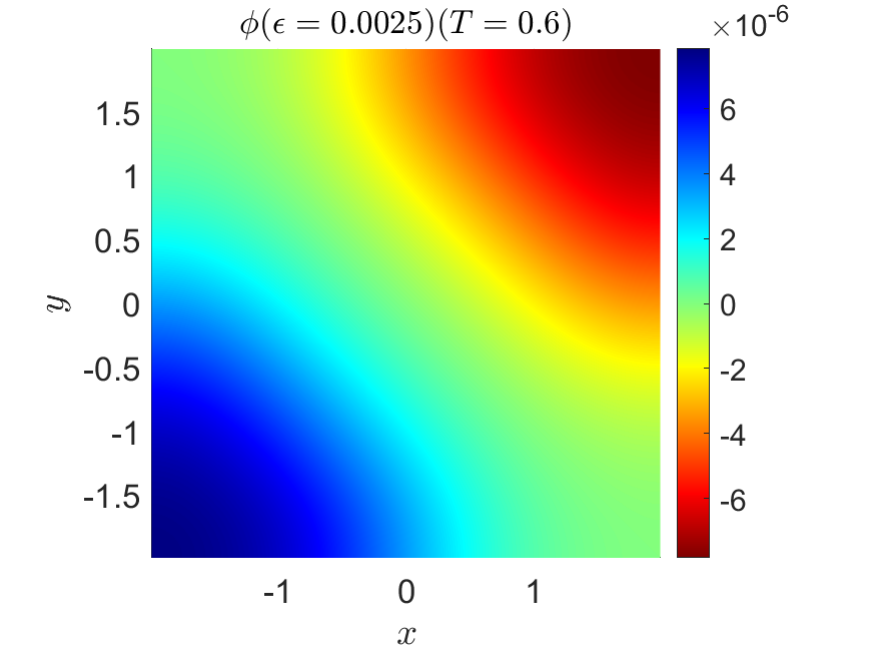}
		\end{minipage}
	}
	\subfigure{
		\begin{minipage}{0.23\textwidth}
			\includegraphics[width=\linewidth]{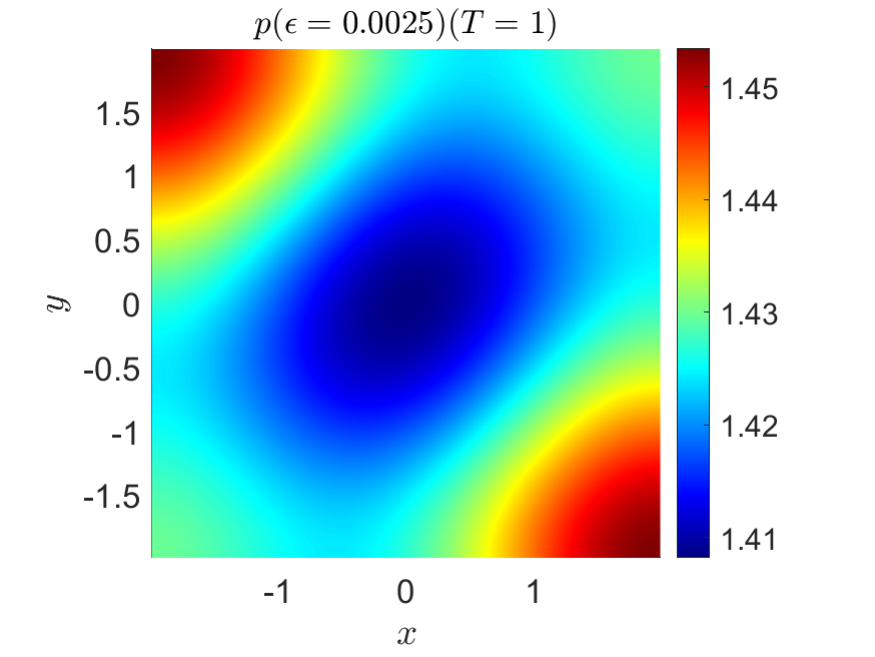}
			\includegraphics[width=\linewidth]{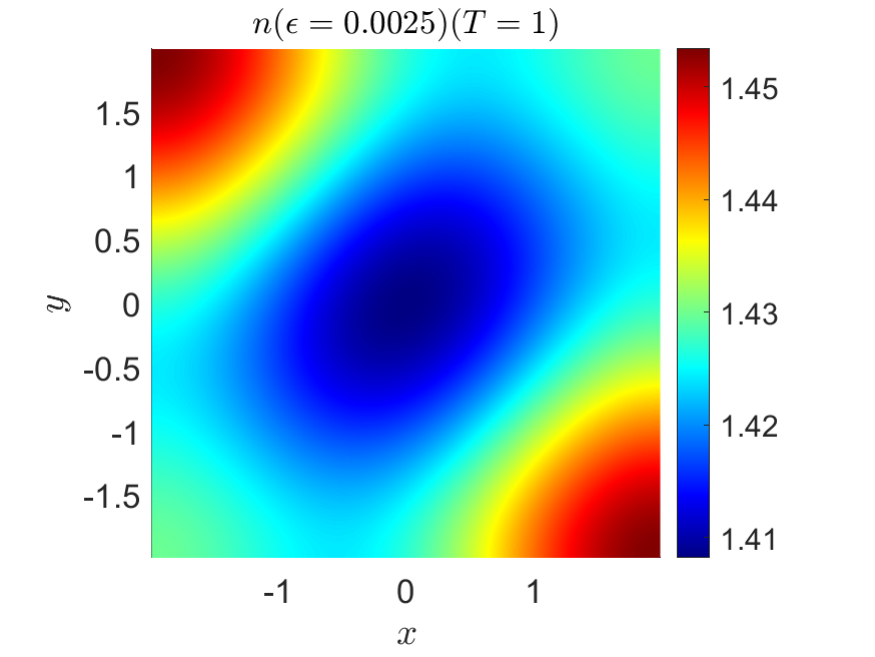}
			\includegraphics[width=\linewidth]{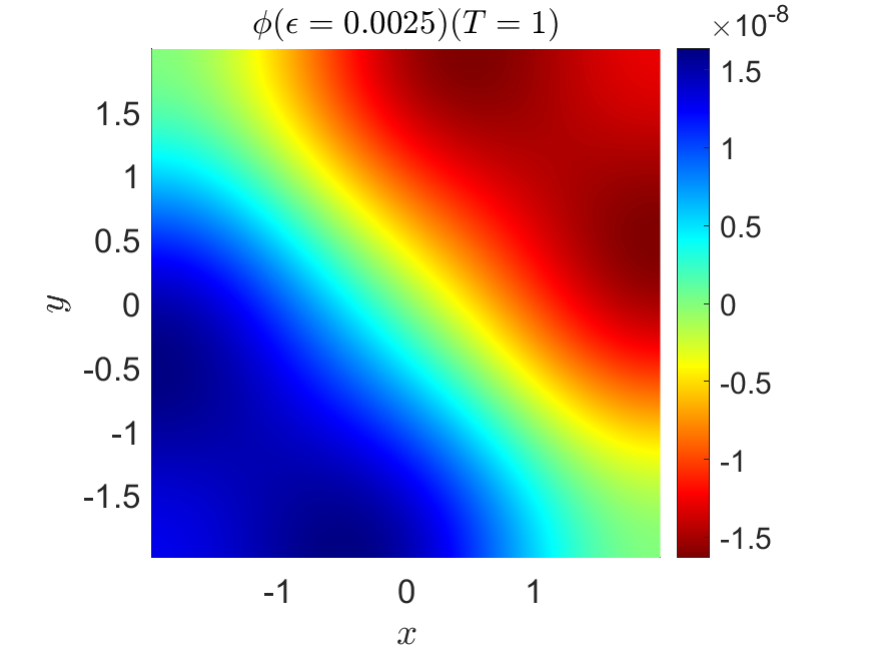}
		\end{minipage}
	}
	\caption{Time evolution snapshots of $p$, $n$, and $\phi$ for the 2D test with initial conditions \eqref{eq:2Dinitial} and Neumann boundary conditions for $\epsilon=0.0025$.}
	\label{fig:2DJKOsmalleps}
\end{figure}

We first compare our Poisson-constrained JKO scheme with other structure-preserving methods for PNP models, specifically the projection method (PJM) \cite{Tong2024PositivityPNP}. We consider a two-dimensional PNP model in a domain $\left[-2,2\right]^{2}$ with homogeneous Neumann boundary conditions on the electrostatic potential with the initial conditions:
\begin{eqnarray}\label{eq:2Dinitial}
	\begin{aligned}
		&(p(x,y,0),n(x,y,0)) =
		\begin{cases}
			(1.0, 0.5),& \quad(x-0.5)^{2}+(y-0.5)^{2}<0.25, \\
			(0.5, 1.0),& \quad(x+0.5)^{2}+(y+0.5)^{2}<0.25, \\
			(1.5, 1.5),& \quad\text{otherwise}.
		\end{cases}
	\end{aligned}
\end{eqnarray}
We solve the Poisson-constrained JKO scheme by PrePD method (with DCT-DS for dual subproblem in Sec.~\ref{sec:DCT-DS}) with $T=1$, $\tau=1/100$, and $\Delta x=\Delta y=4/200$. Fig.~\ref{fig:2DJKOsmalleps}  displays the time evolution of $p$, $n$, and $\phi$ for $\epsilon=0.0025$. With small $\epsilon$, the distributions of electrons $n$ and holes $p$ approaches the local electroneutrality constraint ($p-n+\psi^{0}\approx 0$), while the electrostatic potential $\phi$ approaches a constant field with homogeneous Neumann boundary conditions.

\begin{figure}[htbp]
	\centering
	\subfigure[JKO and PJM at $\epsilon=1$]{
		\begin{minipage}{0.31\textwidth}
			\includegraphics[width=\linewidth]{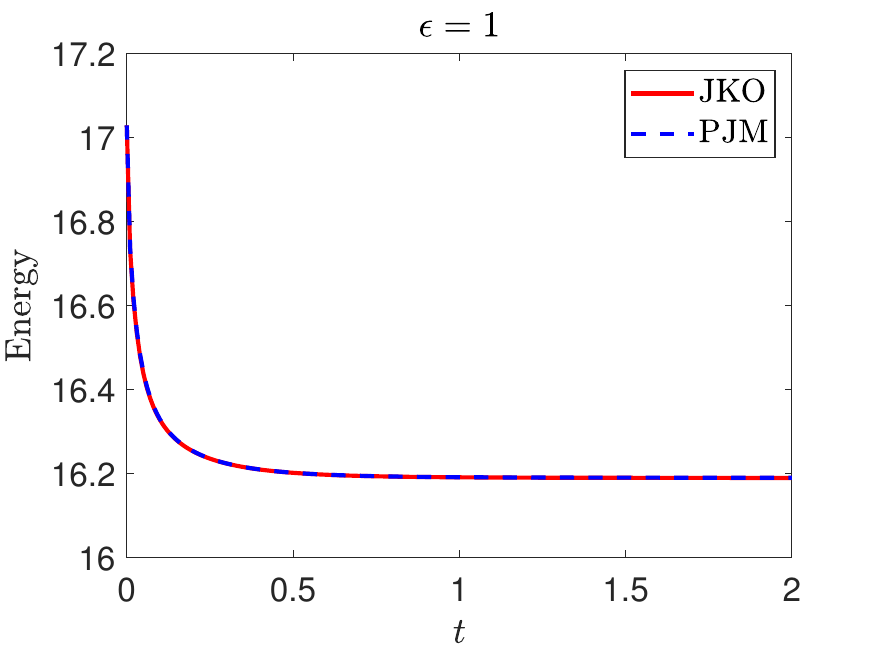}
		\end{minipage}
	}
	\subfigure[JKO solvers and PJM at $\epsilon=0.0025$]{
		\begin{minipage}{0.31\textwidth}
			\includegraphics[width=\linewidth]{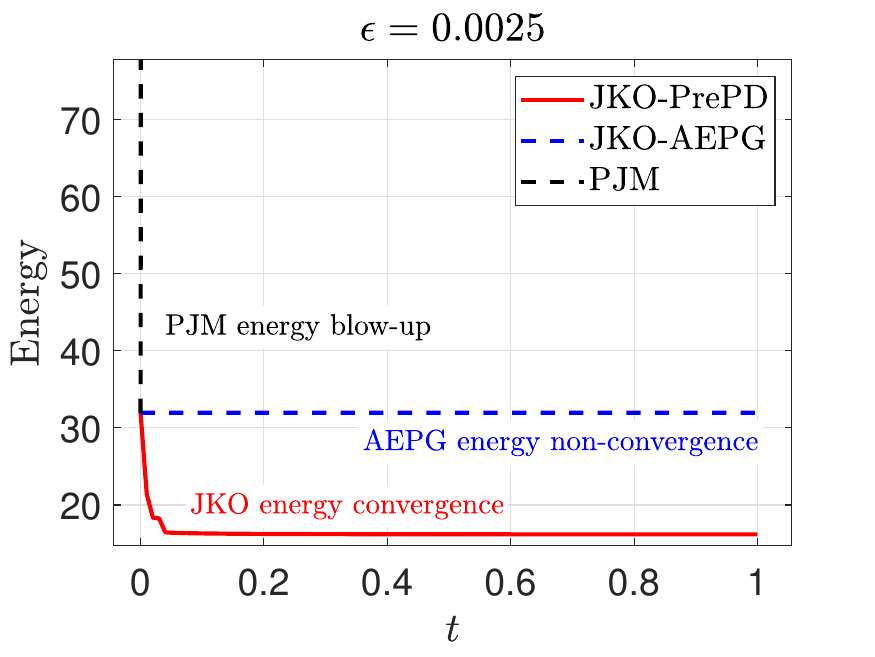}
		\end{minipage}
	}
	\subfigure[JKO at $\epsilon=0.0025$]{
		\begin{minipage}{0.31\textwidth}
			\includegraphics[width=\linewidth]{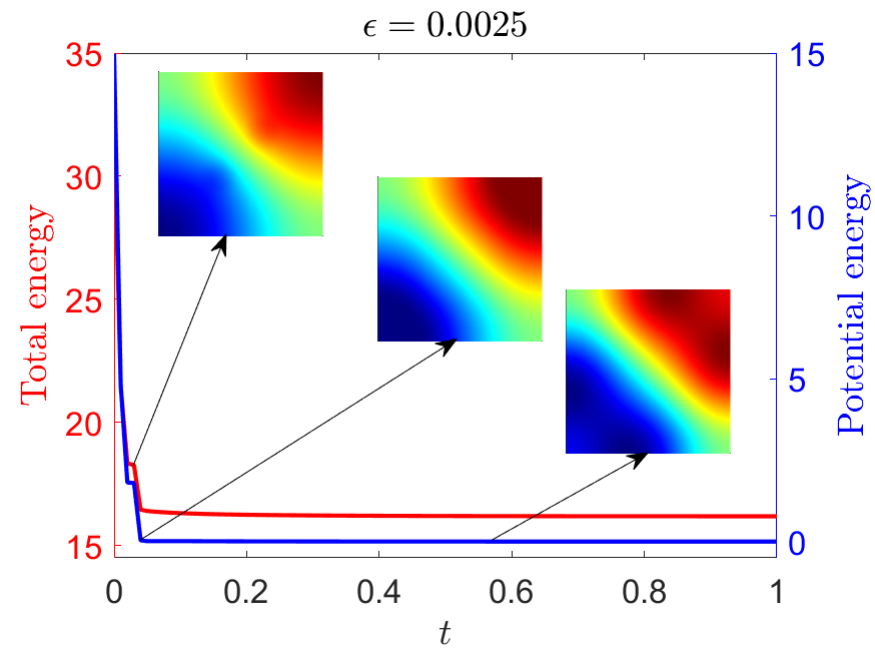}
		\end{minipage}
	}
	\caption{Energy evolution of 2D test with initial condition \eqref{eq:2Dinitial} and Neumann boundary conditions for different numerical schemes. (a) JKO v.s. PJM for $\epsilon=1$. (b) JKO-PrePD v.s. JKO-AEPG v.s. PJM for $\epsilon=0.0025$. (c) Total energy and electrostatic energy evolution computed by JKO-PrePD for $\epsilon=0.0025$.}
	\label{fig:threescheme}
\end{figure}

We compare the performance of the JKO scheme with PJM \cite{Tong2024PositivityPNP}. For $\epsilon=1$, both methods exhibit stable results (Fig.~\ref{fig:threescheme}a); however, for $\epsilon=0.0025$, PJM fails to maintain stability since it does not guarantee unconditionally energy stability while JKO scheme still preserves stability and energy dissipation with $\tau=1/100$ and $\Delta x=\Delta y=4/200$ (Fig.~\ref{fig:threescheme}b,c). PJM can handle the case of $\epsilon=0.0025$ with a smaller time step $\tau=10^{-3}$ and a finer mesh $\Delta x=\Delta y=4/400$, at similar cost (CPU time) with the JKO scheme (Table~\ref{tab:JKOPJM}). This shows that the JKO scheme, while requiring to solve minimization at each time step, is more robust in small-permittivity regimes and does not substantially increase the computational time.

\begin{table}[htbp]
	\centering
	\caption{Comparison between the JKO scheme and the PJM scheme \cite{Tong2024PositivityPNP} for 2D test with Neumann boundary conditions for $\epsilon=0.0025$.}
	\begin{tabular}{ccccc}
		\toprule
		$\text{Method}$ & $\text{Mesh size}$ & $\text{Time size}$ & $\text{CPU time}$ \\ 
		\midrule
		JKO  & 0.02 & $10^{-2}$ & 1739.71 \\
		PJM  & 0.01 & $10^{-3}$ & 1937.66 \\
		\bottomrule
	\end{tabular}
	\label{tab:JKOPJM}
\end{table}

We also compare the PrePD method with other algorithm for the JKO scheme of PNP models, specifically the adaptive energy-based preconditioned gradient (AEPG) method \cite{Li2026OnsagerPNP}. Table~\ref{tal:PrePDAEPG} shows that PrePD requires fewer iterations (around $\tfrac{1}{5}$ to $\tfrac{1}{3}$ iterations) and less total CPU time than AEPG for moderate values of $\epsilon=1, 0.5, 0.1$. For small values of $\epsilon=0.0025$, AEPG somehow fails to converge within the maximum number of iterations, while PrePD still converges and preserves the monotone discrete-energy decay (Fig.~\ref{fig:threescheme}b), demonstrating the robustness of the primal--dual splitting method for this kind of nonsmooth constrained convex optimization problems.	

\begin{table}[htbp]
	\centering
	\tabcaption{Comparison of computational efficiency between PrePD and AEPG \cite{Li2026OnsagerPNP} over 100 JKO steps for 2D test with Neumann boundary conditions.}
	\begin{tabular}{cccccc}
		\toprule
		$\epsilon$ & $\text{Method}$ & $\text{Mean Iter}$ & $\text{Total Iter}$ & $\text{CPU}$ & $\text{CPU/Iter}$ \\
		\midrule
		\multirow{2}{*}{$1$}
		& JKO-PrePD  & 49 & 4948 & \textbf{418.20} & 0.084 \\
		& JKO-AEPG   & 184 & 18412 & 1215.68 & \textbf{0.066} \\
		\midrule
		\multirow{2}{*}{$0.5$}
		& JKO-PrePD  & 53 & 5335 & \textbf{559.30} & 0.104 \\
		& JKO-AEPG   & 199 & 19946 & 1313.93 & \textbf{0.066} \\
		\midrule
		\multirow{2}{*}{$0.1$}
		& JKO-PrePD  & 153 & 15333 & \textbf{1576.80} & 0.103 \\
		& JKO-AEPG   & 749 & 74930 & 4657.27 & \textbf{0.062} \\
		\bottomrule
		\label{tal:PrePDAEPG}
	\end{tabular}
\end{table}

\subsubsection{Performance of BGS and Schur-PCG dual solvers}\label{sec:6.1.3}

We now investigate the performance of the proposed dual solvers (BGS and Schur-PCG) in Sec.~\ref{sec:5} for the Poisson-constrained JKO scheme. We consider a two-dimensional test on $[0,1]^2$ with the initial and Dirichlet boundary conditions for the electrostatic potential \cite{Liu2023DynamicMassPNP}:
\begin{eqnarray}\label{eq:p0n0bcphi}
	\begin{aligned}
		&p^{0}=4x(1-x)+8y(1-y),\quad n^{0}=\sin(\pi x)+\sin(\pi y),\\
		&\phi(0,y)=0,\quad\phi(1,y)=0,\quad\phi(x,0)=0,\quad\phi(x,1)=0,
	\end{aligned}
\end{eqnarray}
where the fixed charge density is given by:
\begin{eqnarray}\label{eq:psi01}
	\begin{aligned}
		&\psi^{0}_{1} =
		\begin{cases}
			10, & \text{if } \dfrac{5}{8}\leq x\leq\dfrac{7}{8}\text{ and }\dfrac{5}{8} \leq y\leq\dfrac{7}{8}, \\
			0, & \text{otherwise}.
		\end{cases}
	\end{aligned}
\end{eqnarray}

\begin{figure}[htbp]
	\centering
	\subfigure{
		\begin{minipage}{0.23\textwidth}
			\includegraphics[width=\linewidth]{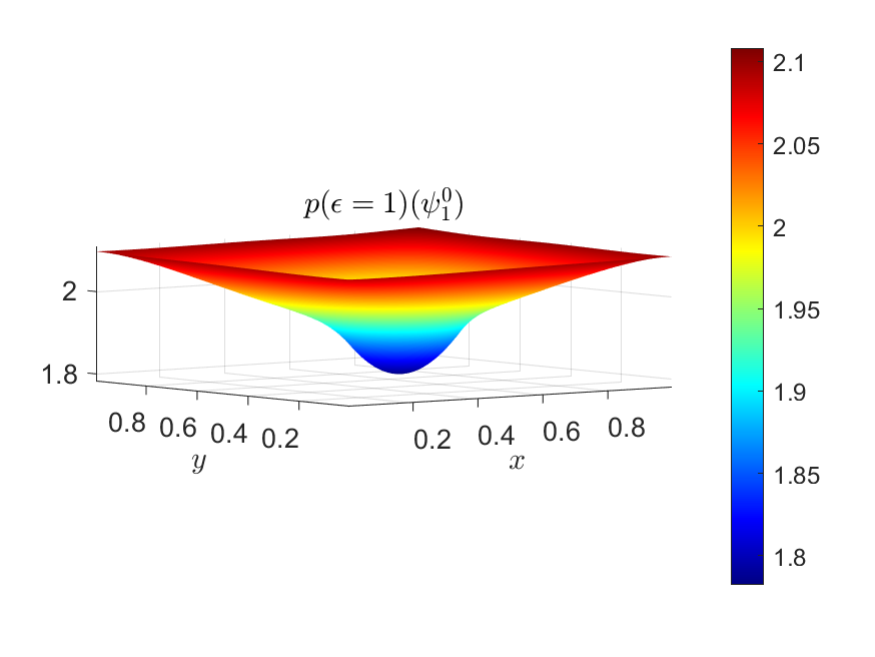}
			\includegraphics[width=\linewidth]{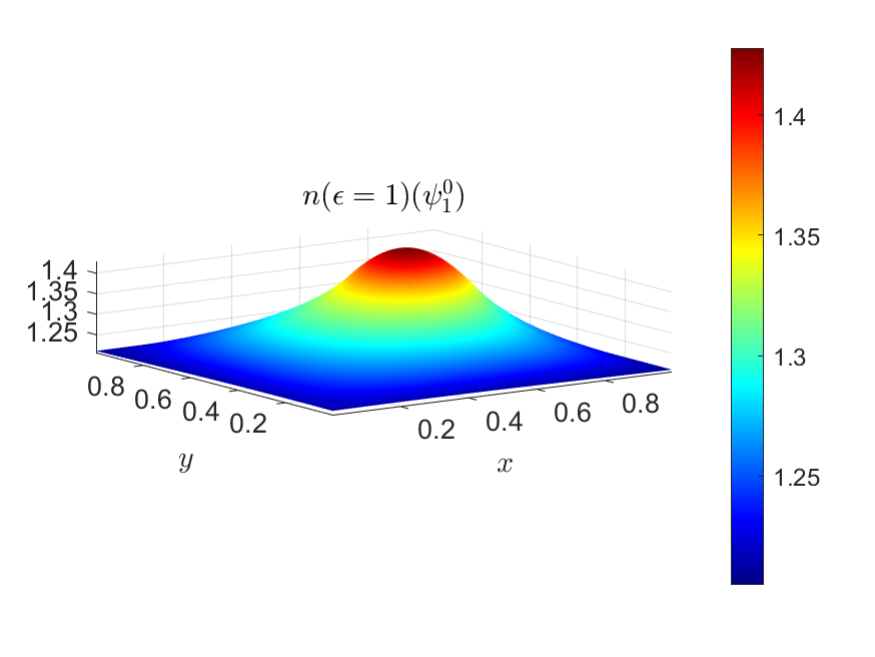}
			\includegraphics[width=\linewidth]{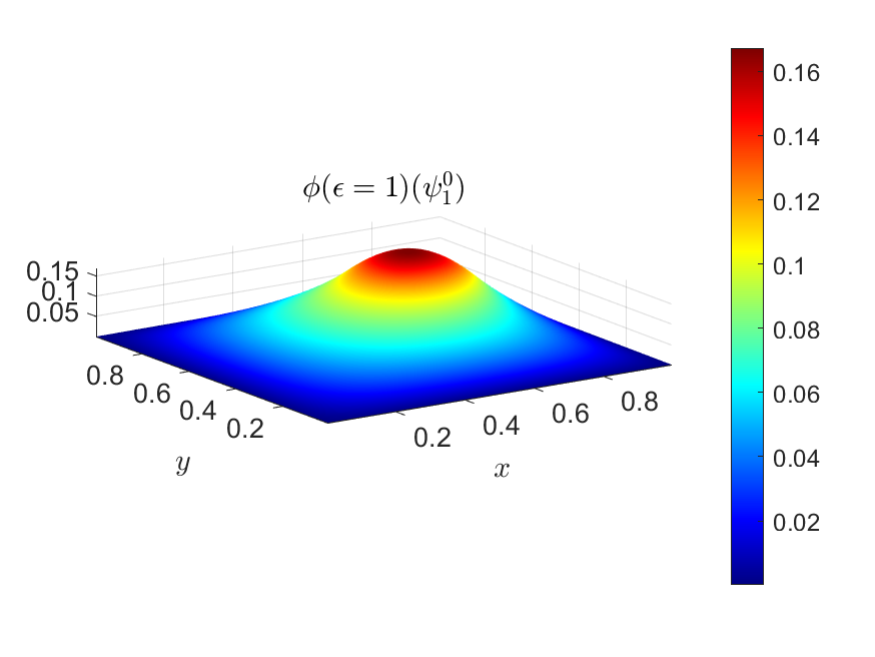}
		\end{minipage}
	}
	\subfigure{
		\begin{minipage}{0.23\textwidth}
			\includegraphics[width=\linewidth]{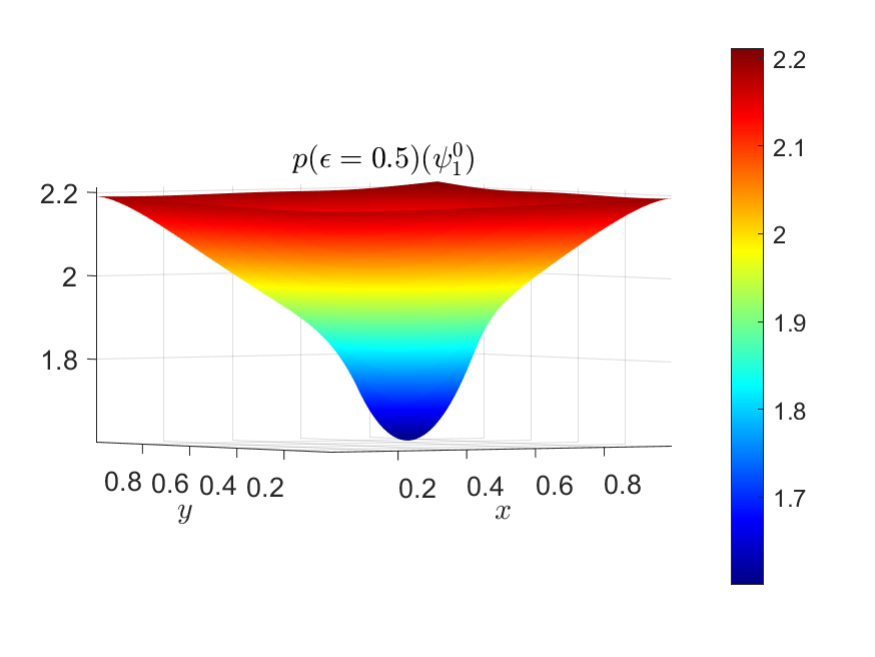}
			\includegraphics[width=\linewidth]{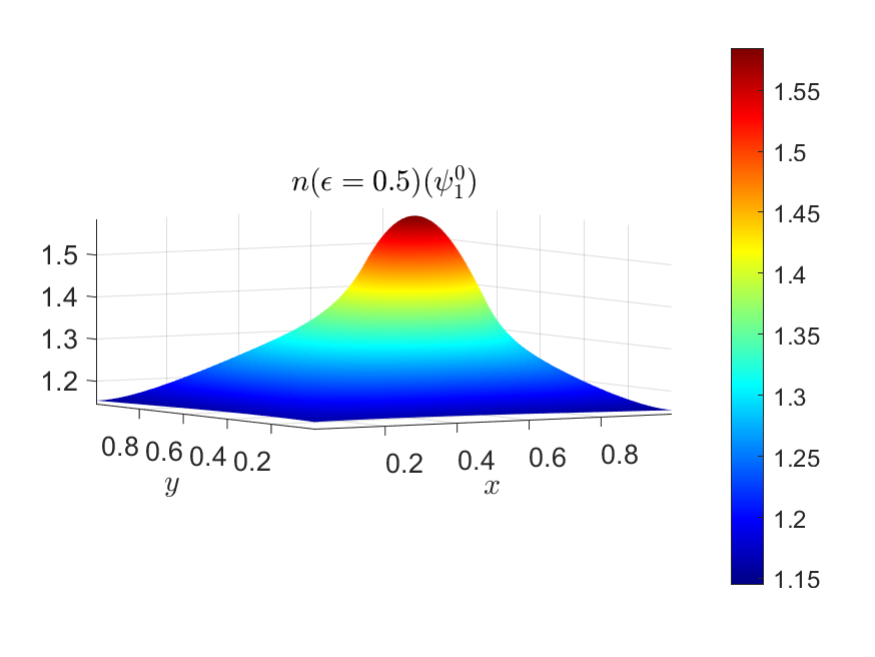}
			\includegraphics[width=\linewidth]{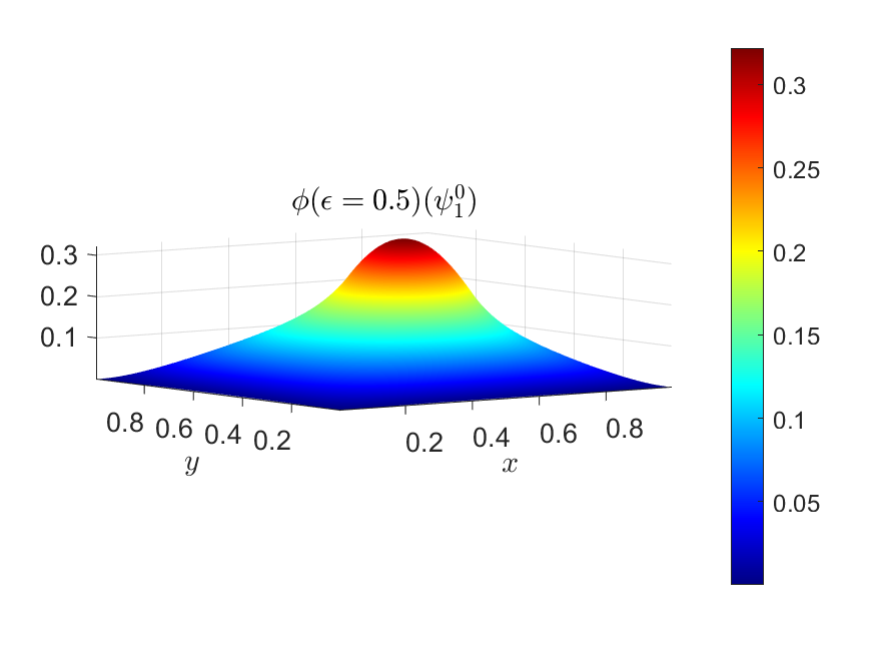}
		\end{minipage}
	}
	\subfigure{
		\begin{minipage}{0.23\textwidth}
			\includegraphics[width=\linewidth]{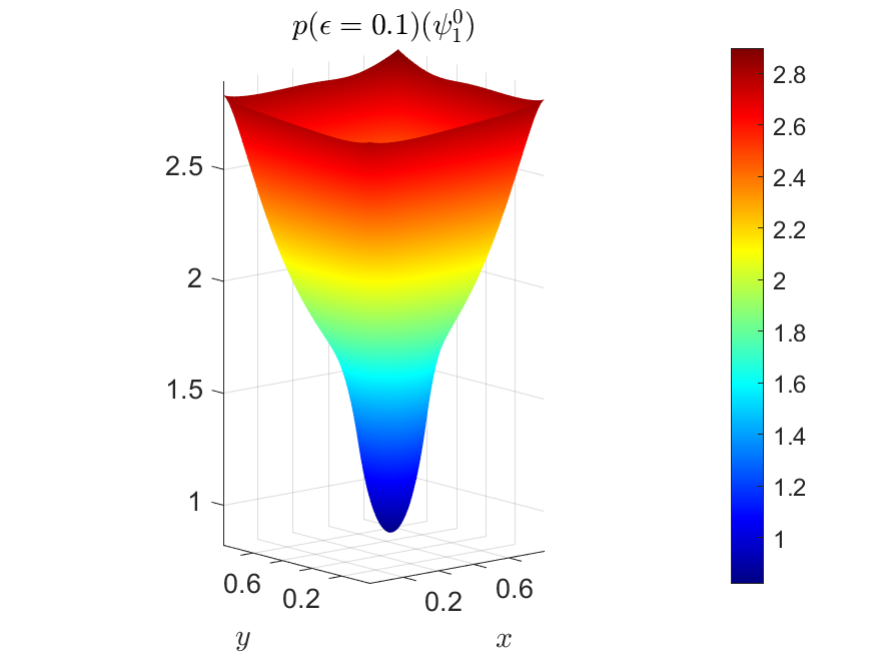}
			\includegraphics[width=\linewidth]{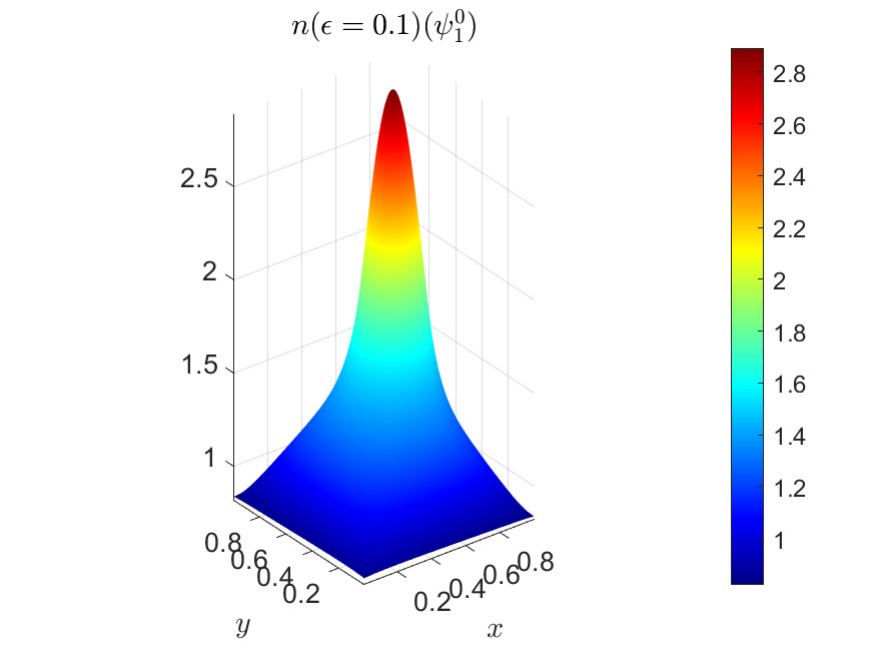}
			\includegraphics[width=\linewidth]{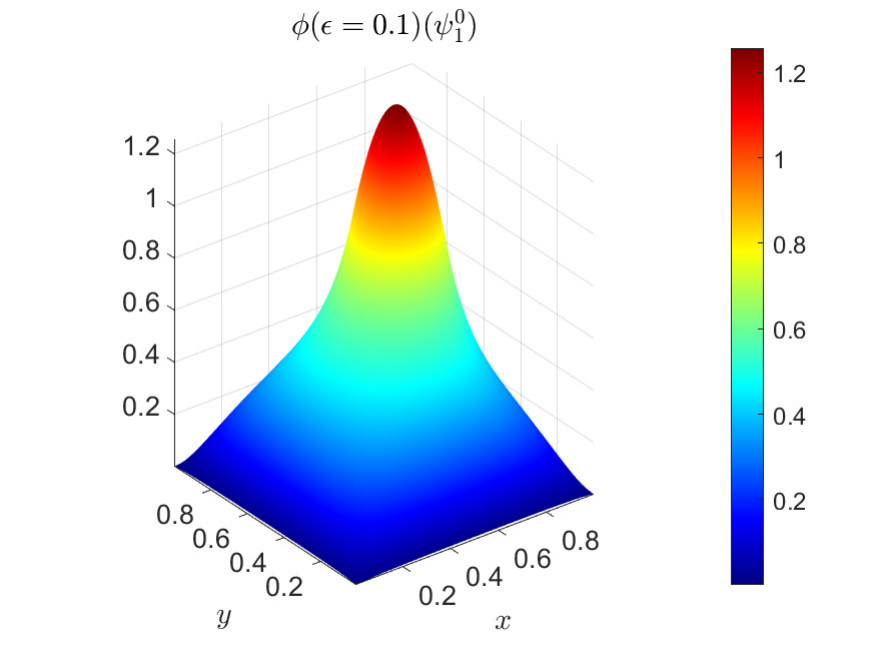}
		\end{minipage}
	}
	\subfigure{
		\begin{minipage}{0.23\textwidth}
			\includegraphics[width=\linewidth]{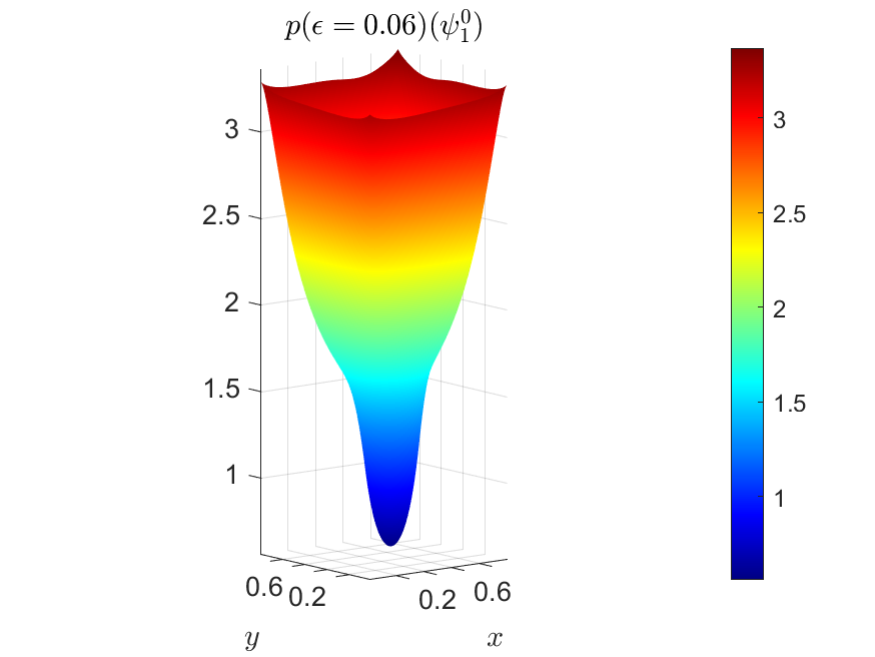}
			\includegraphics[width=\linewidth]{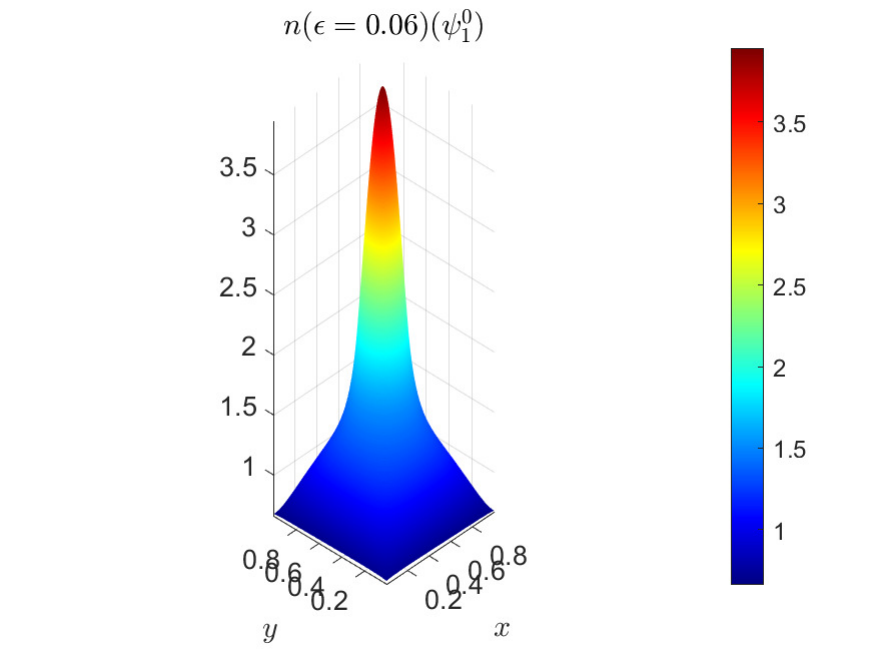}
			\includegraphics[width=\linewidth]{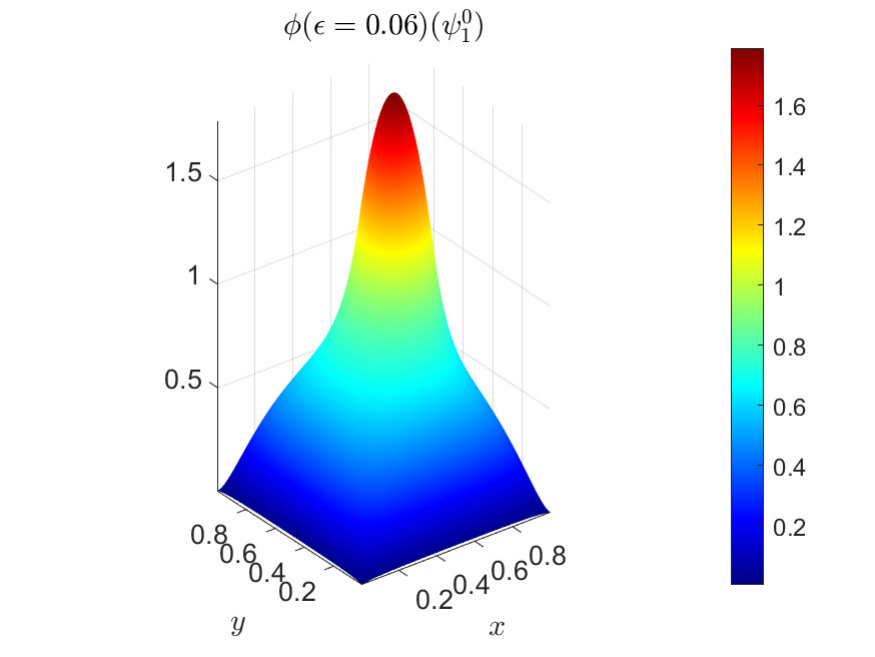}
		\end{minipage}
	}
	\caption{Profiles of the hole concentration $p$, electron concentration $n$, and electrostatic potential $\phi$ at $T=0.2$ for different values of $\epsilon$ with Dirichlet boundary conditions and fixed charge density \eqref{eq:psi01}.}
	\label{fig:diffepspandn}
\end{figure}

We solve the Poisson-constrained JKO scheme \eqref{eq:fullJKO} for the above equation with $T=0.2$, $\tau=1/100$, and $\Delta x=\Delta y=1/150$ by the PrePD method \eqref{eq:PrePD} with BGS or Schur-PCG for the dual subproblem. The profiles of the hole and electron concentrations ($p,n$) and the electrostatic potential ($\phi$) at $T=0.2$ for different values of $\epsilon$ are shown in Fig.~\ref{fig:diffepspandn}. The electrons accumulate near the positive fixed charge, whereas the holes are depleted from this region. The resulting redistribution of mobile charge partially compensates the fixed charge, thereby localizing the potential variation and confining the substantial electric field to a narrow region around the fixed charge. We observe that as $\epsilon$ decreases, the profiles develop sharper spatial variations and become increasingly localized near the charged region. This behavior is consistent with the Debye-length scaling $\epsilon_D\sim\sqrt{\epsilon}$ \cite{Kohonen2000Debye}, according to which a smaller permittivity corresponds to a thinner electrostatic screening layer. 

Table~\ref{tal:exacteps} shows the convergence performance and computational efficiency of the exact PrePD realizations with BGS or Schur-PCG for different values of $\epsilon$, where the dual subproblem is solved until convergence. We observe that the number of outer primal--dual iterations (PD Iter) increases as $\epsilon$ decreases, indicating that the small-$\epsilon$ regime is more difficult to solve. While the outer PD iter is the same for BGS and Schur-PCG, the inner dual solver iterations required for convergence differ. In particular, BGS is faster for moderate $\epsilon$ but the BGS sweeps increase significantly for small $\epsilon$, whereas Schur-PCG becomes more robust and efficient for small $\epsilon$. 

In practice, the dual subproblem need not be solved exactly at each outer PD iteration. Table~\ref{tal:Inexacteps} shows the convergence performance and computational efficiency of the inexact PrePD realizations for different values of $\epsilon$, where the dual subproblem is solved with only one inner iteration at each outer iteration. Intriguingly, we observe that the outer PD iter is the same as in the exact PrePD realizations even for small $\epsilon$, indicating that the inexact dual solver does not affect the convergence of the outer primal--dual iterations. 

\begin{table}
	\centering
	\tabcaption{Convergence performance of PrePD using exact dual solver BGS and Schur-PCG for different $\epsilon$ over 20 JKO steps.}
	\begin{tabular}{cccccc}
		\toprule
		\multirow{2}{*}{$\epsilon$} & \multirow{2}{*}{$\text{Total PD Iter}$}
		& \multicolumn{2}{c}{BGS} & \multicolumn{2}{c}{Schur-PCG}\\
		\cmidrule(lr){3-4}\cmidrule(lr){5-6}
		& & $\text{Total CPU}$ & $\text{BGS Sweeps/PD Iter}$
		& $\text{Total CPU}$ & $\text{PCG Iter/PD Iter}$\\
		\midrule
		$1$    & 36943 & \textbf{3784.31} & 2  & 4336.75 & 3\\
		$0.5$  & 37602 & \textbf{4757.20} & 3  & 5322.08 & 3\\
		$0.1$  & 44493 & 8592.41 & 6  & \textbf{5890.45} & 4\\
		$0.06$ & 49048 & 13510.17 & 10 & \textbf{9222.53} & 4\\
		\bottomrule
		\label{tal:exacteps}
	\end{tabular}
\end{table}

\begin{table}
	\centering
	\tabcaption{Convergence performance of PrePD using inexact dual solver iBGS and iSchur-PCG for different $\epsilon$ over 20 JKO steps.}
	\begin{tabular}{cccccc}
		\toprule
		\multirow{2}{*}{$\epsilon$} & \multirow{2}{*}{$\text{Total PD Iter}$}
		& \multicolumn{2}{c}{iBGS} & \multicolumn{2}{c}{iSchur-PCG}\\
		\cmidrule(lr){3-4}\cmidrule(lr){5-6}
		& & $\text{Total CPU}$ & $\text{BGS Sweeps/PD Iter}$
		& $\text{Total CPU}$ & $\text{PCG Iter/PD Iter}$\\
		\midrule
		$1$    & 36943 & \textbf{2192.75} & 1 & 3193.28 & 1\\
		$0.5$  & 37602 & \textbf{2478.27} & 1 & 3259.69 & 1\\
		$0.1$  & 44493 & \textbf{2716.54} & 1 & 4189.59 & 1\\
		$0.06$ & 49048 & \textbf{3581.89} & 1 & 4315.78 & 1\\
		\bottomrule
		\label{tal:Inexacteps}
	\end{tabular}
\end{table}

\begin{figure}[htbp]
	\centering
	\subfigure{
		\begin{minipage}{0.31\textwidth}
			\includegraphics[width=\linewidth]{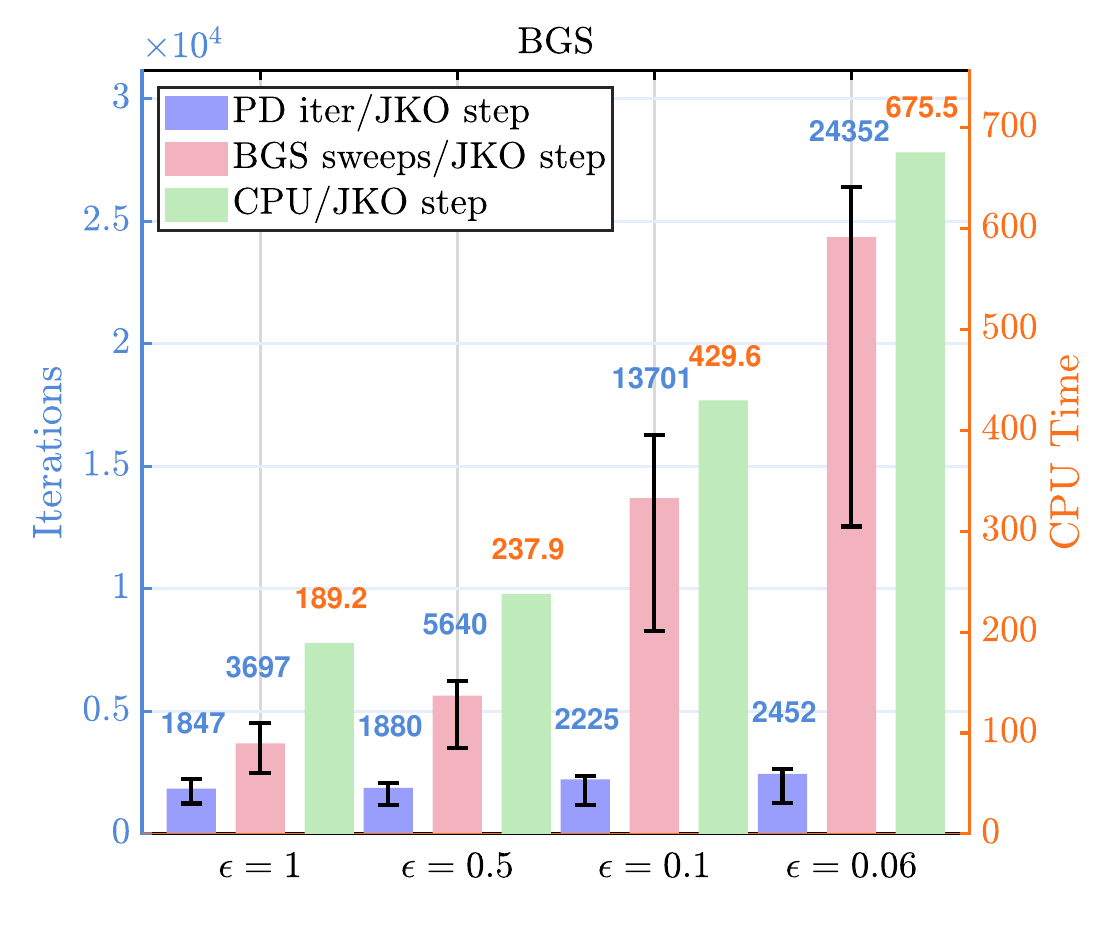}
			\includegraphics[width=\linewidth]{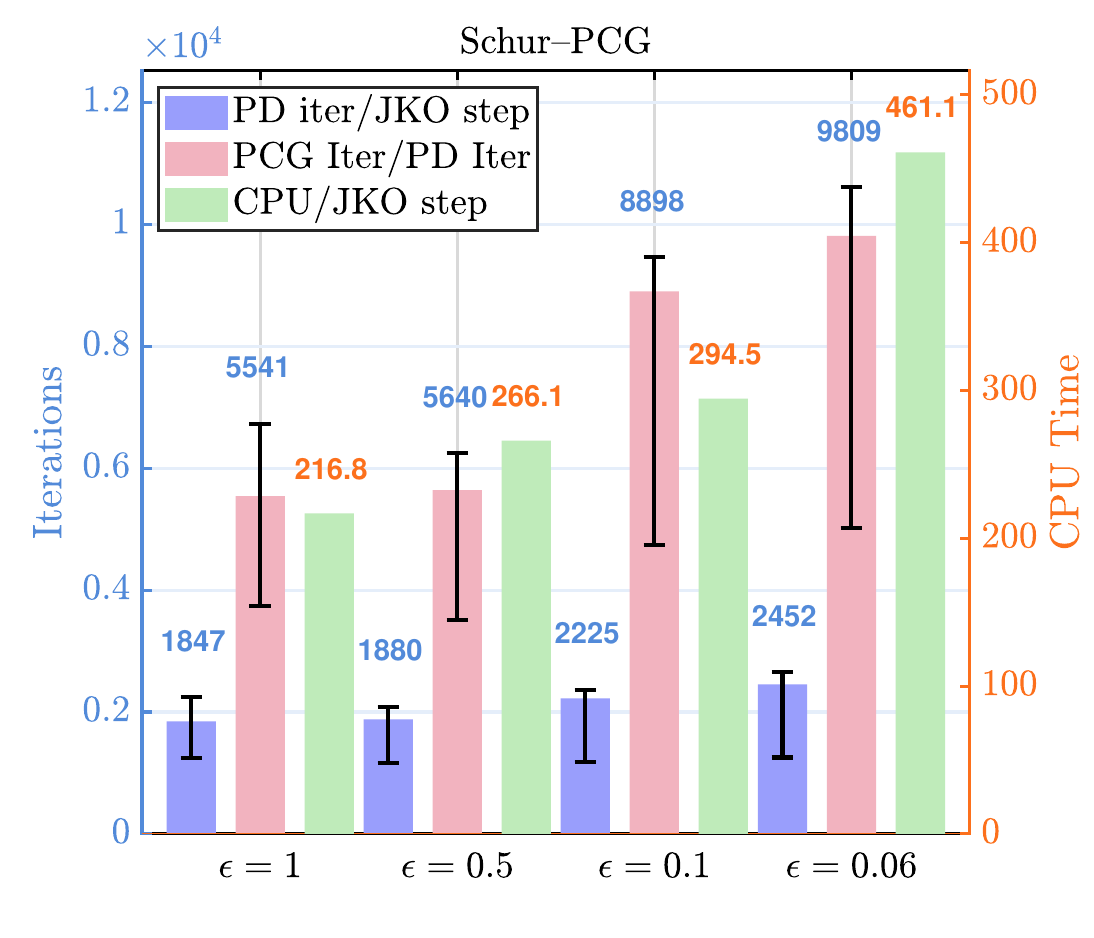}
		\end{minipage}
	}
	\subfigure{
		\begin{minipage}{0.31\textwidth}
			\includegraphics[width=\linewidth]{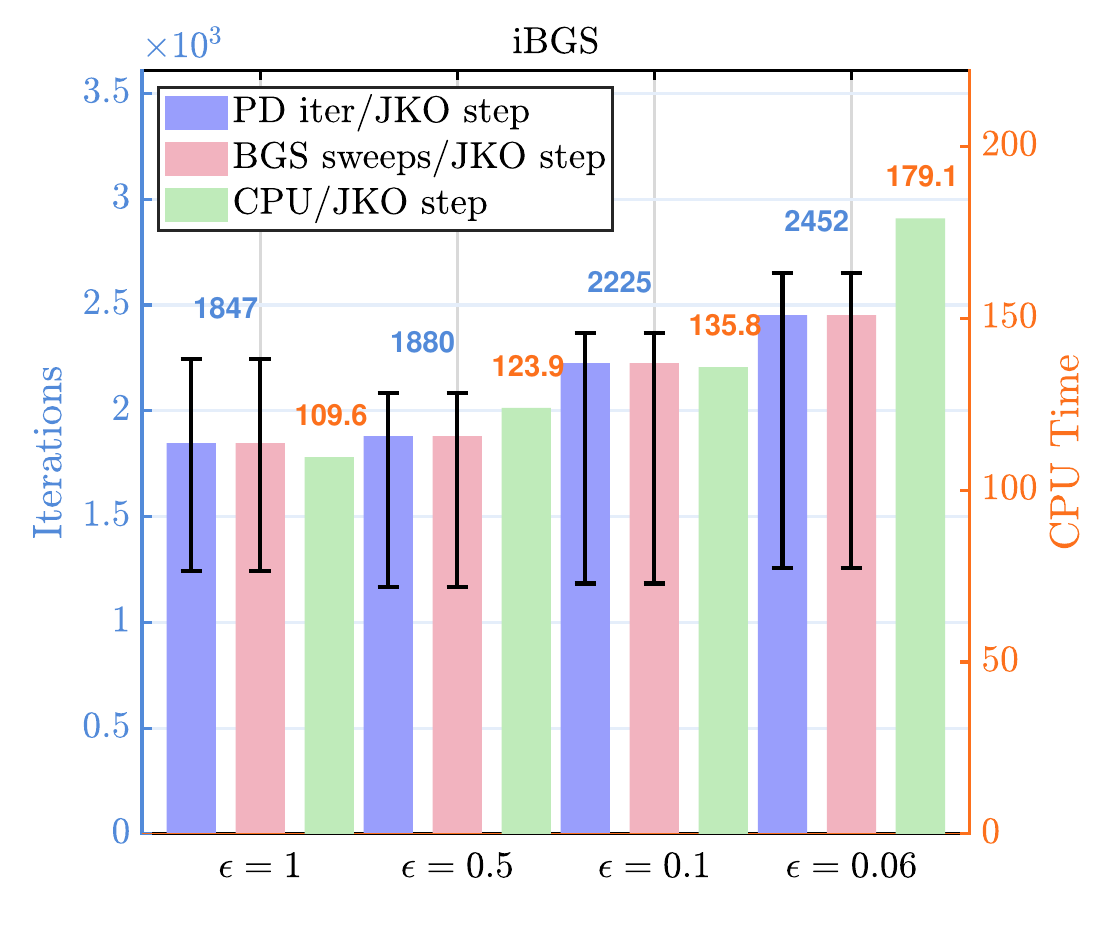}
			\includegraphics[width=\linewidth]{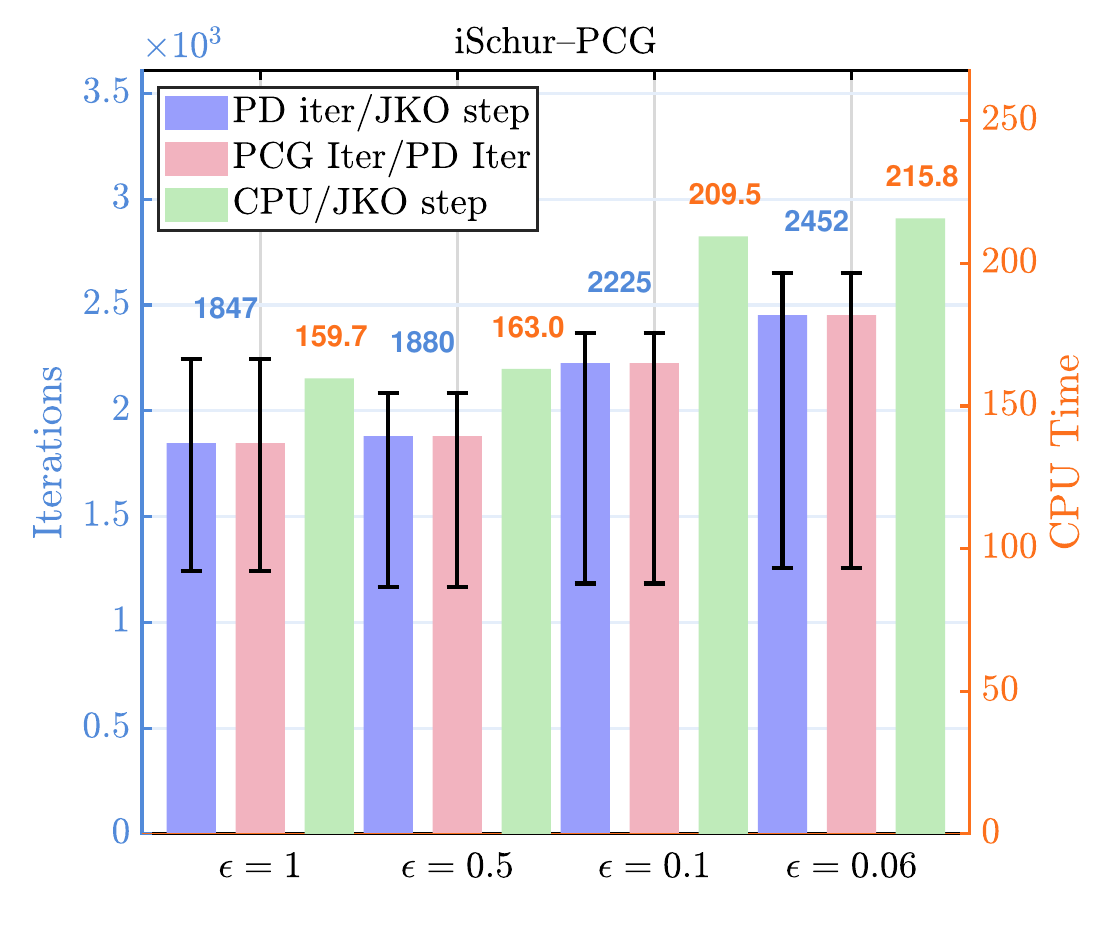}
		\end{minipage}
	}
	\subfigure{
		\begin{minipage}{0.31\textwidth}
			\includegraphics[width=\linewidth]{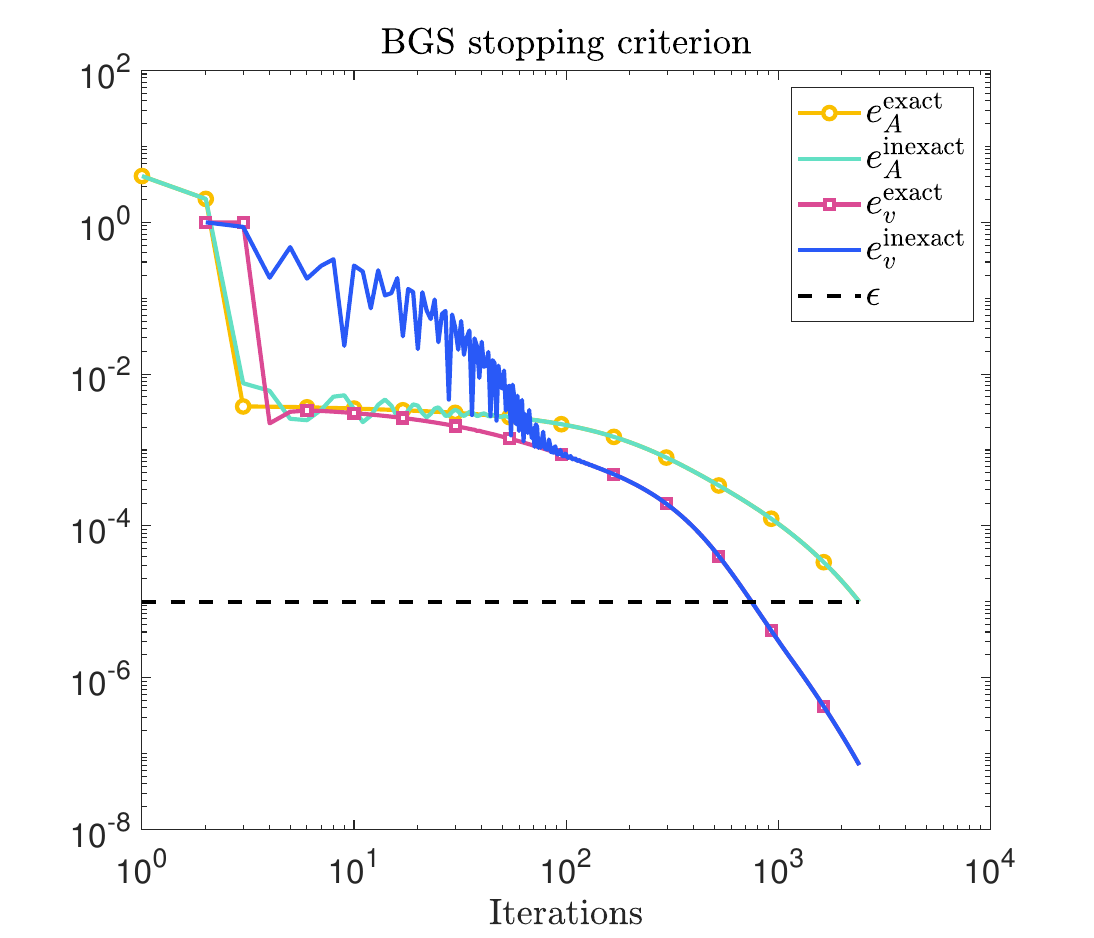}
			\includegraphics[width=\linewidth]{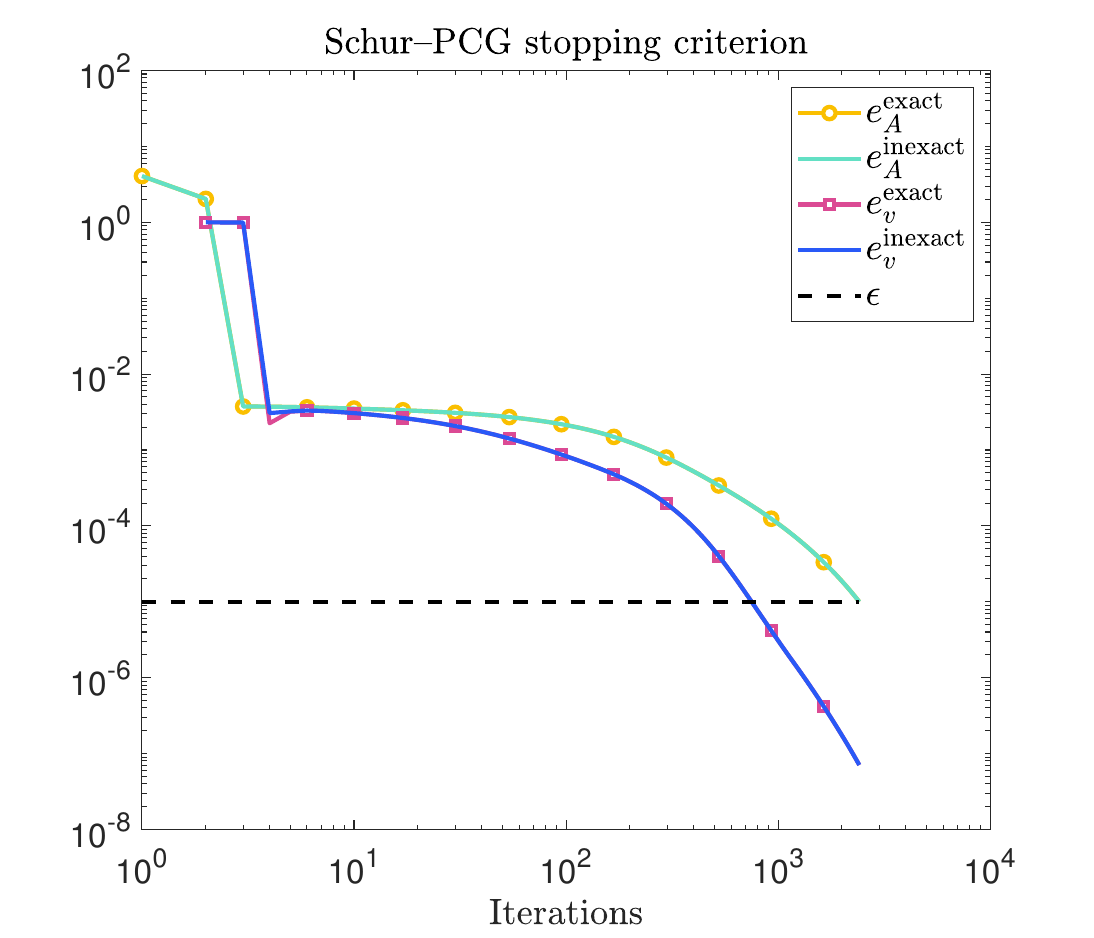}
		\end{minipage}
	}
	\caption{Mean PD iterations, mean inner dual solver iterations, and CPU time per JKO step for exact PrePD (left), inexact PrePD (middle), and the comparison of convergence between exact and inexact PrePD (right) for BGS (top) and Schur-PCG (bottom).}
	\label{fig:BGSSchurPCG}
\end{figure}

We summarize the mean primal--dual iterations, mean inner dual solver iterations, and CPU time for one JKO step in Fig.~\ref{fig:BGSSchurPCG}. The inexact dual solver saves the extra inner iterations (Inner Iter minus PD Iter per JKO step) and significantly reduces the total CPU time for all tested values of $\epsilon$. To understand why the inexact dual solver does not affect the convergence of the outer primal--dual iterations, we plot the convergence behavior of the dual variable ($e_v$) and the linear constraint residual ($e_A$) for exact and inexact dual solvers during one JKO step in Fig.~\ref{fig:BGSSchurPCG}. For both BGS and Schur-PCG, $e_A$ and $e_v$ of the inexact dual solver differ from those of the exact dual solver for the first few PD iterations, whereas they eventually converge to the same trajectory to convergence with the same PD iteration count for one JKO step. 

\subsubsection{Comparison between PrePD and VPTPD methods}\label{sec:6.1.4}
\begin{figure}[htbp]
	\centering
	\subfigure{
		\begin{minipage}{0.23\textwidth}
			\includegraphics[width=\linewidth]{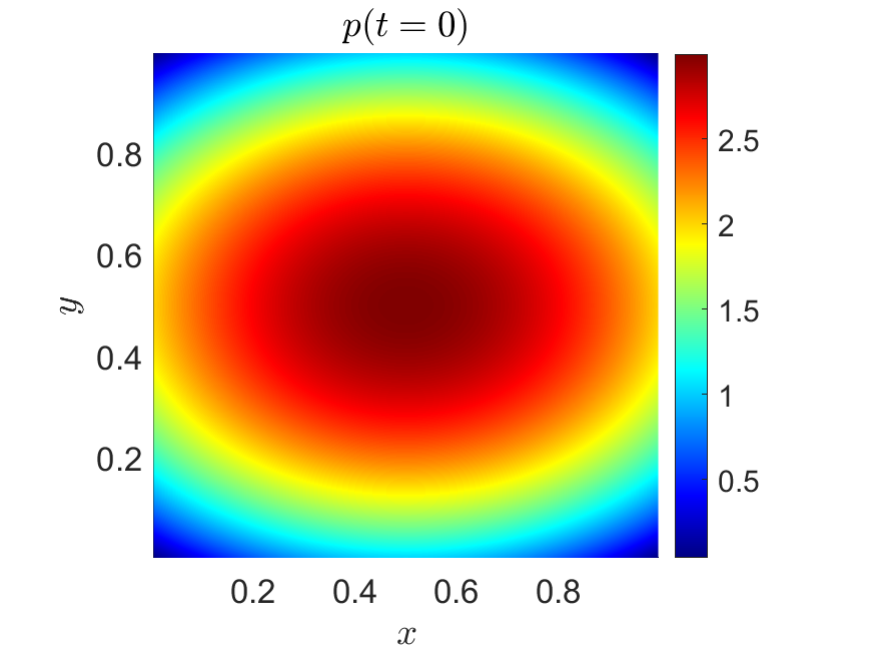}
			\includegraphics[width=\linewidth]{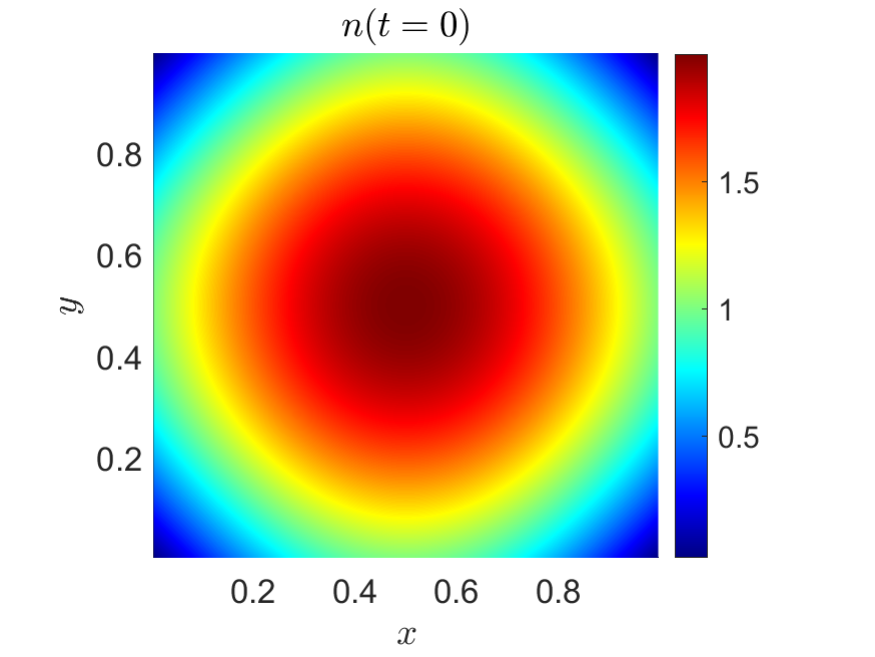}
		\end{minipage}
	}
	\subfigure{
		\begin{minipage}{0.23\textwidth}
			\includegraphics[width=\linewidth]{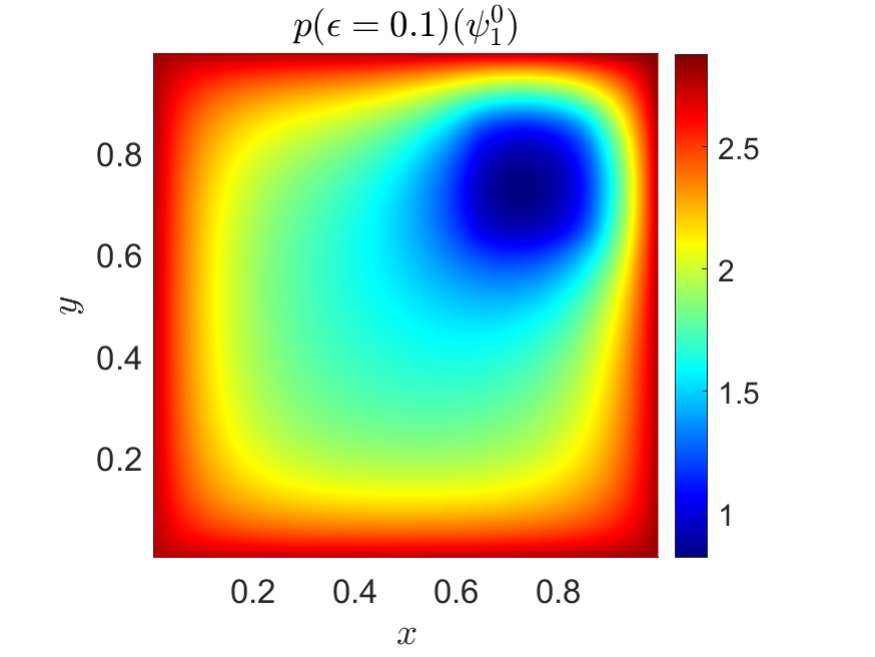}
			\includegraphics[width=\linewidth]{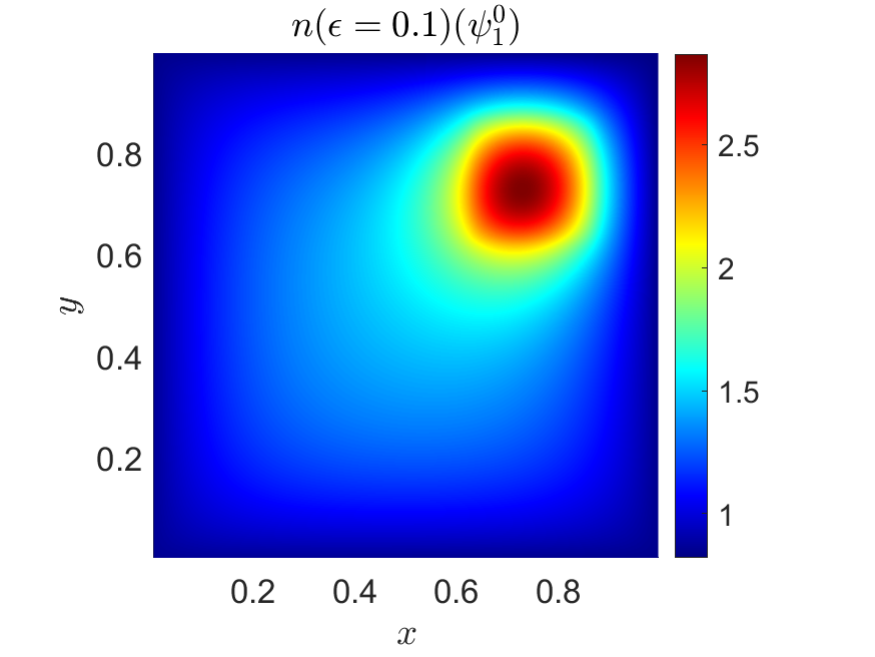}
		\end{minipage}
	}
	\subfigure{
		\begin{minipage}{0.23\textwidth}
			\includegraphics[width=\linewidth]{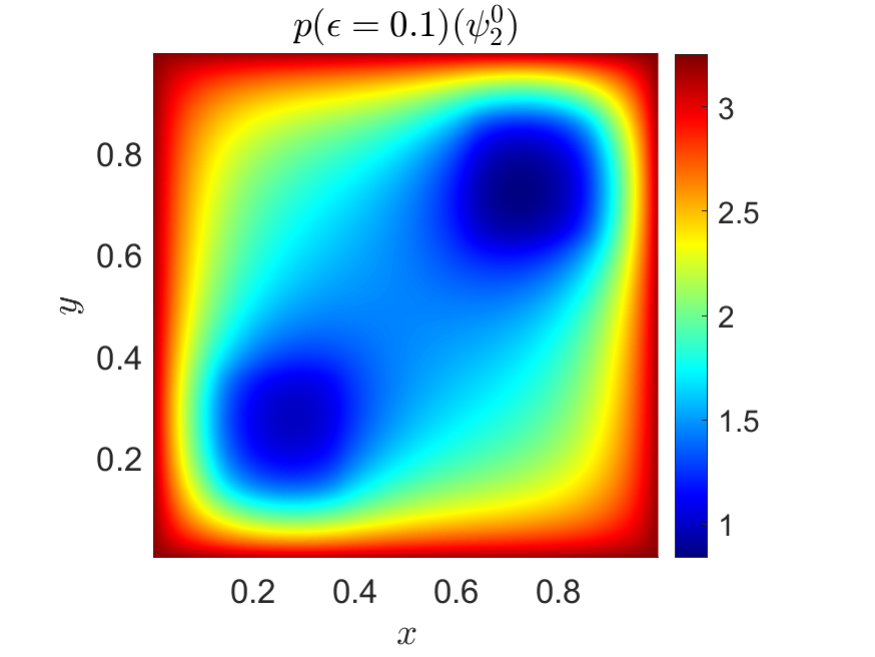}
			\includegraphics[width=\linewidth]{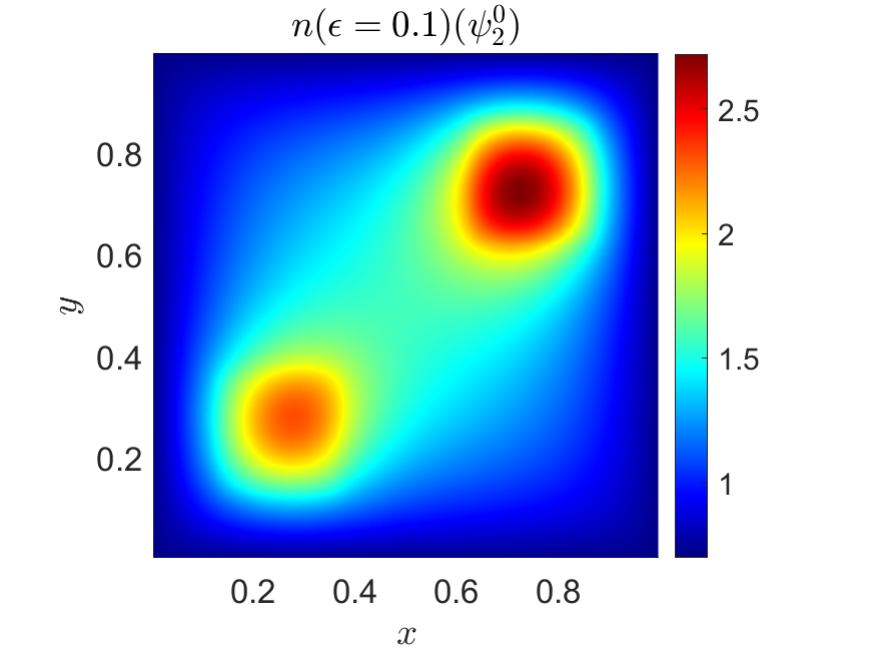}
		\end{minipage}
	}
	\subfigure{
		\begin{minipage}{0.23\textwidth}
			\includegraphics[width=\linewidth]{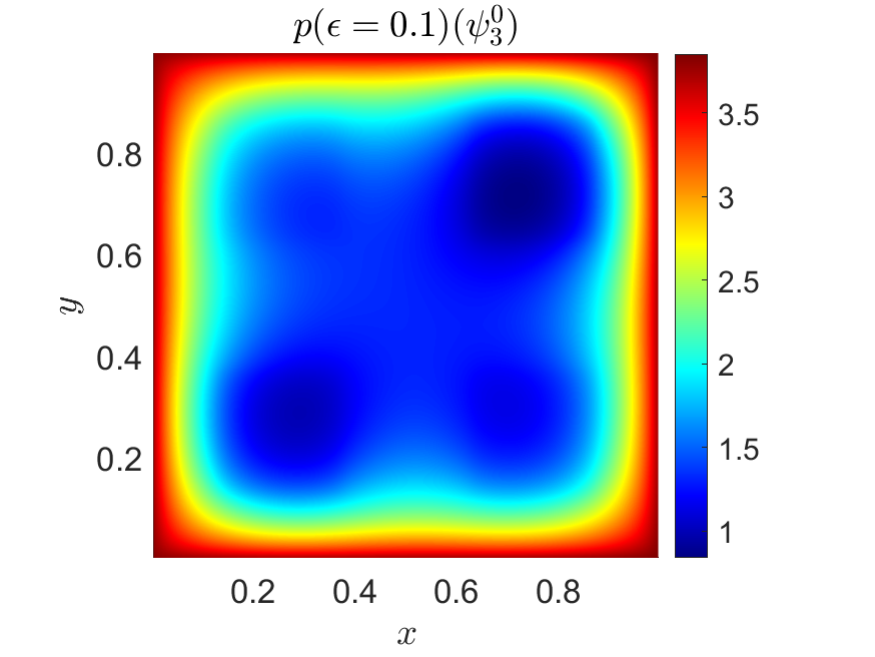}
			\includegraphics[width=\linewidth]{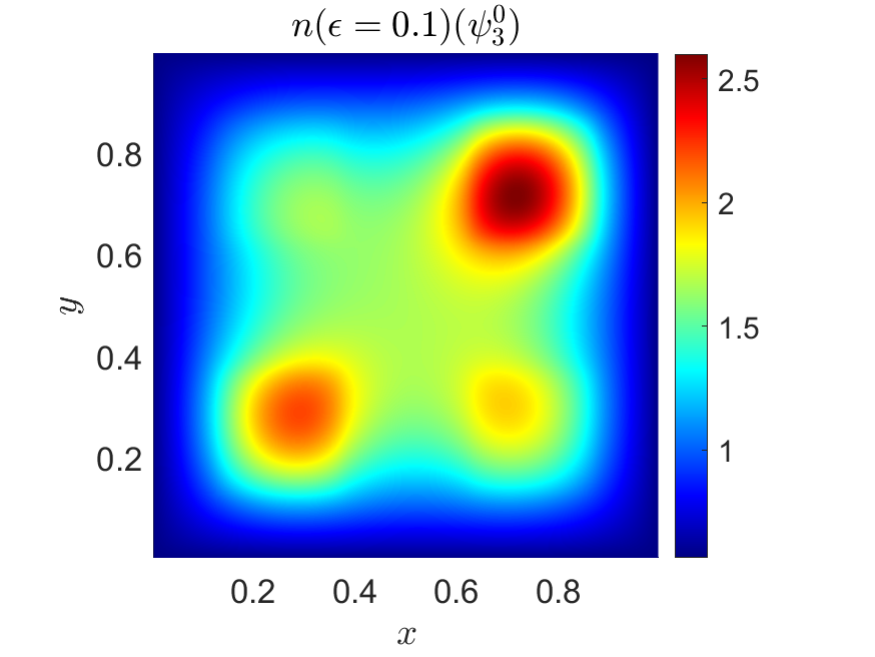}
		\end{minipage}
	}
	\caption{Initial ionic distributions and profiles at $t=0.2$ for different fixed charge densities $\psi^{0}_{1}$ \eqref{eq:psi01}, $\psi^{0}_{2}$, and $\psi^{0}_{3}$ \eqref{eq:diffpsi0}.}
	\label{fig:diffpsipandn}
\end{figure}

We next compare the performance of PrePD and VPTPD methods developed in this work. We consider the 2D PNP model with the initial and Dirichlet boundary conditions in \eqref{eq:p0n0bcphi} and the four-region fixed-charge profile $\psi^{0}_{3}$ defined by:
\begin{eqnarray}\label{eq:diffpsi0}
	\begin{aligned}
		&\psi^{0}_{2} = \psi^{0}_{1}+
		\begin{cases}
			8, & \text{if } \dfrac{1}{8}\leq x\leq\dfrac{3}{8}\text{ and }\dfrac{1}{8} \leq y\leq\dfrac{3}{8}, \\
			0, & \text{otherwise}.
		\end{cases}\\
		&\psi^{0}_{3} = \psi^{0}_{2}+
		\begin{cases}
			6, & \text{if } \dfrac{5}{8}\leq x\leq\dfrac{7}{8}\text{ and }\dfrac{1}{8} \leq y\leq\dfrac{3}{8}, \\
			0, & \text{otherwise}.
		\end{cases}+
		\begin{cases}
			4, & \text{if } \dfrac{1}{8}\leq x\leq\dfrac{3}{8}\text{ and }\dfrac{5}{8} \leq y\leq\dfrac{7}{8}, \\
			0, & \text{otherwise}.
		\end{cases}
	\end{aligned}
\end{eqnarray}
The profiles for ionic concentrations corresponding to the fixed charge density $\psi^{0}_{1}$, $\psi^{0}_{2}$, and $\psi^{0}_{3}$ are shown in Fig.~\ref{fig:diffpsipandn}. We observe that the ionic distributions are strongly influenced by the spatial distribution of the fixed charges, with negative ions accumulating near positive fixed charges and positive ions being depleted in those regions. The degree of accumulation or depletion also depends on the charge magnitude: regions carrying larger positive fixed charges attract more negative ions and induce a stronger reduction in the positive-ion concentration.

\begin{figure}[htbp]
	\centering
	\includegraphics[width=0.4\linewidth]{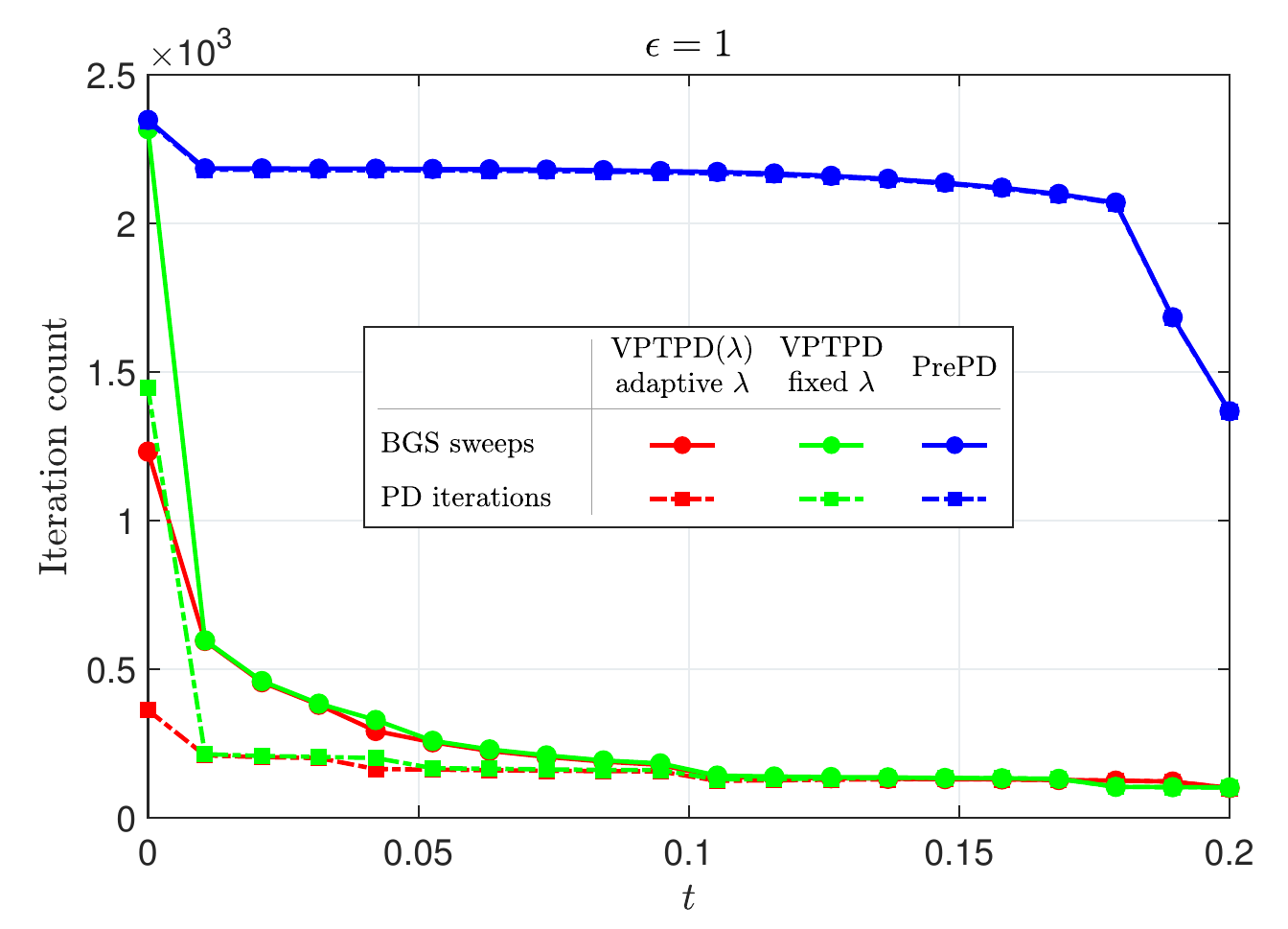}
	\includegraphics[width=0.4\linewidth]{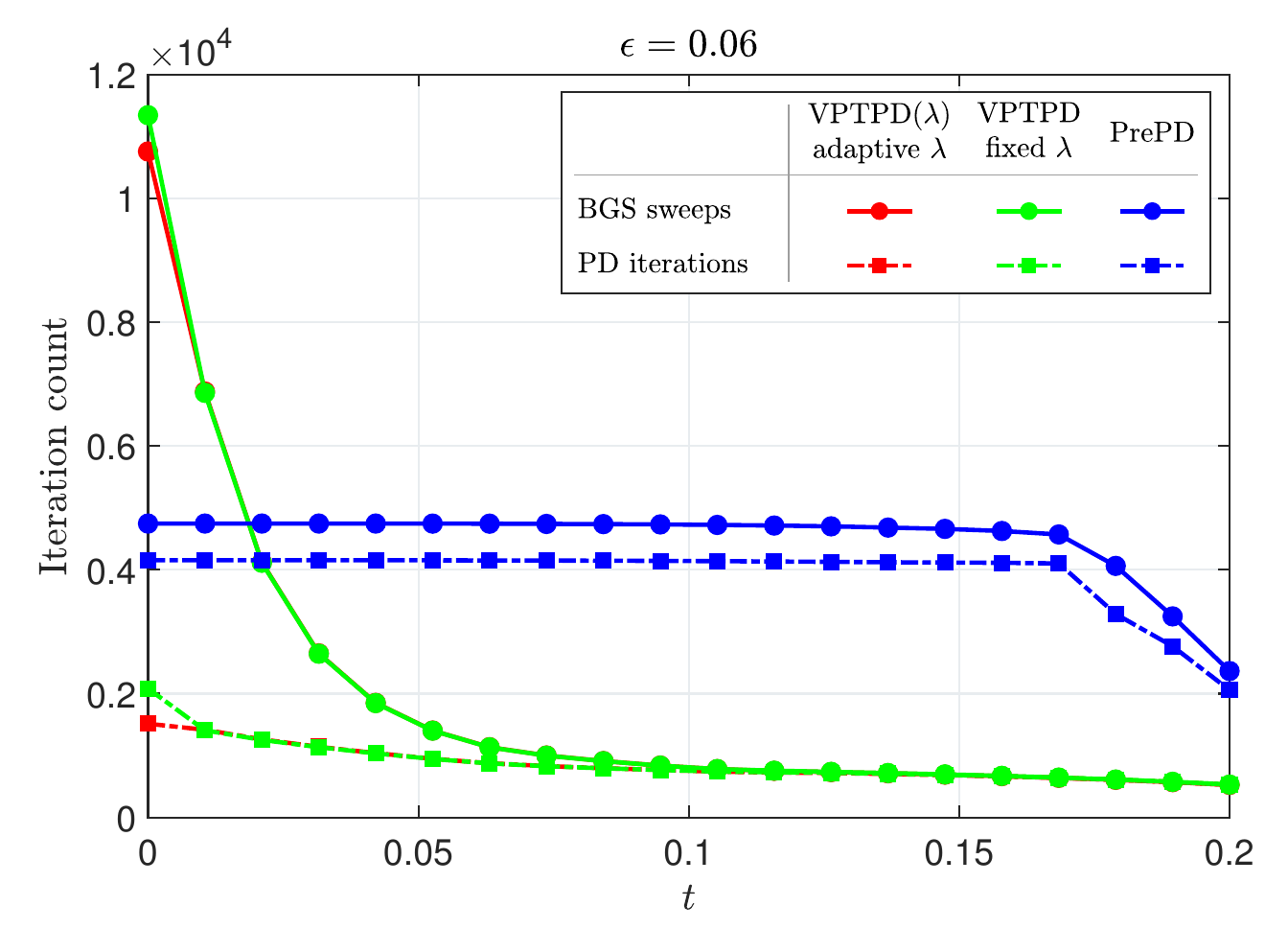}
    \caption{Total outer primal--dual iteration and inner BGS sweeps for PrePD, VPTPD, and VPTPD($\lambda$) for different values of $\epsilon$.}
	\label{fig:PrePDVPTPD}
\end{figure}

We compare the PD iteration and computational efficiency of the PrePD, VPTPD, and VPTPD with adaptive stepsize strategy (VPTPD($\lambda$)) to solve the above four-region fixed charged problem. Table~\ref{tal:PrePDVPTPD} shows that the VPTPD and VPTPD($\lambda$) methods substantially reduce the primal--dual iteration and total CPU time relative to PrePD for both $\epsilon=1$ and $\epsilon=0.06$. VPTPD($\lambda$) achieves better computational efficiency in both tests, although the acceleration effect due to the adaptive stepsize strategy becomes less pronounced for small $\epsilon$. Moreover, we summarize the primal--dual iteration and inner BGS sweeps for PrePD, VPTPD, and VPTPD($\lambda$) in Fig.~\ref{fig:PrePDVPTPD}. We observe that VPTPD and VPTPD($\lambda$) require more extra inner BGS sweeps (i.e., the difference between total BGS sweeps and total PD iterations) than PrePD, especially for small $\epsilon$, mainly because the VPTPD methods involve the inversion of a more complex matrix in the dual subproblem. However, the overall computational efficiency of VPTPD is still better than PrePD, indicating that the reduction in outer PD iterations outweighs the increase in inner BGS sweeps.	

\begin{table}[htbp]
	\centering
	\tabcaption{Comparison of primal--dual iteration and computational efficiency of PrePD, VPTPD, and VPTPD($\lambda$) for $\epsilon=1$ and $\epsilon=0.06$ over 20 JKO steps.}
	\begin{tabular}{clccccc}
		\toprule
		$\epsilon$ & $\text{Method}$ & $\text{Total PD Iter}$ & $\text{CPU}$ & $\text{CPU/JKO}$ & $\text{CPU/PD Iter}$ & $\text{Total BGS sweeps}$ \\
		\midrule
		\multirow{3}{*}{$1$}
		& PrePD & 42014 & 3201.55 & 160.08 & 0.0762 & 42096\\
		& VPTPD & 4303 & 336.44 & 16.82 & 0.0782 & 6403\\
		& VPTPD($\lambda$) & 3154 & \textbf{246.01} & \textbf{12.30} & 0.0780 & \textbf{5255}\\
		\midrule
		\multirow{3}{*}{$0.06$}
		& PrePD & 78483 & 6520.52 & 326.03 & 0.0831 & 89797\\
		& VPTPD & 17737 & 1748.80 & 87.44 & 0.0986 & 38812\\
		& VPTPD($\lambda$) & 17166 & \textbf{1515.84} & \textbf{75.79} & 0.0883 & \textbf{38232}\\
		\bottomrule
		\label{tal:PrePDVPTPD}
	\end{tabular}
\end{table}

\subsection{Extended experiments for modified PNP models}\label{sec:6.2}
\begin{figure*}[t]
	\centering
	\includegraphics[width=0.7\textwidth]{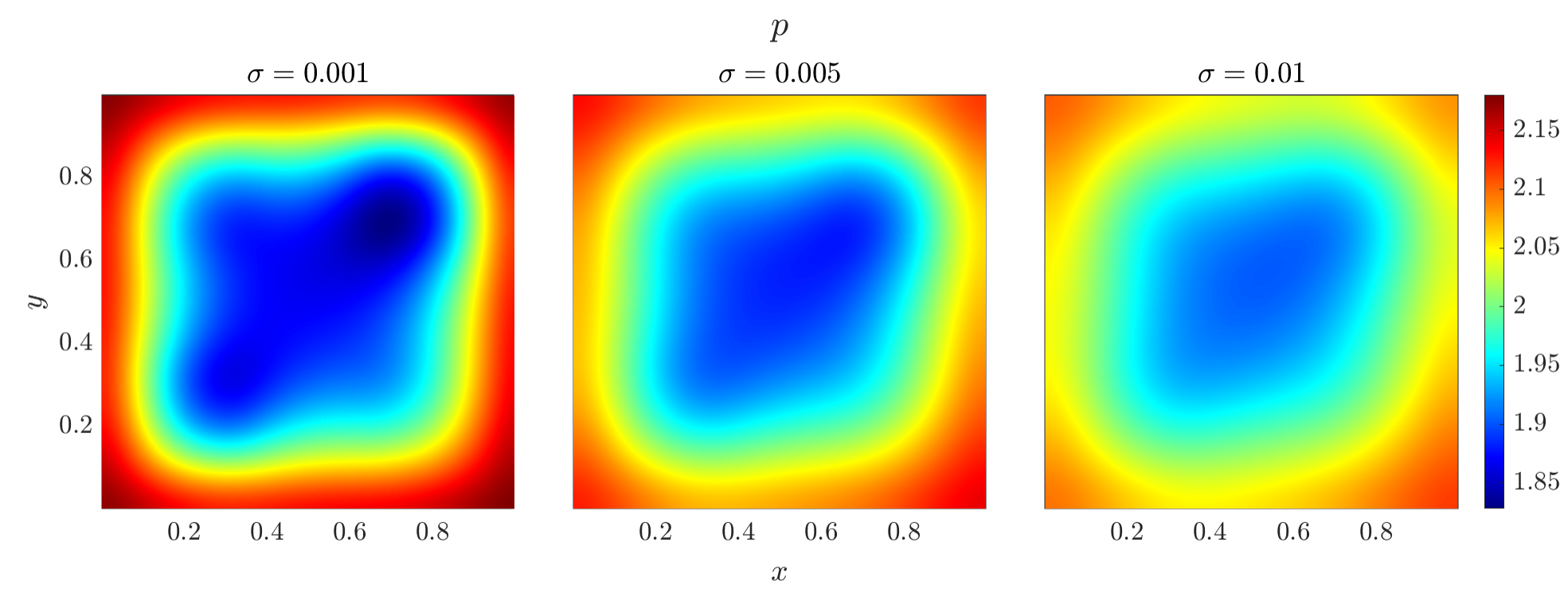}
	\includegraphics[width=0.7\textwidth]{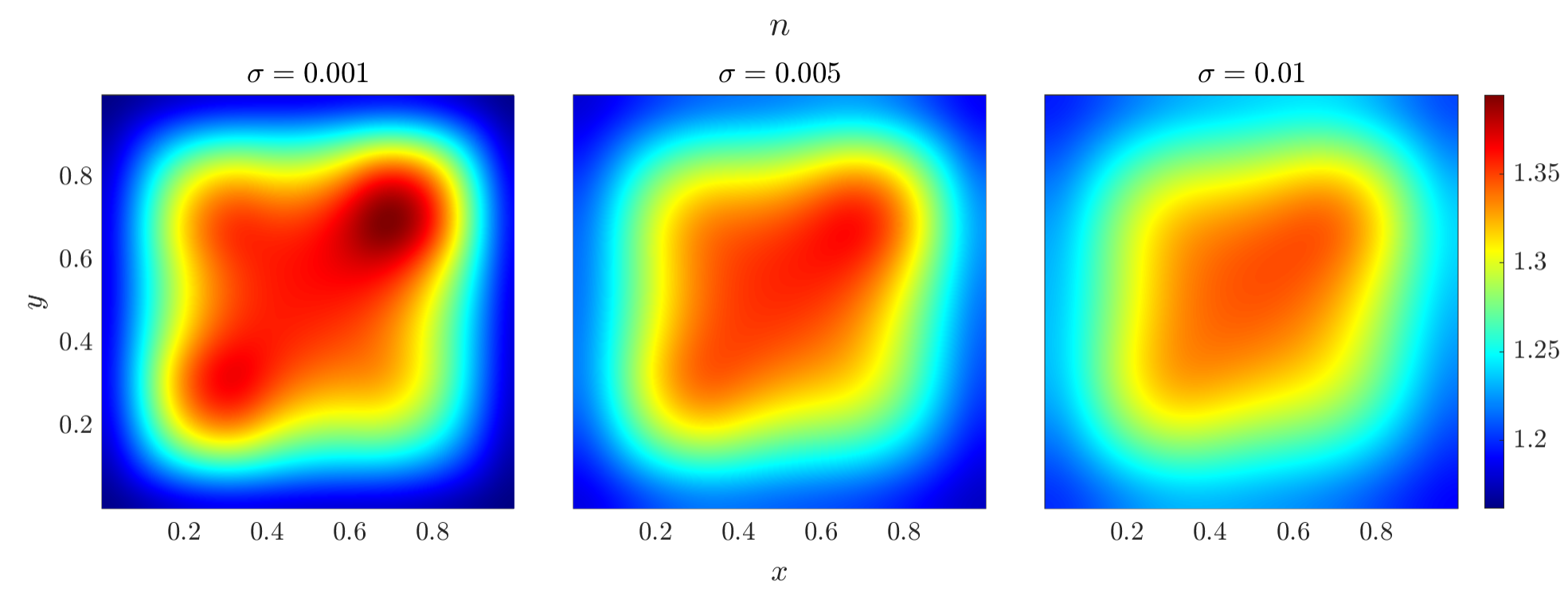}
	\includegraphics[width=0.7\textwidth]{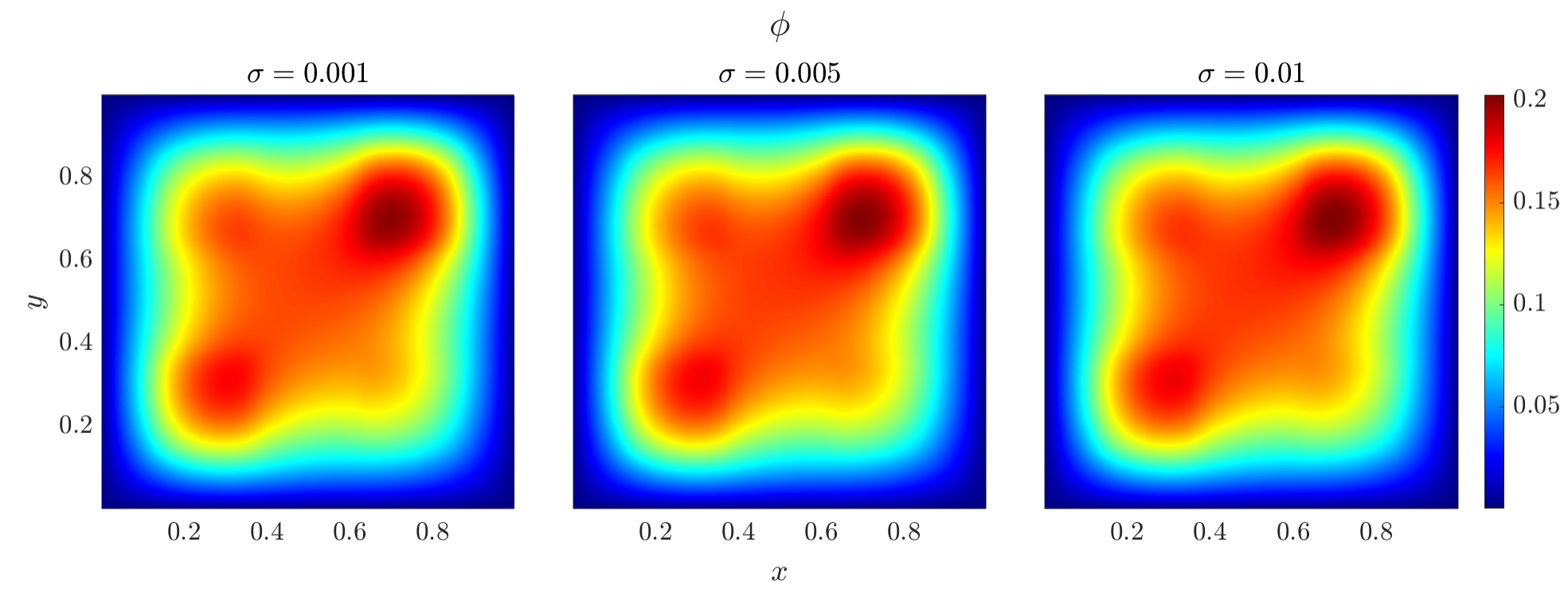}
	\caption{Equilibrium profiles of $p$, $n$, and $\phi$ at $t=1$ for 2D modified PNP models with different concentration-gradient coefficients $\sigma$.}
	\label{fig:diffG}
\end{figure*}

We now investigate how the concentration-gradient energy modifies the ionic redistributions in response to fixed charges. We first consider a 2D modified PNP model with the initial and boundary conditions in \eqref{eq:p0n0bcphi} together with the four-region fixed charge $\psi^{0}_{3}$ in \eqref{eq:diffpsi0}. The equilibrium distributions of holes $p$ and electrons $n$ and the electrostatic potential $\phi$ for different values of the concentration-gradient coefficient $\sigma$ are shown in Fig.~\ref{fig:diffG}, where we set $\sigma_{p}=\sigma_{n}=\sigma\in\left\{0.001,0.005,0.01,0.03\right\}$ and $G=\boldsymbol{0}$. The results show that the concentration-gradient energy produces smoother, more spatially coherent ionic profiles. As $\sigma$ increases, local variations are suppressed and the ionic distributions become less sensitive to individual fixed-charge interfaces. The electrostatic potential changes more moderately, reflecting its nonlocal coupling to the smoothed charge density.

Finally, we consider a three-dimensional modified PNP model on $\Omega=[-0.2,0.2]^{3}$, where the concentration-gradient coefficient and the symmetric interaction matrix are given by:
\begin{align}
	\sigma=0.005,\quad G=
	\begin{pmatrix}
		10 & 1\\
		1 & 10
	\end{pmatrix}.
\end{align}
The initial conditions are given by: 
\begin{eqnarray}\label{eq:3Dinitial}
	\begin{aligned}
		&(p(x,y,z,0),n(x,y,z,0)) =
		\begin{cases}
			(1.5, 2.0),& \quad\sqrt{(x-0.1)^{2}+(y-0.1)^{2}+(z-0.1)^{2}}<0.075, \\
			(2.0, 1.5),& \quad\sqrt{(x+0.1)^{2}+(y+0.1)^{2}+(z+0.1)^{2}}<0.075, \\
			(\varepsilon,\varepsilon), & \quad\text{otherwise},
		\end{cases}
	\end{aligned}
\end{eqnarray}
where $\varepsilon=10^{-6}$ represents the background concentration with a small positive value to avoid the singularity of the logarithmic entropy term. The electrostatic potential is subject to mixed Dirichlet--Neumann boundary conditions:
\begin{align}\label{eq:3Dboundary}
	\phi(a,y,z)=1,\quad
	\dfrac{\partial\phi}{\partial\mathbf{n}}\Big|_{(b,y,z)}=0,\quad\phi(x,c,z)=1,\quad
	\dfrac{\partial\phi}{\partial\mathbf{n}}\Big|_{(x,d,z)}=0,\quad\phi(x,y,e)=1,\quad
	\dfrac{\partial\phi}{\partial\mathbf{n}}\Big|_{(x,y,f)}=0.
\end{align}
The fixed charge is given by a localized distribution:
\begin{align}\label{eq:3Dpsi01}
	\psi^{0}_{3\mathrm{D}} =
	\begin{cases}
		10, & \text{if } \dfrac{1}{20}\leq x\leq\dfrac{3}{20}\text{ and }\dfrac{1}{20}\leq y\leq\dfrac{3}{20}\text{ and }\dfrac{1}{20}\leq z\leq\dfrac{3}{20}, \\
		0, & \text{otherwise}.
	\end{cases}
\end{align}

\begin{figure*}[t]
	\centering
	\includegraphics[width=0.96\textwidth]{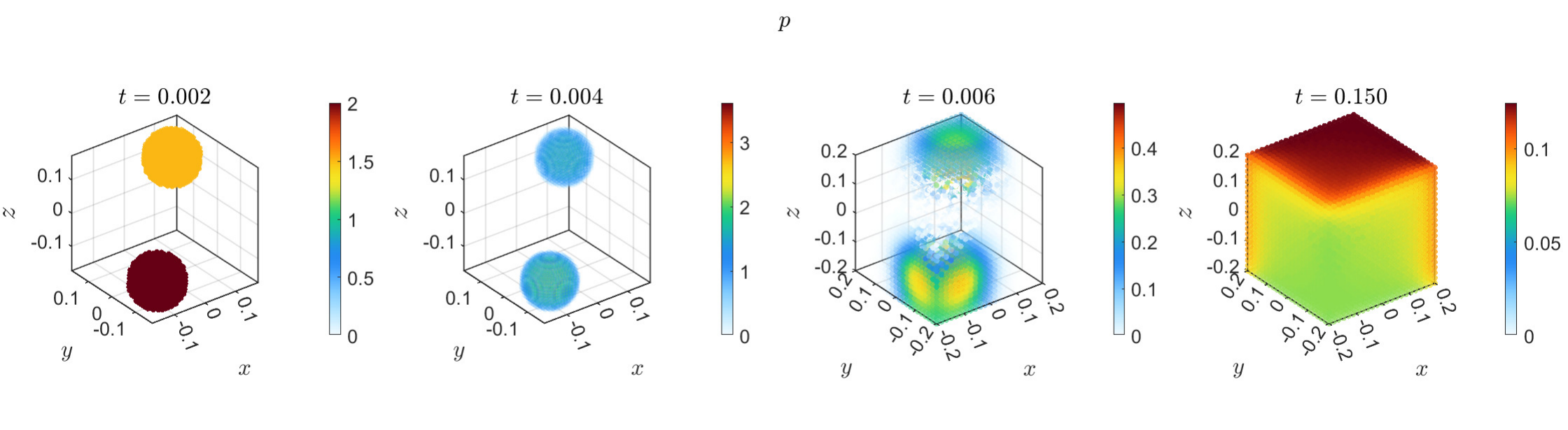}
	\vspace{-2mm}
	\includegraphics[width=0.96\textwidth]{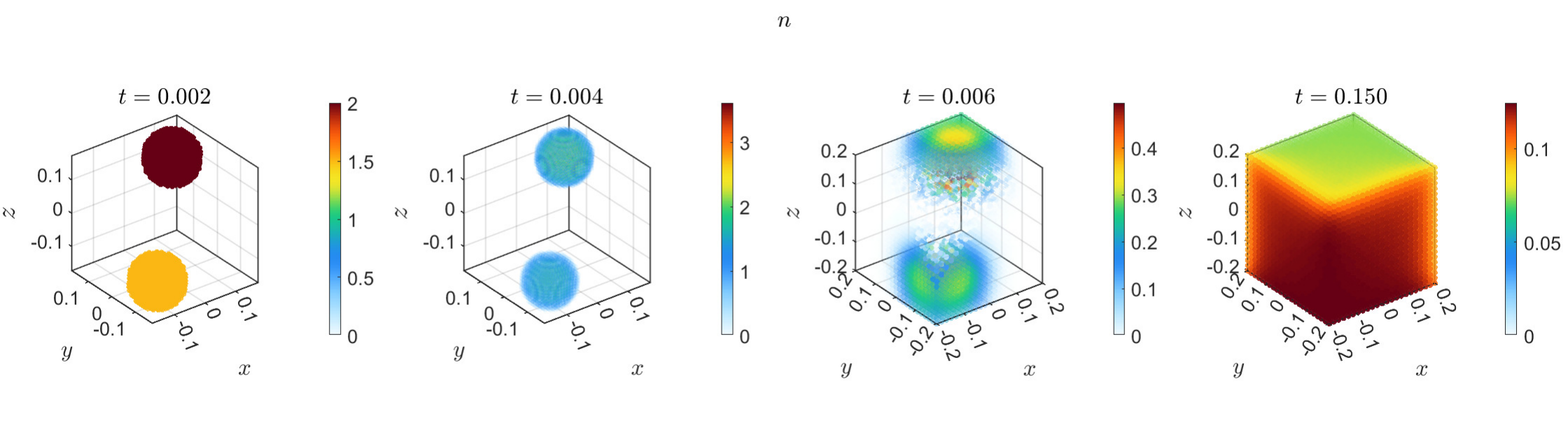}
	\vspace{-2mm}
	\includegraphics[width=0.96\textwidth]{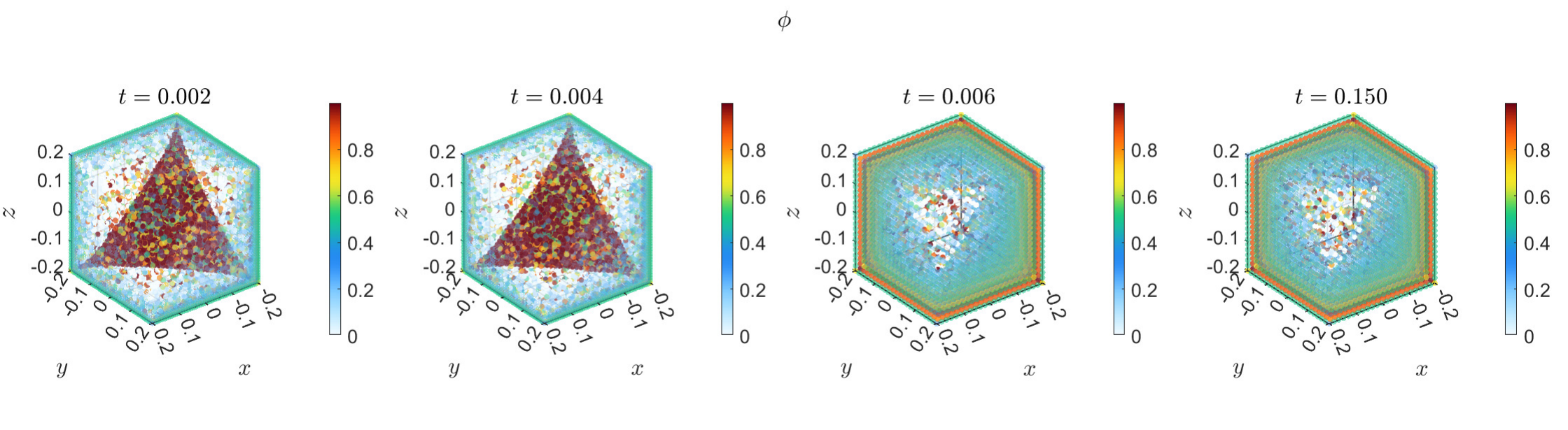}
 	\caption{Evolution of the three-dimensional numerical solution for the modified PNP system with the fixed charge $\psi^{0}_{3D}$.}
    \label{fig:3Dspheres}
\end{figure*}


We solve the above 3D modified PNP model with $N_x=N_y=N_z=45$, $\tau=0.002$ for $T=0.15$. The evolution of the ion distributions $p$ and $n$ is shown in Fig.~\ref{fig:3Dspheres}. At early times, the two ionic concentrations are confined to two separated spherical regions, with their high-concentration zones located at opposite positions. As the system evolves, diffusion and concentration-gradient regularization smooth the sharp interfaces, while the electrostatic field generated by the localized fixed charge produces an anisotropic redistribution of the ions. Consequently, the spherical structures gradually expand, deform, and spread throughout the computational domain. Since $p$ and $n$ carry opposite charges, their electrostatic drift directions are opposite, leading to increasingly complementary spatial profiles. By $t=0.015$, the initially localized structures have developed into boundary-dominated distributions, whereas $\phi$ approaches a smooth quasi-steady configuration governed by the fixed charge and the mixed Dirichlet-Neumann boundary conditions.

\section{Conclusion}
In this paper, we proposed a unified Poisson-constrained JKO scheme for strongly coupled multi-variable PNP models and two efficient primal--dual splitting algorithms (PrePD and VPTPD) equipped with fast dual solvers to solve the resulting constrained optimization problems. This scheme is an extension of the JKO framework for Wasserstein gradient flows \cite{Carrillo2022PrimalDual,Carrillo2024StructurePD,Zeng2026variable} to PNP models, where the Poisson equation is incorporated as an additional linear constraint, thereby preserving the desired properties of energy dissipation, positivity of the ionic concentrations, and mass conservation. The proposed primal--dual splitting algorithms display superior robustness and efficiency compared to some existing methods in numerous experiments with various boundary conditions and fixed charges, especially for small dielectric permittivity. The current proposed scheme is first-order accurate in time, and we will investigate the higher-order variational structure \cite{Cances2026SecondOrder} in future work.

\section*{Declaration of competing interest}
The authors declare that they have no known competing financial interests or personal relationships that could have appeared to influence the work reported in this paper.

\section*{Acknowledgements}
CW is supported by the National Natural Science Foundation of China under grants 12371392, 12622128 and 12431015. 
The work of ZZ is partially supported by National Key R\&D Program of China (2023YFA1011403), the NSFC grant (92470112 and 12426312), and Shenzhen Science and Technology Program (QNXMA20250701095429040).

\bibliographystyle{elsarticle-num}
\bibliography{Manu_refer}
\end{document}